{
\edef\inewcount{\noexpand\csname newcount\endcsname}
\edef\inewdimen{\noexpand\csname newdimen\endcsname}
\edef\inewskip{\noexpand\csname newskip\endcsname}
\edef\inewmuskip{\noexpand\csname newmuskip\endcsname}
\edef\inewbox{\noexpand\csname newbox\endcsname}
\edef\inewhelp{\noexpand\csname newhelp\endcsname}
\edef\inewtoks{\noexpand\csname newtoks\endcsname}
\edef\inewread{\noexpand\csname newread\endcsname}
\edef\inewwrite{\noexpand\csname newwrite\endcsname}
\edef\inewfam{\noexpand\csname newfam\endcsname}
\edef\inewlanguage{\noexpand\csname newlanguage\endcsname}
\edef\inewinsert{\noexpand\csname newinsert\endcsname}
\edef\inewif{\noexpand\csname newif\endcsname}

\inewif\ifscroll 
\inewif\ifsuppressunusedbib 



\def\warning#1{\immediate\write17{! #1}\immediate\write17{l.\the\inputlineno}\immediate\write17{}} 
\tracinglostchars=2

\long\def\gobble#1{}
\long\def\gobblepar#1\par{}
\def\expand#1{\edef\expandmacro{#1}\expandmacro\let\expandmacro\undefined}
\def\setetok#1#2{\expand{\noexpand#1{#2}}}

\def\link#1#2{\lhighlight{#2}}
\def\llink#1{\printlink{llink #1}\link{\ohash#1}}\catcode`\#=11 \def\ohash{#}\catcode`\#=6
\def\anchor#1#2{\printlink{anchor #1 #2}#2}

\def\dumpbox#1#2#3{\shipout\vbox{\unvbox#3}} 

\def\metadata#1#2{}
\def\src{} 

\inewread\epsffilein
\inewif\ifepsfbbfound\inewif\ifepsffilecont
\inewdimen\epsfxsize\inewdimen\epsfysize
\inewdimen\pspoints\pspoints1bp
\def\epsfbox#1{\openin\epsffilein=#1 \ifeof\epsffilein\errmessage{Could not open file #1}\else
	{\def\do##1{\catcode`##1=12}\dospecials\catcode`\ =10\epsffileconttrue
		\epsfbbfoundfalse
		\loop\read\epsffilein to\epsffileline \ifeof\epsffilein\epsffilecontfalse\else\expandafter\epsfaux\epsffileline :. \\\fi\ifepsffilecont\repeat
		\ifepsfbbfound\else\errmessage{No HiResBoundingBox comment found in file #1}\fi}%
	\closein\epsffilein\fi
	\epsfysize\epsfury\pspoints \advance\epsfysize-\epsflly\pspoints
	\epsfxsize\epsfurx\pspoints \advance\epsfxsize-\epsfllx\pspoints
	\hbox{\vbox to\epsfysize{\vfil\hbox to\epsfxsize{\includegraphics{#1}\hfil}}}}
\catcode`\%=12
\let\percent=%
\catcode`\%=14
\def\epsfbblit{\percent\percent HiResBoundingBox}
\def\epsfaux#1:#2\\{\def\testit{#1}\ifx\testit\epsfbblit \epsfgrab #2 . . . \\\epsffilecontfalse\epsfbbfoundtrue\fi}
\def\empty{}
\def\epsfgrab #1 #2 #3 #4 #5\\{\gdef\epsfllx{#1}\ifx\epsfllx\empty\epsfgrab #2 #3 #4 #5 .\\\else\gdef\epsflly{#2}\gdef\epsfurx{#3}\gdef\epsfury{#4}\fi}

\ifx\pdfoutput\undefined\else\ifnum\pdfoutput>0
\pdfcompresslevel=0 
\pdfobjcompresslevel=0
\def\llink#1#2{\lhighlight{\pdfstartlink goto name {#1}#2\pdfendlink}}
\def\link#1#2{\lhighlight{\pdfstartlink user { /Subtype /Link /A << /Type /Action /S /URI /URI (#1) >> }#2\pdfendlink}}
\def\anchor#1#2{\pdfdest name {#1} xyz #2}

\def\dumpbox#1#2#3{\pdfpagewidth#1 \pdfpageheight#2 \shipout\box#3}

\def\metadata#1#2{\pdfinfo{/Title (#1) /Author (#2)}}
\fi\fi

\def\plainfmtname{plain}\ifx\fmtname\plainfmtname\else
\edef\plainoutput{\the\output}
\global\chardef\itfam=4

\outputpenalty=0
\tracingstats=0
\newlinechar=-1
\maxdeadcycles=25
\showboxbreadth=5
\showboxdepth=3
\errorcontextlines=5
\overfullrule=5pt
\hsize=6.5in
\vsize=8.9in
\maxdepth=4pt
\parindent=20pt
\abovedisplayskip=12pt plus 3pt minus 9pt
\belowdisplayskip=12pt plus 3pt minus 9pt
\belowdisplayshortskip=7pt plus 3pt minus 4pt


\font\tensy=cmsy10

\font\tenbf=cmbx10

\catcode"18=12
\catcode`@=11
\global\let\end\@@end
\global\let\input\@@input
\def\eqalign#1{\null\,\vcenter{\openup\jot\m@th
  \ialign{\strut\hfil$\displaystyle{##}$&$\displaystyle{{}##}$\hfil
      \crcr#1\crcr}}\,}
\catcode`@=12
\fi

\font\tenmsa=msam10 \font\sevenmsa=msam7 \font\fivemsa=msam5 \newfam\msafam \textfont\msafam=\tenmsa \scriptfont\msafam=\sevenmsa \scriptscriptfont\msafam=\fivemsa
\font\teneufm=eufm10 \font\seveneufm=eufm7 \font\fiveeufm=eufm5 \newfam\eufmfam \textfont\eufmfam=\teneufm \scriptfont\eufmfam=\seveneufm \scriptscriptfont\eufmfam=\fiveeufm

\font\teneufb=eufb10 \font\seveneufb=eufb7 \font\fiveeufb=eufb5 \newfam\eufbfam \textfont\eufbfam=\teneufb \scriptfont\eufbfam=\seveneufb \scriptscriptfont\eufbfam=\fiveeufb

\font\teneurm=eurm10 \font\seveneurm=eurm7 \font\fiveeurm=eurm5 \newfam\eurmfam \textfont\eurmfam=\teneurm \scriptfont\eurmfam=\seveneurm \scriptscriptfont\eurmfam=\fiveeurm

\font\teneurb=eurb10 \font\seveneurb=eurb7 \font\fiveeurb=eurb5 \newfam\eurbfam \textfont\eurbfam=\teneurb \scriptfont\eurbfam=\seveneurb \scriptscriptfont\eurbfam=\fiveeurb

\font\teneusm=eusm10 \font\seveneusm=eusm7 \font\fiveeusm=eusm5 \newfam\eusmfam \textfont\eusmfam=\teneusm \scriptfont\eusmfam=\seveneusm \scriptscriptfont\eusmfam=\fiveeusm

\font\teneusb=eusb10 \font\seveneusb=eusb7 \font\fiveeusb=eusb5 \newfam\eusbfam \textfont\eusbfam=\teneusb \scriptfont\eusbfam=\seveneusb \scriptscriptfont\eusbfam=\fiveeusb


\font\tenss=cmss10 \font\sevenss=cmss7 \font\fivess=cmss5 \inewfam\ssfam \textfont\ssfam\tenss \scriptfont\ssfam\sevenss \scriptscriptfont\ssfam\fivess

\font\sevenit=cmti7 \scriptfont\itfam=\sevenit

\font\tenbsy=cmbsy10
\font\twelvebsy=cmbsy12

\let\articletitle\seventeenss
\let\chaptertitle\twelvebf
\let\sectiontitle\tenbf
\let\subsectiontitle\tenbfit
\let\subsubsectiontitle\tenit
\let\contchaptertitle\tenbf 
\let\contsectiontitle\tenrm 
\let\contsubsectiontitle\sevenrm 
\let\contsubsubsectiontitle\fiverm 
\let\parnumfont\tenrm 
\let\parbackreffont\fiverm 
\let\proclaimfont\tenbf
\let\prooffont\tenit
\let\mainfont\tenrm


\def\lhighlight{} 
\def\suppresscomments#1#2#3{} 
\def\prohibitcomments#1#2#3{\errmessage{Draft comments are not allowed in the final version}} 

\ifx\format\undefined\else\format\fi 
\ifx\comment\undefined\def\comment{\prohibitcomments}\fi

\inewdimen\totalht
\inewdimen\totalwd
\def\shipbox#1{%
	\totalht\ht#1
	\advance\totalht2in
	\totalwd\wd#1
	\advance\totalwd2in
	\dumpbox\totalwd\totalht{#1}}

\inewtoks\cont
\def\contents{\the\cont}
\let\printcont\gobble 
\def\contlinechapt#1#2{\printcont{#2}\smallskip\everypar{}\noindent{\contchaptertitle\llink{chapter.#1}{#2}}\par}
\def\contlinesect#1#2{\printcont{#1.  #2}\everypar{}\indent{\contsectiontitle\llap{#1.\enskip}\llink{section.#1}{#2}}\par}
\def\contlinesubsect#1#2{\printcont{#1.  #2}\everypar{}\indent{\contsubsectiontitle{#1.\enskip}\llink{paragraph.#1}{#2}}\par}
\def\contlinesubsubsect#1#2{\printcont{#1.  #2}\everypar{}\indent\indent\indent{\contsubsubsectiontitle\llap{#1.\enskip}\llink{paragraph.#1}{#2}}\par}
\def\addcont#1#2#3{%
	\cont\expandafter{\the\cont#1}%
	\append\cont{{#2}}%
	\cont\expandafter{\the\cont{#3}}}

\newif\ifpresec \presecfalse 
\def\chapter#1\par{%
	\def\chapname{#1}%
	\parn0
	\presectrue
	\numbfalse
	\addcont\contlinechapt\chapname{#1}%
	\curverb{Chapter~}%
	\assignlabel{chapter}{\chapname}%
	\tchapter#1\par}
\def\tchapter#1\par{
	\everypar{\let\beforesect\beforesection}
	\chapbreak\bigbreak
	\centerline{\plabel\chaptertitle\anchor{chapter.#1}{#1}}%
	\nobreak\medskip
	\let\beforesect\relax 
}

\let\chapbreak\relax 

\inewcount\secn \secn0
\def\section#1\par{%
	\ifx\sectionid\undefined\advance\secn1 \edef\sectionid{\the\secn}\fi
	\everypar{\numpar}
	\parn0
	\presecfalse
	\numbfalse
	\addcont\contlinesect\sectionid{#1}%
	\curverb{\S}%
	\assignlabel{section}{\sectionid}%
	\tsection#1\par
}
\def\tsection#1\par{
	\beforesect\let\beforesect\beforesection
	\typesetsection{#1}%
	\aftersection	
	\let\sectionid\undefined
}
\def\beforesection{\vskip0pt plus.3\vsize \penalty-250 \vskip0pt plus-.3\vsize \bigskip \vskip\parskip}
\def\aftersection{\nobreak\smallskip}
\def\typesetsection#1{\leftline{\sectiontitle\ifx\sectionid\undefined\indent\else\hbox to \parindent{\hss\plabel\anchor{section.\sectionid}{\sectionid}\enspace\hfill}\fi#1}}
\let\beforesect\beforesection

\def\subsection#1\par{\bigbreak\numbtrue\curverb{\S}\subsectiontitle\noindent#1\/\mainfont\par\nobreak\medskip\addcont\contlinesubsect{\the\secn.\the\parn}{#1}}
\def\subsubsection#1\par{\bigbreak\numbtrue\curverb{\S}\subsubsectiontitle\noindent#1\/\mainfont\par\nobreak\medskip\addcont\contlinesubsubsect{\the\secn.\the\parn}{#1}}

\def\sskip#1{\ifdim\lastskip<\medskipamount \removelastskip\penalty#1\medskip\fi}
\def\slug{\hbox{\kern1.5pt\vrule width2.5pt height6pt depth1.5pt\kern1.5pt}}
\newif\ifqed
\def\qed{\unskip\nobreak\ \slug\global\qedtrue}
\def\proclaim{\numbtrue\gproclaim\proclaimfont{\endgraf\sskip{55}}} 
\def\proof{\gproclaim\prooffont{\ifqed\else\qed\fi\qedfalse\endgraf\sskip{-55}}}
\def\xproclaim#1.{\medbreak{\everypar{}\noindent}{\proclaimfont#1.\enspace}\atendpar{\endgraf\sskip{55}}\ignorespaces} 
\let\abstract\xproclaim 

\inewtoks\atendpar \atendpar{\endgraf} 
\inewtoks\atendbr 
\inewif\ifnumb \numbfalse 
\def\finishpar{\the\atendbr\atendbr{}\the\atendpar\atendpar{\endgraf}\numbfalse}
\def\ppar{\the\atendbr\atendbr{}\endgraf\numbfalse} 
\def\endlist{\endgraf{\parskip\smallskipamount\everypar{}\noindent}} 
\inewcount\parn \parn0
\inewif\ifparbref
\def\nextpar{\ifnumb\advance\parn1 \def\brt{}%
	\ifparbref\ifx\lastlabel\undefined\else\edef\brt{\noexpand\the\noexpand\csname rev\expandafter\string\lastlabel\endcsname}\fi\fi
	\assignlabel{paragraph}{\the\secn.\the\parn}\fi}
\def\numpar{\ifnumb\nextpar\expand{\noexpand\typesetparnum{\the\secn.\the\parn}\noexpand\typesetpbr{\brt}}\fi}
\inewdimen\pindent \pindent\parindent

\ifx\draftnum\undefined 
\def\gproclaim#1#2#3.{\curverb{#3~}
	\medbreak\ifnumb\nextpar\fi{\everypar{}\noindent}\plabel#1#3\ifnumb\ \anchor{paragraph.\the\secn.\the\parn}{}\the\secn.\the\parn\fi.\enspace\mainfont
	\atendpar={}\ifnumb\expand{\noexpand\typesetpbr{\brt}}\fi\atendpar=\expandafter{\the\atendpar#2}\ignorespaces}
\def\typesetparnum#1{\ifnumb{\plabel\parnumfont\anchor{paragraph.#1}{}#1.\enspace}\fi}
\def\typesetpbr#1{\ifnumb\def\brtext{#1}\ifx\brtext\empty\else\setetok\atendbr{{\parbackreffont Used in \noexpand\stripcomma\brtext.}}\fi\fi}
\else 
\let\numbfalse\relax \numbtrue
\def\gproclaim#1#2#3.{\curverb{#3~}\medbreak\atendpar={}\noindent#1#3.\enspace\mainfont\atendpar=\expandafter{\the\atendpar#2}\ignorespaces}
\inewdimen\brwidth \brwidth.6in
\parindent0pt \parskip1ex plus 1ex minus 1ex
\def\typesetparnum#1{\ifnumb\llap{\plabel\anchor{paragraph.#1}{}\parnumfont#1\enspace}\fi}
\def\typesetpbr#1{\ifnumb\def\brtext{#1}\ifx\brtext\empty\else
	\llap{\smash{\vtop{\everypar{}\raggedright\rightskip0pt plus 0pt \leftskip0pt plus 1fill \hsize\brwidth
	\parnumfont \strut \break 
	\parbackreffont\stripcomma#1}}\enspace}\fi\fi}
\fi

\def\hang{\hangindent\pindent}
\def\textindent#1{{\everypar{}\parindent\pindent\indent}\llap{#1\enspace}\ignorespaces}
\def\item{\endgraf\hang\textindent}

\def\filbreak{\endgraf\vfil\penalty-200\vfilneg}

\def\eject{\endgraf\break}
\def\supereject{\endgraf\penalty-20000}
\def\smallbreak{\endgraf\ifdim\lastskip<\smallskipamount
	\removelastskip\penalty-50\smallskip\fi}
\def\medbreak{\endgraf\ifdim\lastskip<\medskipamount
	\removelastskip\penalty-100\medskip\fi}
\def\bigbreak{\endgraf\ifdim\lastskip<\bigskipamount
	\removelastskip\penalty-200\bigskip\fi}

\inewdimen\plht \setbox0\hbox{\char`\_} \plht\ht0 \setbox1\hbox{.} \advance\plht-\ht1 
\catcode`\.\active \def.{\lower\plht\hbox{\char`\_}} \catcode`\.=12 
\def\email{\bgroup\catcode`.\active\xemail}
\def\xemail#1#2{\rlap{\hphantom{#2}@#1}#2\hphantom{@#1}\egroup}

\let\printref\gobble 
\inewcount\backref \backref0
\def\firstletter#1#2\end{#1}
\def\revref#1{\printref{revref \string#1}\expand{\csname rev\string#1\endcsname={\the\csname rev\string#1\endcsname, \noexpand\llink{backreference.\the\backref}{\ifpresec\expandafter\firstletter\chapname\end\else\the\secn\fi.\the\parn\ifnumb\else*\fi}}}}
\def\hinitrev#1#2#3{\printref{hinitrev \string#1 #2 text #3}
	\gdef#1{\printref{invoked \string#1 id #2 text #3}\advance\backref1 \printref{advance \string#1 \the\backref}\revref#1\anchor{backreference.\the\backref}{}\llink{#2}{#3}}} 
\def\pinitrev#1#2#3{\gdef#1{#3}} 
\let\initrev\hinitrev 
\def\xxstripcomma, {}
\def\xstripcomma{\if\ntok,\let\xcont\xxstripcomma\else\let\xcont\relax\fi\xcont}
\def\stripcomma{\futurelet\ntok\xstripcomma}

\let\printlabel\gobble 
\inewtoks\curverb 
\ifx\draftlabel\undefined
\def\plabel{}
\else
\def\labeltext{} 
\def\plabel{\ifx\labeltext\empty\else\smash{\llap{\parbackreffont\labeltext\quad}}\gdef\labeltext{}\fi} 
\fi
\def\assignlabel#1#2{
	\ifx\lastlabel\undefined\else
	\edef\labeltext{\expandafter\expandafter\expandafter\gobble\expandafter\string\lastlabel}
	\edef\tmp{\the\curverb}%
	\ifx\tmp\empty\else 
	\printlabel{verbal label with label text "\labeltext", label command "\expandafter\string\lastlabel", and label prefix "\the\curverb"}%
	\expandafter\edef\csname v\labeltext\endcsname{\the\curverb\expandafter\noexpand\lastlabel}\fi
	\printlabel{label with label text "\labeltext" and label command "\expandafter\string\lastlabel"}%
	\edef\tmp{\noexpand\lastlabel{#1.#2}{#2}}%
	\expandafter\expandafter\expandafter\initrev\tmp
	\fi\let\lastlabel\undefined}
\inewtoks\vlist
\def\label#1{\let#1\undefined\def\lastlabel{#1}\addver#1\verifylabel\numbtrue} 

\inewdimen\bibindent 
\inewtoks\bibt 
\def\tbib#1{\ifx#1\undefined\else\errmessage{Bibliography item \string#1 already defined}\fi
	\edef\key{\expandafter\gobble\string#1}%
	\addver#1\verifybib
	\expandafter\inewtoks\csname rev\string#1\endcsname
	\expand{\noexpand\initrev\noexpand#1{reference.\key}{\key}}%
	\setbox0=\hbox{[\key]}%
	\ifdim\bibindent<\wd0 \bibindent=\wd0\fi
	\append\bibt{\noexpand\typesetbib\noexpand#1\src}%
	\ftbib}
\def\ftbib#1\par{\bibt=\expandafter{\the\bibt#1\par}}
\inewif\iftype
\def\typesetbib#1#2\par{\typetrue
	\ifsuppressunusedbib\edef\tmp{\the\csname rev\string#1\endcsname}\ifx\tmp\empty
	\typefalse\fi\fi
	\iftype
	\edef\key{\expandafter\gobble\string#1}%
	\noindent\hbox to \bibindent{[\anchor{reference.\key}{\key}]\hfil}#2%
	\edef\bibbr{\the\csname rev\string#1\endcsname}%
	\ifx\bibbr\empty\else\expandafter\stripcomma\bibbr.\fi
	\hangindent\bibindent\filbreak\par
	\fi}

\def\addver#1#2{\append\vlist{\noexpand#2\noexpand#1{\the\inputlineno}}}
\def\verifyref#1#2#3{\edef\tmp{\the\csname rev\string#1\endcsname}\ifx\tmp\empty\warning{#3 \string#1 at line #2}\fi}
\def\verifylabel#1#2{\verifyref#1{#2}{Unused label}}
\def\verifybib#1#2{\verifyref#1{#2}{Unused reference}}

\let\printurl\gobble 
\inewtoks\urltext
\inewtoks\urlt
\inewif\ifpunct
\def\urldash{-}
\def\urltilde{{\tensy^^X}} 
\def\ndash{\def\urldash{--}}
\def\http://{\hfil\penalty900\hfilneg\urltext={http://}\urlt={http:/\negthinspace/}\punctfalse\urlgrab}
\def\https://{\hfil\penalty900\hfilneg\urltext={https://}\urlt={https:/\negthinspace/}\punctfalse\urlgrab}
\def\urlgrab{\catcode`\#=11 \futurelet\ntok\urldispatch}
\def\urldispatch{%
	\ifx\ntok~\let\proceed\urlcont\else
	\ifcat\noexpand\ntok\space\let\proceed\urlfinish\else
	\ifcat\noexpand\ntok\relax\let\proceed\urlfinish\else
	\let\proceed\urlcont
	\fi\fi\fi\proceed}
\def\append#1#2{\setetok#1{\the#1#2}}
\def\urlcont#1{\ifpunct\append\urltext\punctc\append\urlt\punctc\punctfalse\fi
	\ifx\ntok~\append\urltext{\noexpand~}\append\urlt\urltilde
	\else\if\ntok\ohash\append\urltext\ohash\append\urlt\#%
	\else\if\ntok_\append\urltext_\append\urlt\_%
	\else\if\ntok-\append\urltext-\append\urlt\urldash
	\else\if\ntok.\puncttrue\def\punctc{.}%
	\else\if\ntok,\puncttrue\def\punctc{,}%
	\else\if\ntok;\puncttrue\def\punctc{;}%
	\else\append\urltext{#1}\append\urlt{#1}%
	\fi\fi\fi\fi\fi\fi\fi\urlgrab}
\def\urlfinish{\catcode`\#=6 \hbox{\printurl{\the\urltext}\link{\the\urltext}{\the\urlt}}\ifpunct\punctc\punctfalse\fi\def\urldash{-}}
\def\idgrab{\futurelet\ntok\iddispatch}
\def\iddispatch{\ifcat\noexpand\ntok\space\let\proceed\urlfinish\else\if\ntok,\let\proceed\urlfinish\else\let\proceed\idcont\fi\fi\proceed}
\def\idcont#1{\ifpunct\append\urltext.\append\urlt.\punctfalse\fi
	\if\ntok.\puncttrue\def\punctc{.}\else\append\urltext{#1}\append\urlt{#1}\fi\idgrab}

\let\printgrab\gobble 
\inewtoks\grabname
\inewtoks\grabtoks 
\inewtoks\grabcseq 
\inewtoks\subsuptoks 
\inewtoks\dtoks 
\inewcount\grabsize
\inewif\ifgrabsubscript 
\inewif\ifgrabsupscript 
\def\grabsequence{\bgroup 
	\grabsubscripttrue\grabsupscripttrue 
	\grabstring}
\def\grabalpha{\bgroup 
	\grabsubscriptfalse\grabsupscriptfalse 
	\grabstring}
\def\grabingroup{\ifinfont\errmessage{Already inside a math token}\fi\append\grabtoks{\bgroup}\infonttrue}
\def\graboutgroup{\ifinfont\append\grabtoks{\egroup}\infontfalse\fi}
\def\grabstring#1#2#3{
	\aftergroup#3 
	\inewif\ifdefine \inewif\ifinfont
	\printgrab{}\printgrab{grab a string of type #1, typeset using font \string#2, with postcommand \string#3}%
	\grabname{#1}\def\grabfont{#2}\grabsize0 \grabtoks={}\grabingroup \append\grabtoks{#2}\grabcseq={}%
	\futurelet\ntok\grabdeflookahead}
\def\grabdeflookahead{\if=\noexpand\ntok 
	\definetrue\printgrab{defining}\expandafter\grabgobblefuturelet
	\else\printgrab{referencing}\definefalse\expandafter\grablookahead\fi}
\def\grabgobblefuturelet#1{\grabfuturelet} 
\def\grabfuturelet{\futurelet\ntok\grablookahead}
\def\grablookahead{\printgrab{futurelet token meaning: \meaning\ntok}%
	\let\ncom\grabfinish
	\if\bgroup\noexpand\ntok \printgrab{left brace, terminating}%
	\else \if\egroup\noexpand\ntok \printgrab{right brace, terminating}%
	\else \if\space\noexpand\ntok \printgrab{blank space, terminating}%
	\else \let\ncom\grabexamine \fi\fi\fi \ncom}
\def\grabexamine#1{\printgrab{grabexamine argument: \string#1, meaning \meaning#1}%
	\def\ncom{\grabfinish#1}%
	\ifcat$\ifcat*\string#1\fi$
		\ifcat\noexpand~\noexpand#1 \printgrab{active character \string#1, continuing}%
			\advance\grabsize1 \grabtoks\expandafter{\the\grabtoks#1}\append\grabcseq{\string#1}\def\ncom{\grabfuturelet}%
		\else \ifcat _\noexpand#1 \ifgrabsubscript\printgrab{subscript, continuing}%
			\graboutgroup \grabtoks\expandafter{\the\grabtoks#1}\subsuptoks{#1}\def\ncom{\grabsubsupfuturelet}%
						\else\printgrab{subscript, terminating}\fi
		\else \ifcat ^\noexpand#1 \ifgrabsupscript\printgrab{superscript, continuing}%
			\graboutgroup \grabtoks\expandafter{\the\grabtoks#1}\subsuptoks{#1}\def\ncom{\grabsubsupfuturelet}%
						\else\printgrab{superscript, terminating}\fi
		\else \ifcat a\noexpand#1 \printgrab{letter #1, continuing}%
			\advance\grabsize1 \grabtoks\expandafter{\the\grabtoks#1}\grabcseq\expandafter{\the\grabcseq#1}\def\ncom{\grabfuturelet}%
		\else \printgrab{nonactive character \string#1, terminating}\fi\fi\fi\fi
	\else \printgrab{command sequence \string#1, terminating}\fi
	\ncom}
\def\grabsubsupfuturelet{\futurelet\ntok\grabsubsuplookahead}
\inewcount\dig
\inewif\ifdigit
\def\grabsubsuplookahead{\printgrab{subsup futurelet token meaning: \meaning\ntok}%
	\if\bgroup\noexpand\ntok \printgrab{left brace, continuing}\let\ncom\grabentiresubsup%
	\else \if\egroup\noexpand\ntok \printgrab{right brace, continuing}\let\ncom\grabentiresubsup%
	\else \if\space\noexpand\ntok \errmessage{Blank space after \the\subsuptoks}%
	\else \let\ncom\grabsubsupexamine \fi\fi\fi \ncom}
\def\grabentiresubsup#1{\printgrab{subsup entire group added}\grabtoks\expandafter{\the\grabtoks{#1}}\grabfuturelet}
\def\grabremainderfuturelet{\futurelet\ntok\grabremainderlookahead}
\def\grabremainderlookahead{\printgrab{futurelet token meaning: \meaning\ntok}%
	\let\ncom\grabfinish
	\if\bgroup\noexpand\ntok \printgrab{left brace, terminating}%
	\else \if\egroup\noexpand\ntok \printgrab{right brace, terminating}%
	\else \if\space\noexpand\ntok \printgrab{blank space, terminating}%
	\else \let\ncom\grabremainderexamine \fi\fi\fi \ncom}
\def\grabremainderexamine#1{\printgrab{grabremainderexamine argument: \string#1, meaning \meaning#1}%
	\def\ncom{\grabfinish#1}%
	\ifcat$\ifcat*\string#1\fi$
		\ifcat\noexpand~\noexpand#1 \printgrab{active character \string#1, continuing}%
			\grabtoks\expandafter{\the\grabtoks#1}\def\ncom{\grabremainderfuturelet}%
		\else \ifcat _\noexpand#1 \ifgrabsubscript\printgrab{subscript, continuing}%
			\graboutgroup \grabtoks\expandafter{\the\grabtoks#1}\subsuptoks{#1}\def\ncom{\grabsubsupfuturelet}%
						\else\printgrab{subscript, terminating}\fi
		\else \ifcat ^\noexpand#1 \ifgrabsupscript\printgrab{superscript, continuing}%
			\graboutgroup \grabtoks\expandafter{\the\grabtoks#1}\subsuptoks{#1}\def\ncom{\grabsubsupfuturelet}%
						\else\printgrab{superscript, terminating}\fi
		\else \ifnum"8000=\the\mathcode`#1 \printgrab{math active character \string#1, continuing}%
			\grabtoks\expandafter{\the\grabtoks#1}\def\ncom{\grabremainderfuturelet}%
		\else \ifcat a\noexpand#1 \printgrab{letter #1, continuing}%
			\grabtoks\expandafter{\the\grabtoks#1}\def\ncom{\grabremainderfuturelet}%
		\else \printgrab{nonactive character \string#1, terminating}\fi\fi\fi\fi\fi
	\else \printgrab{command sequence \string#1, terminating}\fi
	\ncom}
\def\grabsubsupexamine#1{\printgrab{examining subsup argument \string#1, meaning \meaning#1}%
	\ifcat$\ifcat*\string#1\fi$
		\ifcat\noexpand~\noexpand#1 \printgrab{active character \string#1, continuing}%
			\grabingroup\append\grabtoks{\grabfont\noexpand#1}\let\ncom\grabremainderfuturelet
		\else\ifnum"8000=\the\mathcode`#1 \printgrab{math active character \string#1, continuing}%
			\grabingroup\append\grabtoks{\grabfont\noexpand#1}\let\ncom\grabremainderfuturelet
		\else\ifcat a\noexpand#1 \printgrab{letter #1, checking whether single or not}%
			\dtoks{#1}\def\ncom{\futurelet\ntok\grabsubsupsecondlookahead}%
		\else\printgrab{something else, inserting a single-character sub/superscript, continuing}
			\advance\grabsize2 \append\grabtoks{\bgroup\grabfont\noexpand#1\egroup}\append\grabcseq{\the\subsuptoks\string#1}%
			\def\ncom{\grabfuturelet}\fi\fi\fi
	\else \printgrab{command sequence \string#1, continuing}%
		\append\grabtoks{\noexpand#1}\let\ncom\grabfuturelet\fi
	\ncom}
\def\grabsubsupsecondlookahead{\def\ncom{\append\grabtoks{\the\dtoks}\grabfuturelet}%
	\ifcat a\noexpand\ntok \printgrab{not a single letter, grabbing the entire subsupscript}%
		\advance\grabsize1 \append\grabcseq{\expandafter\string\the\subsuptoks}
		\grabingroup
		\edef\tmp{\the\grabtoks\grabfont\the\dtoks}%
		\grabtoks\expandafter{\tmp}%
		\append\grabcseq{\the\dtoks}%
		\def\ncom{\grabfuturelet}%
	\else \printgrab{single letter, continuing}\fi\ncom}
\def\grabfinish{\printgrab{grabfinish}\graboutgroup
	\printgrab{grabsize=\the\grabsize, grabtoks=\the\grabtoks, grabcseq=\the\grabcseq}%
	\def\ncom{\the\grabtoks}
	\ifnum\the\grabsize=0 \errmessage{No string to grab}\fi
	\ifnum\the\grabsize>1 
		\ifdefine
			\ifparbref 
				\def\ncom{\anchor{\the\grabname.\the\grabcseq}{\the\grabtoks}}%
			\else
				\expandafter\ifx\csname\the\grabname.\the\grabcseq\endcsname\relax
					\expandafter\xdef\csname\the\grabname.\the\grabcseq\endcsname{\the\grabcseq}%
				\else
					\warning{Duplicate definition of a reference of type \the\grabname\space \the\grabcseq\space to \the\grabtoks}%
				\fi
			\fi
		\else
			\ifparbref 
				\expandafter\ifx\csname\the\grabname.\the\grabcseq\endcsname\relax
					\warning{Undefined reference of type \the\grabname\space \the\grabcseq\space to \the\grabtoks}%
				\else
					\def\ncom{\llink{\the\grabname.\the\grabcseq}{\the\grabtoks}}%
				\fi
			\fi
		\fi
	\fi
	\ncom
	\egroup}


\def\gen:{\http://libgen.io/get.php?md5=}

\def\arXiv:{\urltext={https://arxiv.org/abs/}\urlt={arXiv:}\punctfalse\idgrab}
\def\MR{\urltext={http://www.ams.org/mathscinet-getitem?mr=}\urlt={MR}\punctfalse\idgrab}
\def\Zbl:{\urltext={https://zbmath.org/?q=an:}\urlt={Zbl:}\punctfalse\idgrab}
\def\doi:{\ndash\urltext={http://dx.doi.org/}\urlt={doi:}\punctfalse\urlgrab}

\def\matrix#1{\null\,\vcenter{\normalbaselines
	\ialign{\hfil$##$\hfil&&\enspace\hfil$##$\hfil\crcr
		\mathstrut\crcr\noalign{\kern-\baselineskip}
		#1\crcr\mathstrut\crcr\noalign{\kern-\baselineskip}}}\,}
\def\sqmatrix#1{\null\,\vcenter{\normalbaselines
	\ialign{\hfil$##$&\enspace\hfil$##$\hfil\thinspace&$##$\hfil\crcr
		\mathstrut\crcr\noalign{\kern-\baselineskip}
		#1\crcr\mathstrut\crcr\noalign{\kern-\baselineskip}}}\,}
\def\cdbl{\def\normalbaselines{\baselineskip20pt \lineskip3pt \lineskiplimit3pt }}
\def\cd{\cdbl\let\dskip\relax\matrix}
\def\sqcd{\cdbl\let\dskip\enskip\sqmatrix}

\newcount\arrowsize \arrowsize3
\def\ltoarr#1{\mathop{\count0=#1 \loop\ifnum\count0>0 \smash-\mkern-7mu \advance\count0 -1 \repeat \mathord\rightarrow}\limits} 
\def\lto#1#2{\mathrel{\ltoarr{#1}^{#2}}} 
\def\lgetsarr#1{\mathop{\mathord\leftarrow \count0=#1 \loop\ifnum\count0>0 \mkern-7mu\smash-\advance\count0 -1 \repeat}\limits} 
\def\mapright#1{\smash{\lto\arrowsize{#1}}}

\def\mapdown#1{\dskip\Big\downarrow\rlap{$\vcenter{\hbox{$\scriptstyle#1$}}$}\dskip}
\def\mapup#1{\dskip\Big\uparrow\rlap{$\vcenter{\hbox{$\scriptstyle#1$}}$}\dskip}

\def\rx#1#2{\rlap{\kern #1pt \raise#1pt \hbox{#2}}}
\def\dottednearrow{\rx{-8}. \rx{-6}. \rx{-4}. \rx{-2}. \rx0. \rx2. \rx4. \kern6pt \raise7.7pt \hbox{$\nearrow$}}

\inewcount\forno \forno0
\def\arrno#1#2{\global\advance\arr1 \edef\eeqnno{\the\arr}%
	\global\advance\forno1 \edef\eforno{\the\forno}%
	\xdef#2{\noexpand\llink{equation.\eforno}{\eeqnno}}
	#1{(\anchor{equation.\eforno}{\eeqnno})}}

\inewbox\mdiag
\def\wrapdiagram{%
	\setbox\mdiag\vtop\bgroup
	\null 
	\vskip\baselineskip
	\inewcount\arr \arr0
	\baselineskip0pt
	\lineskip4pt
	\lineskiplimit4pt
	\let\par\cr
	\obeylines
	\halign\bgroup\hfil$\displaystyle##$\hfil\cr
	\ewrapdiagram}

\def\ewrapdiagram#1{#1
	\egroup
	\egroup
	\vskip0pt plus \dp\mdiag \penalty-250 \vskip0pt plus-\dp\mdiag 
	\hangafter-\dp\mdiag
	\divide\hangafter\baselineskip
	\advance\hangafter-2
	\hangindent-\wd\mdiag
	\advance\hangindent-2em
	\hbox to\hsize{\hfil\dp\mdiag0pt \box\mdiag}%
	\ignorespaces}

\let\printextract\gobble 
\long\def\extract#1\label#2{\ifx#2\relax\let\ncom\relax\else
	\ifx#2\undefined
		\expandafter\inewtoks\csname rev\string#2\endcsname
		\def#2{\blah\advance\backref1 \printextract{extract \string#2 \the\backref}\revref#2}%
		\expandafter\def\csname v\expandafter\gobble\string#2\endcsname{\blah\advance\backref1 \printextract{extract \string#2 \the\backref}\revref#2}%
		\let\ncom\extract
	\else\errmessage{Label \string#2 already defined}\fi\fi\ncom}

\long\def\extractlabels#1\par#2 {\extract#2\label\relax} 

\def\preprocess#1{ 
	\expandafter\extractlabels\input#1 
	\expandafter\gobblepar\input#1\relax 
	\par\vfill\supereject 
}

\def\importlabels#1{ 
	\let\initrev\pinitrev 
	\def\bib##1{\let##1\relax\append\bibt{\let##1\noexpand\undefined}} 
	\preprocess{#1}%
	\let\initrev\hinitrev
	\let\bib\tbib
	\cont={} 
	\vlist={} 
	\the\bibt 
	\bibt={} 
	\let\chapname\undefined \secn0 \backref0 
}

\def\blah{blah} 
\output{\setbox0\box255\setbox0\box\footins\deadcycles0} 
\let\bib\tbib 
\let\refs\relax 
\everypar{\numpar}
\def\par{\finishpar} 
\dimen0\hsize
\dimen1\vsize
\hsize\maxdimen
\vsize\maxdimen
\hbadness10000
\preprocess\jobname 
\let\blah\undefined 
\hsize\dimen0 
\vsize\dimen1
\hbadness1000

\def\importlabels#1{} 
\ifscroll\vsize\maxdimen\inewbox\abox\output{\setbox\abox\vbox{\unvbox255\unskip}\shipbox\abox}\else\output{\plainoutput}\fi 
\def\bib#1\par{} 
\def\refs{\raggedright\rightskip0em plus \maxdimen \advance\bibindent1em \everypar{}\the\bibt \vfil\eject} 
\def\revref#1{} 
\def\addver#1#2{} 
\def\addcont#1#2#3{} 
\def\prepass#1{} 
\parbreftrue 
\everypar{\numpar}
\let\chapname\undefined \secn0 \backref0 
\expandafter\gobblepar\input\jobname\relax
\par\vfill\supereject
\ifsuppressunusedbib\def\verifybib#1#2{}\fi
\the\vlist 

\egroup 
\end
}

\def\contlinesubsect#1#2{}
\def\contlinesubsubsect#1#2{}

\iftrue
\pretolerance-1
\tolerance59
\hbadness60
\vbadness150
\binoppenalty10000
\relpenalty10000
\fi

\hyphenation{Grothen-dieck}

\def\longleftrightarrows{\mathbin{\vcenter{\hbox{$\longrightarrow$}\vskip-1.5ex\hbox{$\longleftarrow$}}}}
\def\ii{{\edef\a{\let\noexpand\cfont\the\font}\a\ifx\cfont\tenbf\tenbsy\else\ifx\cfont\tfont\twelvebsy\else\tensy\fi\fi^^q}} 

\def\Z{{\bf Z}} 
\def\R{{\bf R}} 
\def\C{{\bf C}} 
\def\Sph{{\bf S}} 
\def\Zs{{\underline\Z}} 

\let\mcfont\bf 
\let\qcfont\it 
\let\dgfont\eurm 

\def\Fil{\mathop{\rm Fil}} 
\def\Gr{\mathop{\rm Gr}} 
\def\Ch{\mathop{\rm Ch}} 
\def\colim{\mathop{\rm colim}} 
\def\Fun{\mathop{\rm Fun}} 
\def\im{\mathop{\rm im}} 
\def\Ho{\mathop{\rm Ho}} 
\def\holim{\mathop{\rm holim}} 
\def\hocolim{\mathop{\rm hocolim}} 
\def\hocofib{\mathop{\rm hocofib}} 
\def\cofib{\mathop{\rm cofib}} 
\def\hofib{\mathop{\rm hofib}} 
\def\fiber{\mathop{\rm fib}} 
\def\Seq{\mathop{\rm Seq}} 
\def\bSeq{\mathop{\mcfont Seq}\nolimits} 
\def\iSeq{\mathop{\qcfont Seq}\nolimits} 
\def\dgSeq{\mathop{\dgfont Seq}\nolimits} 
\def\bFil{\mathop{\mcfont Fil}} 
\def\iFil{\mathop{\qcfont Fil}} 
\def\dgFil{\mathop{\dgfont Fil}} 
\def\dgMod{\mathop{\dgfont Mod}\nolimits} 
\def\dashdgMod{\hbox{-}{\dgfont Mod}} 
\def\comp{\mathop{\rm comp}} 
\def\quot{\mathop{\rm quot}} 
\def\tot{\mathop{\rm tot}} 
\def\End{\mathop{\rm End}\nolimits} 
\def\PSh{\mathop{\rm PSh}} 
\def\Map{\mathop{\rm Map}\nolimits} 
\def\Hom{\mathop{\rm Hom}\nolimits} 
\def\Mor{\mathop{\rm Mor}\nolimits} 
\def\cst{\mathop{\rm const}} 
\def\Gap{\mathop{\qcfont Gap}} 
\def\res{\mathop{\rm res}} 
\def\fres{\mathop{\rm fres}} 

\let\To\lto

\mathchardef\boxtimes="2802 

\mathchardef\cof="381A 
\mathchardef\fib="3810 
\def\afib{\mathrel{\tilde\fib}} 

\def\tC{{{\dgfont C}}} 
\def\tD{{{\dgfont D}}} 
\def\tE{{{\dgfont E}}} 

\def\abcat{{\tt A}} 
\def\dg{{\rm dg}} 
\def\Ndg{\N_\dg} 

\def\bA{{\mcfont A}} 
\def\bV{{\mcfont V}} 
\def\bMod{{\mcfont Mod}} 
\def\fMod{{\mcfont fMod}} 
\def\bC{{\mcfont C}} 
\def\bD{{\mcfont D}} 
\def\bW{{\mcfont W}} 
\def\bM{{\mcfont M}} 
\def\bone{{\mcfont 1}} 

\def\cC{{\cal C}} 
\def\cD{{\cal D}} 
\def\WG{{\cal W}_G} 
\def\P{{\cal P}} 
\def\cU{{\cal U}} 
\def\udg{{\rm Z^0}} 
\def\hodg{{\rm H^0}} 

\def\Mod{{\rm Mod}} 
\def\Alg{{\rm Alg}} 
\def\Oper{{\rm Oper}} 
\def\N{{\rm N}} 
\def\dgCat{{\rm dgCat}} 
\def\crf{{\rm Q}} 

\def\whX{\widehat X} 
\def\whY{\widehat Y} 
\def\whcrfX{\widehat{\crf X}} 
\def\whf{\widehat f} 
\def\ctp{\mathbin{\widehat{\otimes}_s}} 

\def\cO{{\cal O}} 
\def\D{{\rm D}} 
\def\Dfil{\D^{\rm fil}} 
\def\id{{\rm id}} 
\def\op{{\rm op}} 
\def\sSp{{\rm Sp}^\Sigma} 
\def\ss{\Delta} 
\def\sb{\partial\ss} 
\def\O{{\rm O}} 
\def\U{{\rm U}} 
\def\M{{\rm M}} 
\def\E{{\rm E}} 
\def\B{{\rm B}} 
\def\Gen{{\rm Gen}} 
\def\pt{{\rm pt}} 
\def\d{{\rm d}} 
\def\bL{{\bf L}} 
\def\loc{{\rm L}} 
\def\wr{{\rm wr}} 
\def\repr{{\rm repr}} 

\def\sSet{{\rm sSet}} 

\bib\Fact
Aldridge K. Bousfield.
Constructions of factorization systems in categories.
Journal of Pure and Applied Algebra 9:2 (1977), 207--220.
\MR0478159, \Zbl:0361.18001, \doi:10.1016/0022-4049(77)90067-6.

\bib\StarAut
Michael Barr.
*-autonomous categories.
Lecture Notes in Mathematics 752 (1979).

\bib\Calc
William~G.~Dwyer, Daniel~M.~Kan.
Calculating simplicial localizations.
Journal of Pure and Applied Algebra 18:1 (1980), 17--35.
\MR0578563, \Zbl:0485.18013, \doi:10.1016/0022-4049(80)90113-9.

\bib\Perv
Alexander~A.~Beilinson, Joseph Bernstein, Pierre Deligne.
Faisceaux pervers.
Ast\'erisque 100 (1982).
\MR0751966, \Zbl:0536.14011, \gen:1826AE6F2EEC89E9403ABB85AF30E52D.

\bib\DmodF
G\'erard Laumon.
Sur la cat\'egorie d\'eriv\'ee des $\cal D$-modules filtr\'es.
Lecture Notes in Mathematics 1016 (1983), 151--237.
\MR0726427, \Zbl:0551.14006, \doi:10.1007/bfb0099964.

\bib\AccCat
Michael Makkai, Robert Par\'e.
Accessible categories: the foundations of categorical model theory.
Contemporary Mathematics 104 (1989).
\MR1031717, \Zbl:0703.03042, \doi:10.1090/conm/104, \gen:B942999D06BE8D844CE4755DA2CCAF0D.

\bib\HMC
Takashi Kimura, Jim Stasheff, Alexander Voronov.
Homology of moduli of curves and commutative homotopy algebras.
The Gelfand Mathematical Seminars, 1993--1995 (1996), 151--170, Birkh\"auser Boston.
\MR1398921, \Zbl:0858.18010, \arXiv:alg-geom/9502006v2, \doi:10.1007/978-1-4612-4082-2_9.

\bib\ModCat
Mark Hovey.
Model categories.
Mathematical Surveys and Monographs 63 (1999).
\MR1650134, \Zbl:0909.55001, \doi:10.1090/surv/063, \gen:229756633EF1A320BED3109B1AC1AB52.

\bib\SymSpec
Mark Hovey, Brooke Shipley, Jeff Smith.
Symmetric spectra.
Journal of the American Mathematical Society 13:1 (1999), 149--208.
\MR1695653, \Zbl:0931.55006, \arXiv:math/9801077v2, \doi:10.1090/S0894-0347-99-00320-3.

\bib\MonAx
Stefan Schwede, Brooke~E.~Shipley.
Algebras and modules in monoidal model categories.
Proceedings of the London Mathematical Society 80:2 (2000), 491--511.
\MR1734325, \Zbl:1026.18004, \arXiv:math/9801082v1, \doi:10.1112/s002461150001220x.

\bib\Pres
Daniel Dugger.
Combinatorial model categories have presentations.
Advances in Mathematics 164:1 (2001), 177--201.
\MR1870516, \Zbl:1001.18001, \arXiv:math/0007068v1, \doi:10.1006/aima.2001.2015.

\bib\OpMay
G.~Maxwell~Kelly.
On the operads of J.~P.~May.
Reprints in Theory and Applications of Categories 13 (2005).
\MR2177746, \Zbl:1082.18009, \http://tac.mta.ca/tac/reprints/articles/13/tr13abs.html.

\bib\DerLim
Jan-Erik Roos.
Derived functors of inverse limits revisited.
Journal of the London Mathematical Society 73:1 (2006), 65--83.
\MR2197371, \Zbl:1089.18007, \doi:10.1112/s0024610705022416.

\bib\HoDG
Bertrand To\"en.
The homotopy theory of {\it dg\/}-categories and derived Morita theory.
Inventiones Mathematicae 167:3 (2007), 615--667.
\MR2276263, \Zbl:1118.18010, \arXiv:math/0408337v7, \doi:10.1007/s00222-006-0025-y.

\bib\LecDG
Bertrand To\"en.
Lectures on dg-categories.
Lecture Notes in Mathematics 2008 (2011), 243--302.
\MR2762557, \Zbl:1216.18013, \doi:10.1007/978-3-642-15708-0_5.

\bib\DQPB
Vladimir Dotsenko.
An operadic approach to deformation quantization of compatible Poisson brackets.
Journal of Generalized Lie Theory and Applications 1:2 (2007), 107--115.
\MR2320771, \Zbl:1120.53048, \arXiv:math/0611154v4, \doi:10.4303/jglta/s070203.

\bib\HTT
Jacob Lurie.
Higher topos theory.
Annals of Mathematics Studies 170 (2009).
\MR2522659, \Zbl:1175.18001, \arXiv:math/0608040v4, \doi:10.1515/9781400830558, \http://math.harvard.edu/~lurie/papers/HTT.pdf (April 9, 2017).

\bib\DCSH
Gunnar Carlsson.
Derived completions in stable homotopy theory.
Journal of Pure and Applied Algebra 212:3 (2008), 550--577.
\MR2365333, \Zbl:1146.55006, \arXiv:0707.2585v1, \doi:10.1016/j.jpaa.2007.06.015.

\bib\LR
Clark Barwick.
On left and right model categories and left and right Bousfield localizations.
Homology, Homotopy and Applications 12:2 (2010), 245--320.
\MR2771591, \Zbl:1243.18025, \arXiv:0708.2067v2, \doi:10.4310/hha.2010.v12.n2.a9.

\bib\HTHC
Carlos T. Simpson.
Homotopy theory of higher categories.
New mathematical monographs 19 (2012).
\MR2883823, \Zbl:1232.18001, \arXiv:1001.4071v1, \doi:10.1017/cbo9780511978111, \gen:97D4B9B748CF3370651A0054791C5D22.

\bib\Traces
Kate Ponto, Michael Shulman.
Traces in symmetric monoidal categories.
Expositiones Mathematicae 32:3 (2014), 248--273.
\MR3253568, \Zbl:1308.18008, \arXiv:1107.6032v2, \doi:10.1016/j.exmath.2013.12.003.

\bib\GVdual
Mitya Boyarchenko, Vladimir Drinfeld.
A duality formalism in the spirit of Grothendieck and Verdier.
Quantum Topology 4:4 (2013), 447--489.
\MR3134025, \Zbl:06256955, \arXiv:1108.6020v2, \doi:10.4171/qt/45.

\bib\DAGXII
Jacob Lurie.
Derived Algebraic Geometry XII: Proper morphisms, completions, and the Grothendieck existence theorem.
November 8, 2011.
\http://math.harvard.edu/~lurie/papers/DAG-XII.pdf.

\bib\AKT
Clark Barwick.
On the algebraic $K$-theory of higher categories.
Journal of Topology 9 (2016), 245--347.
\MR3465850, \Zbl:06563157, \arXiv:1204.3607v6, \doi:10.1112/jtopol/jtv042.

\bib\CddR
Bhargav Bhatt.
Completions and derived de Rham cohomology.
\arXiv:1207.6193v1.

\bib\CellCat
Michael Makkai, Ji\v r\'\i~Rosick\'y.
Cellular categories.
Journal of Pure and Applied Algebra 218:9 (2014), 1652--1664.
\MR3188863, \Zbl:1285.18007, \arXiv:1304.7572v1, \doi:10.1016/j.jpaa.2014.01.005.

\bib\DCFO
Pierre Schapira, Jean-Pierre Schneiders.
Derived category of filtered objects.
\arXiv:1306.1359v1.

\bib\DGstable
Lee Cohn.
Differential graded categories are k-linear stable \ii-categories.
\arXiv:1308.2587v2.

\bib\KelDG
Bernhard Keller.
On differential graded categories.
Proceedings of the International Congress of Mathematicians, Madrid, August 22--30, 2006, Volume~II, 151--190.
\MR2275593, \Zbl:1140.18008, \arXiv:math/0601185v5, \doi:10.4171/022-2/8.

\bib\Enhanced
Alexei Bondal, Mikhail Kapranov.
Enhanced triangulated categories.
Mathematics of the USSR Sbornik 70:1 (1991), 93--107.
\MR1055981, \Zbl:0729.18008, \doi:10.1070/SM1991v070n01ABEH001253.

\bib\DCIC
Saul Glasman.
Day convolution for \ii-categories.
Mathematical Research Letters 23:5 (2016), 1369--1385.
\MR3601070, \Zbl:06686299, \arXiv:1308.4940v4, \doi:10.4310/mrl.2016.v23.n5.a6.

\bib\PETS
Bhargav Bhatt, Peter Scholze.
The pro-\'etale topology for schemes.
Ast\'erisque 369 (2015), 99--201.
\MR3379634, \Zbl:06479630, \arXiv:1309.1198v2.

\bib\Enrich
Bertrand Guillou, J.~P.~May.
Enriched model categories and presheaf categories.
\arXiv:1110.3567v3.

\bib\SNIC
Giovanni Faonte.
Simplicial nerve of an ${\cal A}_\infty$-category.
Theory and Applications of Categories 32:2 (2017), 31--52.
\MR3607208, \arXiv:1312.2127v1.

\bib\OpCom
James~T.~Griffin.
Operadic comodules and (co)homology theories.
\arXiv:1403.4831v1.

\bib\NTT
Domenico Fiorenza, Fosco Loregi\`an.
t-structures are normal torsion theories.
Applied Categorical Structures 24:2 (2015), 181--208.
\MR3474895, \Zbl:1345.18011, \arXiv:1408.7003v2, \doi:10.1007/s10485-015-9393-z.

\bib\HA
Jacob Lurie.
Higher algebra.
September 18, 2017.
\http://math.harvard.edu/~lurie/papers/HA.pdf.

\bib\Operads
Dmitri Pavlov, Jakob Scholbach.
Admissibility and rectification of colored symmetric operads.
Journal of Topology (to appear).
\arXiv:1410.5675v3, \doi:10.1112/topo.12008.

\bib\Spectra
Dmitri Pavlov, Jakob Scholbach.
Symmetric operads in abstract symmetric spectra.
Journal of the Institute of Mathematics of Jussieu (to appear).
\arXiv:1410.5699v2, \doi:10.1017/S1474748017000202.

\bib\EnhGroth
Alberto Canonaco and Paolo Stellari.
Uniqueness of dg enhancements for the derived category of a Grothendieck category.
\arXiv:1507.05509v4.

\bib\Uniq
Valery~A.~Lunts, Dmitri~O.~Orlov.
Uniqueness of enhancement for triangulated categories.
Journal of the American Mathematical Society 23:3, 853--853.
\MR2629991, \arXiv:0908.4187, \doi:10.1090/s0894-0347-10-00664-8

\bib\Rot
Jacob Lurie.
Rotation invariance in algebraic K-theory.
September 2, 2015.
\http://math.harvard.edu/~lurie/papers/Waldhaus.pdf.

\bib\HTSP
Dmitri Pavlov, Jakob Scholbach.
Homotopy theory of symmetric powers.
Homology, Homotopy, and Applications 20:1 (2018), 359--397.
\arXiv:1510.04969v3, \doi:10.4310/HHA.2018.v20.n1.a20.

\bib\MdR
Dmitri Pavlov, Jakob Scholbach.
Filtered D-modules and filtered Omega-modules.
In preparation.

\inewtoks\title
\title{Enhancing the filtered derived category}
\centerline{\articletitle\the\title}
\medskip
\halign{&#\hfil\cr
Owen Gwilliam, \https://people.mpim-bonn.mpg.de/gwilliam/, Max Planck Institute for Mathematics, Bonn\cr
Dmitri Pavlov, \https://dmitripavlov.org/, Department of Mathematics and Statistics, Texas Tech University\cr
}
\metadata{\the\title}{Owen Gwilliam; Dmitri Pavlov}

\tchapter Contents

\the\cont

\abstract Abstract.
The filtered derived category of an abelian category has played a useful role in subjects
including geometric representation theory, mixed Hodge modules, and the theory of motives.
We develop a natural generalization using current methods of homotopical algebra,
in the formalisms of stable \ii-categories, stable model categories, and pretriangulated, idempotent-complete dg categories.
We characterize the filtered stable \ii-category $\iFil(\cC)$ of a stable \ii-category~$\cC$
as the left exact localization of sequences in~$\cC$ along the \ii-categorical version of completion
(and prove analogous model and dg category statements).
We also spell out how these constructions interact with spectral sequences and monoidal structures.
As examples of this machinery, we construct a stable model category of filtered $\cD$-modules
and develop the rudiments of a theory of filtered operads and filtered algebras over operads.

\section Introduction

The filtered derived category plays an important role in the setting of constructible sheaves, $\cD$-modules, mixed Hodge modules, and elsewhere.
Our goal is to revisit this construction using current technology for homotopical algebra,
extending it beyond the usual setting of chain complexes (or sheaves of chain complexes).
Put simply, we would like to describe the filtered version of a stable \ii-category,
generalizing this classical situation.
We also develop explicitly this machinery in the setting of model categories and dg categories,
so that it can be deployed in highly structured contexts and used in concrete computations.

There are several places where this kind of machinery would be useful.
For instance, it applies to filtered spectra (and sheaves of spectra).
In a different direction, this work allows one to work with filtered \ii-operads and
filtered algebras over \ii-operads in a clean way, as explained in~\vfiltoperads.

Several results here are undoubtedly well-known but seem to be unavailable for convenient reference.
In this introduction, we begin by describing the classical construction and the basic problem we pursue.
We then describe our main results and the structure of the paper and finish with a comparison to other work.

\subsection The classical construction

Let $\Z$ denote the integers, equipped with the usual total ordering by~$<$.
Let $\Zs$ denote the associated category, whose objects are integers and where $\Zs(m,n)$ is empty if $m>n$ and is a single element if $m\le n$.

Let $\abcat$ denote an abelian category.

\proclaim Definition.
The category of {\it sequences in~$\abcat$\/} is the functor category $\Fun(\Zs,\abcat)$.
We denote it by~$\Seq(\abcat)$.

\proclaim Definition.
The {\it filtered category of~$\abcat$}, denoted~$\Fil(\abcat)$, is the full subcategory of the category of sequences
in which an object $X\colon\Zs\to\abcat$ satisfies the condition that $X(m\to n)\colon X(m)\to X(n)$ is a monomorphism for every $m\le n$.

Given an object $X$ in $\Fil(\abcat)$, we view it as equipping the object $X(\infty) = \colim X$
with the filtration whose $n$th component is $X(n)$.
Thus, we only consider filtrations that are ``exhaustive'' in the classical terminology.

\label\tower
\proclaim Remark.
Often, people are interested in a {\it tower}, i.e., a sequence where one is interested in the limit $X(-\infty)$ (or homotopy limit) rather than colimit.
In algebra, one often studies examples where each structure map is an epimorphism, such as
$$\cdots \to \Z/p^n \Z \to \Z/p^{n-1}\Z \to \cdots\to\Z/p^2 \Z \to \Z/p\Z \to0.$$
In stable homotopy theory, one studies towers of spectra, such as the chromatic tower.
It is possible to view a tower as a sequence in our sense by working in the opposite category (see \vtowertoo).

It is well-known that $\Fil(\abcat)$ is additive but not abelian, which increases the complexity of homological algebra in this setting.

\proclaim Definition.
The {\it associated graded\/} functor $\Gr\colon\Seq(\abcat)\to\prod_\Z\abcat$ sends a sequence~$X$ to $\Gr X$ where $$(\Gr X)_n=X(n)/\im X(n-1\to n)$$ for all $n\in\Z$.

\proclaim Remark.
In \vGrasadjoint, we show that $\Gr$ is left adjoint to the functor that turns a list of objects $(A_n)_{n \in \Z}$ into a sequence where every structure map is zero.

We now consider the abelian category $\Ch(\abcat)$ of unbounded chain complexes in~$\abcat$.
We equip it with the quasi-isomorphisms as its class of weak equivalences.

\label\fildercat
\proclaim Definition.
The {\it filtered derived category of~$\abcat$}, denoted $\Dfil(\abcat)$, is the localization of $\Ch(\Fil(\abcat))$
(equivalently, $\Fil(\Ch(\abcat))$) with respect to the {\it filtered weak equivalences},
which are maps of sequences $f\colon X \to Y$ such that $\Gr f\colon\Gr X\to\Gr Y$ is an indexwise quasi-isomorphism.

A~core objective of this paper is to formulate and analyze a version of this construction that takes as input a stable \ii-category or stable model category.
We show, of course, that our construction applied to $\Ch(\abcat)$ recovers this classical construction at the level of homotopy categories.
We also explore how this lift of the classical construction interacts with, e.g., monoidal structures and spectral sequences.

\proclaim Remark.
This kind of construction is also useful in the setting of stable homotopy theory,
as the spectral sequence arising from a sequence of spectra $X = \cdots \to X(n) \to X(n+1) \to \cdots$ only depends on~$X$ up to the spectral version of filtered weak equivalence.
(We prove this assertion in Propositions~\qcspseq\ and~\mcspseq.)

\subsection Overview of the paper and its main results

The statements below admit both
a \ii-categorical interpretation (with ``category'' replaced by ``\ii-category'' and terms like ``colimit'' interpreted \ii-categorically), as explained in~\vfiltquasi,
and a model-categorical interpretation (with ``category'' replaced by ``model category'', ``colimit'' by ``homotopy colimit'', and so on), as explained in~\vfiltmodel.
As shown in \S\seccomparison, these two approaches are equivalent.
In \S\dgcat, we also construct a dg category of filtered objects and show it is equivalent to its quasicategorical and model-categorical cousins,
so we refrain from formulating the dg category version of our main result.

\proclaim Remark.
The reader will notice, just by counting pages, that the \ii-categorical treatment is noticeably quicker, thanks to the convenient machinery developed by Lurie in~[\HTT] and~[\HA].
The relative length in the model category section arises because we need to introduce and develop model categorical versions of some concepts, such as t-structures, in a convenient fashion.
We expect that these tools may be helpful for other purposes.

Consider the category of functors from~$\Zs$ to a locally presentable stable category~$A$, which we call {\it sequences in~$A$}.
A~crucial role is played, in both approaches, by a {\it completion\/} functor, which an endofunctor on sequences in~$A$.

\proclaim Theorem.
The following categories of {\it filtered objects in~$A$\/} are equivalent:
\item{(1)} sequences localized at graded equivalences;
\item{(2)} sequences localized at completion maps;
\item{(3)} the essential image of the completion functor;
\item{(4)} complete sequences;
\item{(5)} sequences reflectively localized at maps from the zero object to constant sequences.
\endlist
The {\it associated graded functor\/} is a continuous and cocontinuous functor to $\Z$-graded objects in~$A$.
It creates equivalences.

One reason for the importance of these filtered categories is the following relationship with spectral sequences.

\proclaim Proposition.
Given a presentable stable category~$A$ with a t-structure,
the functor that constructs a spectral sequence from a sequence in~$A$
factors through the category of filtered objects in~$A$.

The constructions play well with a symmetric monoidal structure.
If $A$ is closed symmetric monoidal, then so is the category of sequences in~$A$,
with the monoidal structure given by the Day convolution with respect to the addition on~$\Zs$.

\proclaim Theorem.
The Day convolution monoidal structure descends to a closed symmetric monoidal structure on filtered objects.
The associated graded functor is strong monoidal and strong closed.

We discuss issues of duals and dualizability in this context as well.

\proclaim Remark.
If one restricts to the subcategory of nonnegative integers and hence to sequences of the form $X(0) \to X(1) \to \cdots$,
then the Dold--Kan correspondence (and its generalization to stable quasicategories, see Theorem~1.2.4.1 in~[\HA]) provide an efficient approach to the kind of results we develop here.
But our interest includes unbounded sequences, which is not subsumed by those results.

Here is an outline of the paper's structure.
In \S\filtquasi\ of the paper, we rely on the formalism articulated by Lurie in~[\HTT] and~[\HA] to develop this notion;
a key point is showing that the completion functor is the localization functor for graded equivalences.
(It begins by filling in and fixing the assertions of Example~1.2.2.11 of Lurie~[\HA].)
In \S\filtmodel, we approach these questions via model categories,
in hopes of making it as easy as possible to apply these notions in concrete situations.
(These two sections can be read independently, according to the taste of the reader.
Comparing them, though, can be illuminating.)
In \S\seccomparison, we provide a precise comparison and methods for converting between the quasicategorical and model categorical settings.
In \S\dgcat, we also discuss differential graded (dg) categories, develop a version of the filtered dg category for a nice class of dg categories, and verify both that our dg construction agrees with the quasicategorical and model categorical versions and that it is the unique dg enhancement of the filtered derived category.
We finish by developing some quick applications:
we provide convenient model categories for filtered operads and filtered algebras over filtered operads,
and we construct a stable model category of filtered $\cD$-modules.

\subsection Notation

We work in the settings of model categories, dg categories, and \ii-categories
(which will mean quasicategories here, although nearly all arguments are model-independent).
To help make clear which setting applies,
we use italic letters like $\cC$ to denote \ii-categories,
bold letters like $\bA$ to denote model categories,
and calligraphic letters like $\tC$ to denote dg categories.
When working with \ii-categories, we will use terms like ``functor'' or ``colimit'' and mean the \ii-categorical notions (which are the only ones that make sense in that setting).
When working with model or dg categories, we will carefully distinguish between, e.g., ``colimit'' (the 1-categorical notion) and ``homotopy colimit.''

We work with categories equipped with various kinds of enhancements,
so that there are different types of morphisms (a plain set or an enriched hom or derived hom and so on).
For an ordinary category~$C$, we use $\Mor_C(X,Y)$ to denote the underlying set of morphisms from object~$X$ to object~$Y$.
If $C$ is enriched over a symmetric monoidal category $V^\otimes$, then we use $C(X,Y)$ to denote the morphisms that live in~$V$.
If $C^\otimes$ is a closed symmetric monoidal category, then $\Hom_C(X,Y)$ denotes the internal hom adjoint to $\otimes$.
For $\bC$ a relative category or a model category, we use $\Map(X,Y)$ to denote the derived mapping space from~$X$ to~$Y$.
If $\cC$ is an \ii-category, then $\cC(X,Y)$ denotes the space of morphisms.
Finally, if $\cC^\otimes$ is a closed symmetric monoidal \ii-category,
then we use $\Hom_\cC$ to denote the internal hom adjoint to~$\otimes$.

\subsection Relation to other work

There is a substantial literature using filtered derived categories,
and filtrations themselves are used wherever people want to do nontrivial homological algebra.
We thus do not provide a detailed survey of the literature here;
\vmainresult\ shows that we are enhancing the classical situation.
Our construction using stable model categories, due to its 1-categorical flavor,
is easy to compare to approaches via classical homological algebra (e.g., by Laumon~[\DmodF]).
For instance, the recent work of Schapira and Schneiders~[\DCFO] finds an abelian category
presenting the filtered derived category in much the way our model category construction does.

Carlsson~[\DCSH] studies a version of derived completion in the setting of stable homotopy theory.
His focus is on base-change of spectral modules along morphisms of ring spectra, and
his derived completion is a spectral generalization of $I$-adic completion in classical commutative algebra,
as he explains carefully.
(Bhatt~[\CddR] and Bhatt--Scholze~[\PETS] employ techniques similar to Carlsson's.)
Using the \ii-categorical Dold--Kan correspondence of Lurie, Theorem~1.2.4.1 of~[\HA],
one can convert Carlsson's cosimplicial constructions into constructions with sequences (in our terminology).
Lurie in~\S4 of~[\DAGXII] also develops a similar formalism.
Our constructions use the term ``derived completion'' in a more abstract setting
but apply in these algebraic settings,
although we do not examine here whether our results are helpful in those contexts.
Our work bears the same relationship to that of Carlsson or Bhatt and Scholze
as an abstract spectral sequence bears to $I$-adic completion.

More recently, Lurie (\S1.2.2, \S1.2.3 in~[\HA]) introduced filtered objects in stable quasicategories,
and in \S3.1, \S3.2 of~[\Rot] he studies some properties of filtered spectra, in particular, the strong monoidality of the associated graded functor.
Our focus here includes many aspects not covered there, such as natural symmetric monoidal structures and dualizability,
in addition to the comparisons with model categories.

Another treatment of filtrations in higher categories occurs in Barwick~[\AKT],
which provides a context in which algebraic $K$-theory obtains a characterizing universal property.
There, Barwick studies objects with finite-length filtrations in the setting of Waldhausen quasicategories
(i.e., quasicategories that have a zero object and are equipped with a class of ``ingressive morphisms'' that are closed under cobase changes
and contain all morphisms from the zero object --- and hence behave like inclusions).
From the point of view of $K$-theory, an extension (and more generally, finitely-filtered object) ought to be understood as equivalent to its constituents (respectively, associated-graded object).
By declaring all morphisms to be ingressive, any stable quasicategory provides a Waldhausen quasicategory.
Hence, Barwick's work intersects with ours in that setting, namely finite filtrations in a stable quasicategory;
however, our work explicitly includes filtrations extending infinitely in both directions,
and Barwick moves beyond the stable setting.

\subsection Acknowledgments

We thank Thomas Nikolaus and Benjamin Hennion for helpful conversations about localization,
Tobias Barthel and Nathaniel Stapleton for interesting examples of filtered spectra,
Fosco Loregi\`an for pointing out Bousfield's paper~[\Fact] to us,
Fran\c cois Petit and Lee Cohn for help with the literature,
and Jakob Scholbach for providing extensive feedback on this paper.
We are grateful for an anonymous reviewer's feedback, which led to improved articulation and arrangement of our results, and to the anonymous referee, whose questions and suggestions clarified the paper and led us to develop the dg results.
Our gratitude extends to the \link{http://ncatlab.org/nlab/show/HomePage}{$n$Lab} for being a wonderful resource.
A~part of this work was done while both authors participated in the ``Homotopy theory, manifolds, and field theories'' trimester at the Hausdorff Research Institute for Mathematics.
The second author was also supported by the Max Planck Institute for Mathematics in Bonn and the SFB 1085 (Higher Invariants) in Regensburg.

\label\filtquasi
\section Filtered objects in the language of stable \ii-categories

In this section, our arguments will be model-independent in nature, never relying on particular features of quasicategories.
Hence here we will use the term ``\ii-category'' and not ``quasicategory.''
In \S\seccomparison, we use ``quasicategory'' to emphasize that there we use quasicategories as an explicit model for \ii-categories.

Let $\cC$ be a stable \ii-category.
For $\phi\colon A \to B$ in $\cC$, let $\cofib(\phi)$ denote the pushout of $\phi$ and the zero map $0\colon A\to0$.
Let $\iSeq(\cC)$ denote the functor category $\Fun(\Zs,\cC)$.
We call such a functor a {\it sequence\/} in $\cC$.
The \ii-category $\iSeq(\cC)$ is automatically stable, as limits and colimits in functor categories are computed objectwise.

In this section, for a sequence $X$, we will use $X(n)/X(n-1)$ to denote
$$\cofib(X(n-1\to n) \colon X(n-1) \to X(n)),$$
for simplicity.
The functor $\cofib \colon \Fun(\Delta^1,\cC) \to \cC$ is left adjoint to the functor $A \mapsto (0 \to A)$,
which is the left Kan extension induced by the inclusion $\{1\} \to \Delta^1$.
(See Remarks 1.1.1.7 and 1.1.1.8 of~[\HA].)

Our focus in this paper is on the interplay between sequences and their associated graded objects.

\label\idefgraded
\proclaim Definition.
The {\it associated graded functor\/} $\Gr \colon \iSeq(\cC) \to \prod_{\Z} \cC$ sends a sequence $X$ to the graded object $(X(n)/X(n-1))_{n \in \Z}$.

The functor $\Gr$ is the composite
$$\iSeq(\cC) \to \prod_\Z \Fun(\Delta^1, \cC) \to \prod_\Z \cC,$$
where $X$ maps to $(X(n-1\to n) )_{n \in \Z}$ and then to $(\cofib(X(n-1\to n) ))_{n \in \Z}$.

\proclaim Definition.
A~morphism $f \colon X \to Y$ is a {\it graded equivalence\/} if $\Gr f$ is an equivalence.

The usual terminology is {\it filtered equivalence\/} but we find {\it graded\/} more descriptive.
In words, a morphism between sequence $f\colon X \to Y$ is a graded equivalence
if the induced map $$X(n)/X(n-1) \to Y(n)/Y(n-1)$$ is an equivalence for every $n$.

Let $\WG$ denote the collection of graded equivalences.

\proclaim Definition.
For $\cC$ a stable \ii-category, the {\it filtered \ii-category of~$\cC$\/} is the localization $\iSeq(\cC)[\WG^{-1}]$.
We denote it by $\iFil(\cC)$.

\proclaim Remark.
The word ``localization'' above might seem ambiguous: it can mean either inverting morphisms up to a homotopy,
or it can mean passing to the full subcategory of local objects.
As we show in \vcompisloc, the two notions coincide in our case.

\proclaim Theorem.
If $\cC$ possesses sequential limits, then the localization $\iFil(\cC)$ exists and is a stable \ii-category.

\proof Proof.
In \vicomp\ we construct $\iFil(\cC)$ by giving an explicit equivalence with the image of a left exact localization functor on $\iSeq(\cC)$,
namely the completion functor.
By Lemma~1.4.4.7 of Lurie~[\HA], such a localization is stable.

This terminology of ``filtered \ii-category'' is justified by the next result.
Let $\cD(\abcat)$ denote the derived \ii-category of a Grothendieck abelian category $\abcat$.
(See \S1.3.5 of Lurie~[\HA] for its construction.)

\label\mainresult
\proclaim Theorem.
The homotopy category of the filtered \ii-category $\iFil(\cD(\abcat))$ is equivalent to the classical filtered derived category $\Dfil(\abcat)$.

This result is a folklore theorem, undoubtedly well-known to many people, but it is not available in the literature so far we can tell.
One goal of this paper is to make the result available for reference.
To ease the comparison with more classical approaches, we use model categories as an intermediary,
so the proof is deferred to \S\seccomparison,
where it follows immediately from \vbFilvsiFil\ and \vthmstmodcat.

\label\icomp
\subsection Complete filtered objects

We will now construct the localization of $\iSeq(\cC)$ at the graded equivalences~$\WG$ (see \videfgraded).
The key notion will be that of completion of a sequence.
{\it We assume for the remainder of this section that $\cC$ possesses sequential limits.}

\proclaim Definition.
Let $X$ be a sequence, and let $X(-\infty) = \lim X$.
The {\it completion of~$X$\/} is the sequence $\whX$ with $\whX(n)=\cofib(X(-\infty) \to X(n))$.
Let $\comp \colon \iSeq(\cC) \to \iSeq(\cC)$ denote the {\it completion functor\/} sending $X$ to~$\whX$.

For $X$ a sequence, let $X(\infty)$ denote $\colim X$, which we view as the ``underlying object'' whose filtration is given by the sequence.
The completion of a constant sequence $X$ is thus always zero, since the map $X(-\infty)\to X(n)$ is an equivalence for all~$n$.
This observation means that completion typically changes the ``underlying object'' of a filtration~$X$:
indeed, a sequence is {\it complete\/} precisely when the map $X(\infty)\to\whX(\infty)$ is an equivalence.

This notion of completion is not the obvious generalization of the ``classical'' approach to completion,
so we will show that both approaches do coincide.
The classical completion $\overline{X}$ of an ordinary sequence (say, of chain complexes) goes as follows.
Let $\overline{X}(\infty)$ denote $\lim_{n \in \Zs} X(\infty)/X(n)$.
The classical completion of~$X$ is the sequence
$\overline{X}(n) = \fiber (\overline{X}(\infty) \to X(\infty)/X(n))$.

\label\xminusinf
\proclaim Lemma.
The fiber of the canonical map $X(\infty) \to \overline{X}(\infty)$ is $X(-\infty)$.
Hence $$\overline{X}(n) \simeq \cofib(X(-\infty) \to X(n)) = \whX(n),$$
so that $\overline{X} \simeq \whX$.
In particular, $\colim \overline{X}(n) \simeq \overline{X}(\infty) \simeq \whX(\infty)$.

Moreover, this result justifies {\it a posteriori\/} our use of the notation $\overline{X}(\infty)$,
as it is the colimit of the sequence~$\overline{X}$.

\proof Proof.
The first assertion is by the following computation:
$$\eqalign{
\fiber(X(\infty) \to \overline{X}(\infty)) &\simeq \fiber(\lim X(\infty) \to \lim X(\infty)/X(n)) \cr
&\simeq \lim \fiber(X(\infty) \to X(\infty)/X(n)) \cr
&\simeq \lim X(n)=X(-\infty), \cr
}$$
as desired.
\ppar
The second claim is then
$$\eqalign{
\overline{X}(n) & \simeq \fiber(\overline{X}(\infty) \to X(\infty)/X(n))\cr
&\simeq \fiber(\lim_{m<n} X(\infty)/X(m) \to \lim_{m<n} X(\infty)/X(n))\cr
&\simeq \lim_{m<n} \fiber (X(\infty)/X(m) \to X(\infty)/X(n))\cr
&\simeq \lim_{m<n} \cofib(X(m) \to X(n)),\cr
&\simeq \cofib(\lim_{m<n} X(m) \to X(n))\cr
&\simeq \cofib(X(-\infty) \to X(n)),\cr
}$$
as asserted. The final claim follows by analogous manipulations.

By the lemma, for every~$n$, there is a canonical morphism $X(n) \to \whX(n)$,
and these are compatible with the structure maps of the sequences $X$ and $\whX$.
In other words, we have the following.

\proclaim Corollary.
There is a natural transformation $\gamma\colon\id_{\iSeq(\cC)} \Rightarrow \comp$ via the natural map $X \to \whX$.

The following is now immediate.

\label\gammagreq
\proclaim Lemma.
For every sequence $X$, the map $\gamma(X)$ is a graded equivalence.
Indeed, for each pair $n < m$, the induced map
$$X(m)/X(n) \to \whX(m)/\whX(n)$$
is an equivalence.

\proclaim Definition.
A~sequence~$X$ is {\it complete\/} if $\gamma(X) \colon X \to \whX$ is an equivalence.

A~sequence~$X$ is complete if and only if $X(-\infty) \simeq0$.
For example, the ``step-sequence'' $\langle n, A\rangle$ of an object~$A$---with $\langle n, A\rangle(m)=0$ for $m < n$ and $\langle n, A\rangle(m) = A$ for $m \ge n$,
with the identity map for each morphism $\langle n, A\rangle(n \to n+1)$---is complete.

\label\powerseries
\proclaim Remark.
One should be careful here.
When working in a classical situation, such as with a filtered chain complex, one must bear in mind that $X(-\infty)$ is a homotopy limit and often not a limit.
Hence $X(-\infty) \simeq0$ should {\it not\/} be identified with the classical notion of a ``Hausdorff'' or ``separated'' filtration.
For instance, anticipating \vpoly\ below,
we note that the filtration of the polynomial ring $k[t]$ by powers of the ideal $(t)$ is separated in the classical sense,
but this sequence $X$ is not complete.
Its completion is $\whX \simeq k[[t]]$, so that
$$X(-\infty) \simeq \fiber(X(\infty) \to \whX(\infty)) \simeq (k[[t]]/k[t])[-1],$$
using the fact that fibers and cofibers agree up to suspension.

\proclaim Proposition.
The completion functor is a left exact localization of $\iSeq(\cC)$.

\proof Proof.
Any localization of a stable \ii-category is left exact because left adjoints between stable \ii-categories preserve finite limits.
We now show that completion is a localization.
Let $X$ be a sequence and $Y$ a complete sequence.
If $f \colon X \to Y$ is a map of sequences, then there is a natural map $f(-\infty) \colon X(-\infty) \to Y(-\infty) \simeq0$.
Hence, for every $n$, the map~$f(n)$ factors through $X(n)/X(-\infty)$, and so
$f$ factors as $$X \To4{\gamma(X)} \whX \To4{\whf} Y.$$
Thus $\comp$ is left adjoint to the inclusion of the full subcategory of complete sequences into $\iSeq(\cC)$.

We now want to relate this localization with localizing~$\iSeq(\cC)$ at~$\WG$.

\label\compisloc
\proclaim Lemma.
The functor $\comp$ realizes localization with respect to the graded equivalences~$\WG$.

\proof Proof.
A~map $f \colon X \to Y$ is a graded equivalence if and only if $\cofib(\Gr f) \simeq0$.
Now observe that $\Gr$ is an exact functor: in functor categories, like $\iSeq(\cC)$, colimits are computed objectwise,
and since colimits commute, it is manifest that $\Gr$ is right-exact and hence exact.
Thus, $f$ is a graded equivalence if and only if $\Gr( \cofib f) \simeq \cofib(\Gr f) \simeq0$,
which is equivalent to $\cofib f$ being a constant sequence.
\ppar
Similarly, as $\comp$ is a left adjoint,
we see that $\comp(f)$ is an equivalence if and only if $$\comp(\cofib f) \simeq \cofib(\comp f) \simeq0,$$
which is equivalent to $\cofib f$ being a constant sequence (we need $\cofib f(n) \simeq \cofib f(-\infty)$ for all n).
Thus, a map $f$ is an equivalence after completion if and only if it is a graded equivalence.

\proclaim Remark.
The above shows that $\iFil(\cC)$ can also be defined as the localization of $\iSeq(\cC)$ with respect to the completion maps $X\to\comp(X)$ for all objects~$X$.

\subsection Spectral sequences

Given a t-structure on a stable \ii-category $\cC$, its {\it heart\/} $\cC^\heartsuit$ is an abelian category.
In \S1.2.2 of~[\HA], Lurie explains how to construct a spectral sequence in $\cC^\heartsuit$
for each sequence $X \in \iSeq(\cC)$.
We obtain the following \ii-categorical generalization of the classical fact
that a filtered complex and its completion have the same spectral sequence.

\label\qcspseq
\proclaim Proposition.
For any sequence $X$, the completion map $\gamma(X) \colon X \to \whX$ induces an isomorphism of spectral sequences.

In other words, the construction of spectral sequences factors through $\iFil(\cC)$.

\proof Proof.
We freely use Lurie's terminology from \S1.2.2 of~[\HA].
The spectral sequence for $X$ is determined by the $\Zs$-complex associated to $X$,
which assigns to each $i\le j$ the cofiber $\cofib(X(i) \to X(j))$.
The $\Zs$-complex of~$\whX$ is precisely the same, by \vgammagreq, and
$\gamma(X)$ produces the isomorphism.

\proclaim Remark.
We want to emphasize that a t-structure on $\cC$ also provides two natural sequences for each object $X$
via the shifted truncations.
The first sequence is $$\cdots\to\tau_{\ge k+1}X\to\tau_{\ge k}X\to\tau_{\ge k-1}X\to\cdots,$$
which generalizes Whitehead towers of spectra to the setting of stable \ii-categories.
The second sequence is $$\cdots\to\tau_{<k+1}X\to\tau_{<k}X\to\tau_{<k-1}X\to\cdots,$$
which generalizes Postnikov towers
and can be obtained by taking the homotopy cofiber of the canonical map of the Whitehead filtration on~$X$ to the constant sequence on~$X$.
Thus every t-structure associates to each object two spectral sequences (though they only differ by indexing).

\label\towertoo
\proclaim Remark.
In \vtower, we mentioned that one might care about towers;
in other words, one might view a sequence $X$ as a way of describing its limit $\lim_n X$, rather than its colimit.
It is easy to accomplish this redirection, as follows.
The opposite category $\cC^\op$ of a stable \ii-category $\cC$ is also stable.
Consider the functor $\iota \colon \Zs \to \Zs^\op$ sending $n$ to $-n$.
Given a sequence $X \in \iSeq(\cC)$, let $X^\circ \in \iSeq(\cC^\op)$ denote the composite $X^\op \circ \iota$,
where $X^\op \colon \Zs^\op \to \cC^\op$ is $X$ viewed in the opposite categories.
Then $\colim X^\circ =X^\circ(\infty) \simeq \lim X$.
The spectral sequence for $X^\circ$ is the natural spectral sequence for a tower.
(This spectral sequence agrees with that for~$X$, up to some reindexing,
since cofibers in $\cC^\op$ agree, up to a shift, with cofibers in $\cC$.)

We would like to further clarify the relationship between Lurie's approach to spectral sequences and ours here,
although we have depended on his work in our discussion.
In \S1.2.2 of~[\HA], a key role is played by a correspondence between sequences in $\cC$ and ``chain complexes,''
whose category is denoted by $\Gap(\cC)$.
For the precise notion of a $\Zs$-complex, see Definition~1.2.2.2 and Remark~1.2.2.3 of~[\HA],
but roughly speaking, such a chain complex is a functor $\Zs\to\cC$ equipped with null homotopies of each composition $n\to n+1\to n+2$
(which encodes the idea that $d^2 = 0$)
along with appropriate higher coherences.
Using our completed category, we now establish a useful variant of Lurie's Lemma 1.2.2.4,
which characterizes the \ii-category $\Gap(\cC)$ of $\Zs$-complexes.

\proclaim Lemma.
Let $\Gap(\Zs,\cC)$ denote the \ii-category of $\Zs$-chain complexes in the stable \ii-category~$\cC$.
The restriction functor $\res \colon \Gap(\Zs,\cC) \to \iSeq(\cC)$
sending a complex to the underlying sequence
post-composes with localization to produce a functor
$\fres \colon \Gap(\Zs,\cC) \to \iFil(\cC)$,
and $\fres$ is an equivalence of \ii-categories.

\proclaim Remark.
Construction~1.2.2.6 in~[\HA] extracts a spectral sequence from an object in $\iFil(\cC)$
by passing to an object in $\Gap(\Zs,\cC)$ using Lemma~1.2.2.4.
Our variant can be used in a similar role.
A~crucial difference between our lemma and Lurie's is that our chain complexes are indexed by~$\Zs$ and not by $\hat\Zs:=\{-\infty\}\cup\Zs$.
Lemma~1.2.2.4 uses the element $-\infty$ to account for the homotopy limit of a possibly noncomplete sequence.

\proof Proof.
By Lemma~1.2.2.4 in~[\HA] we have an equivalence $\Gap(\hat\Zs,\cC)\to\iSeq(\cC)$ induced by the restriction map.
Furthermore, $\Gap(\Zs,\cC)$ is a (reflective) localization of $\Gap(\hat\Zs,\cC)$
whose local objects are precisely those $\hat\Zs$-chain complexes~$F$ for which $0\simeq F(-\infty,-\infty)\to\lim_{i\in\Zs}F(-\infty,i)$ is an equivalence.
Similarly, $\iFil(\cC)$ is a (reflective) localization of $\iSeq(\cC)$ whose local objects are complete sequences, i.e., the limit is the zero object.
We now observe that equivalences in both directions preserves these classes of local objects:
the restriction functor clearly preserves locality,
whereas the functor in the other direction yields a chain complex whose values at $(-\infty,k)$ give the original sequence,
so again locality is clearly preserved.
\ppar
Thus we have an equivalence between the localization of $\Gap(\hat\Zs,\cC)$ and $\iFil(\cC)$.
We complete the proof by showing that the restriction functor~$L$ induces an equivalence from the localization of $\Gap(\hat\Zs,\cC)$ to $\Gap(\Zs,\cC)$.
A~functor~$R$ in the other direction extends a $\Zs$-chain complex~$F$
to a $\hat\Zs$-chain complex~$\hat F$ by setting $$\hat F(-\infty,k)=\lim_{i\in\Zs}F(i,k)$$ and $\hat F(-\infty,-\infty)=0$.
(This is formalized in the obvious fashion using right Kan extensions.)
This construction lands in chain complexes because sequential limits preserve (co)fibers
and it lands in local objects because the diagonal $\Zs\to\Zs^{[1]}$ is a final functor.
The induced endofunctor~$LR$ on $\Zs$-chain complexes is equivalent to the identity functor.
It remains to show that the induced endofunctor~$RL$ on $\hat\Zs$-chain complexes is also equivalent to the identity functor.
The unit map is $\hat F\to RL(\hat F)$ and it is an equivalence on $(i,j)$ if $i\ne-\infty$ by construction
and also on $(-\infty,-\infty)$ (trivially).
Thus it remains to show it is an equivalence on $(-\infty,k)$.
The definition of a chain complex (Definition~1.2.2.2 in~[\HA]) immediately implies that the restriction of the unit to the full subcategory on $(-\infty,k)$ ($k\in\Zs$)
is a graded equivalence.
The source is complete by assumption and the target is complete by the above observation,
therefore it is also an objectwise equivalence, which completes the proof.

\subsection Symmetric monoidal structures

If $\cC$ is a presentable closed symmetric monoidal stable \ii-category,
then so is the \ii-category $\iSeq(\cC)$ with its natural symmetric monoidal product~$\otimes_s$ given by the Day convolution.
(See Proposition~2.9 of Glasman~[\DCIC].)
In essence, we have
$$(X \otimes_s Y)(n) \simeq \colim_{p+q \le n} X(p) \otimes Y(q),$$
which assembles all the elements of the ``double complex'' $X \otimes Y \colon \Zs^2 \to \cC$
below the line $p + q = n$,
so that \ii-categorical Day convolution is a direct generalization of
the classical formula from homological algebra.
The colimit functor $X \mapsto X(\infty)$ is a strong monoidal functor from $\iSeq(\cC)$ to $\cC$,
and thus the ``underlying object'' of a tensor product of sequences goes to the tensor product of the underlying objects.

\label\generator
\proclaim Remark.
One way to understand this symmetric monoidal product is by examining how it behaves on a collection of generators of $\iSeq(\cC)$.
A~natural choice for generators are the ``step-sequences'':
let $\langle m,A\rangle$ denote the sequence which is zero for $n<m$ and the fixed object~$A$ for $n\ge m$, with identity as the structure maps.
These generators form a \ii-subcategory $\Gen$ of all sequences, and
there is a relative Yoneda embedding~$\iSeq(\cC) \to \PSh(\Gen)$, which is a fully faithful functor from sequences to the \ii-category of presheaves on these generators.
Indeed, the representable functors for the $\langle m,A\rangle$ and the maps $\langle m,A\rangle\to\langle m,A'\rangle$
detect the $m$th component, whereas the Yoneda image of $\langle m,A\rangle\to\langle m+1,A\rangle$
detects the transition map from the $m$th to the $(m+1)$st component.
(Any map between generators is either zero or can be presented as the composition of these two maps.)
The monoidal product $\langle m,A\rangle\otimes_s\langle n,B\rangle$ is just $\langle m+n,A\otimes B\rangle$,
and this produces a monoidal product on presheaves on generators.
By the monoidal Yoneda embedding, this presheaf monoidal structure determines the tensor product of sequences.

Let $\iota \colon \iFil(\cC) \to \iSeq(\cC)$ denote the forgetful functor right adjoint to $\comp$.
We use it to induce a symmetric monoidal structure on the filtered \ii-category.

\proclaim Theorem.
The filtered \ii-category $\iFil(\cC)$ is a symmetric monoidal \ii-category via the {\it completed tensor product\/} $\ctp$, where
$$X \ctp Y = \comp(\iota X \otimes_s \iota Y)$$
for any $X$ and $Y$ in $\iFil(\cC)$.

\proof Proof.
By Proposition~2.2.1.9 and Example~2.2.1.7 of Lurie~[\HA], to prove our claim,
we need only to show that for any graded equivalence $f \colon X \to Y$ and any sequence $Z$,
$\id_Z \otimes_s f \colon Z \otimes_s X \to Z \otimes_s Y$ is a graded equivalence.
\ppar
Observe the following simple fact: given a commuting square
$$\cd{
X(n-1) &\mapright{}& X(n)\cr
\mapdown{f(n-1)}&&\mapdown{f(n)}\cr
Y(n-1) &\mapright{}& Y(n)\cr
}$$
in a stable \ii-category, the map $X(n)/X(n-1) \to Y(n)/Y(n-1)$ is an equivalence
if and only if the map $\cofib f(n-1) \to \cofib f(n)$ is an equivalence.
Thus we have another mechanism for identifying graded equivalences: check if the sequence $\cofib f$ is constant.
\ppar
Suppose $f \colon X \to Y$ is a graded equivalence and $Z$ is arbitrary.
We compute
$$\eqalign{
\cofib(\id_Z \otimes_s f)(n) &\simeq \colim_{p + q \le n} \cofib(Z(p)\otimes X(q) \To{13}{\id_{Z(p)} \otimes f(q)} Z (p) \otimes Y(q))\cr
&\simeq \colim_{p + q \le n} (Z(p) \otimes \cofib f(q))\cr
&\simeq \colim_{p + q \le n} (Z(p) \otimes \cofib f(q-1))\cr
&\simeq \colim_{p + q \le n} \cofib(Z(p)\otimes X(q-1)\To{13}{\id_{Z(p)} \otimes f(q-1)} Z (p)\otimes Y(q-1))\cr
&\simeq \cofib(\id_Z \otimes_s f)(n-1).
}$$
Hence, $\id_Z \otimes_s f$ is also a graded equivalence.

The following proposition is proved by Lurie for the case of filtered spectra as Proposition~3.2.1 in~[\Rot].

\label\strongmonoidal
\proclaim Proposition.
The associated graded functor $\Gr$ is strong monoidal, intertwining $\ctp$ on $\iFil(\cC)$ with the tensor product $\otimes_{\Gr}$ given by Day convolution on $\prod_\Z \cC$.

\proof Proof.
Recall from \vgenerator\ that the monoidal structure is determined by its behavior on the step-sequences $\langle m, A\rangle$,
due to the fully faithful relative Yoneda embedding~$Y$.
For formal reasons, the functor~$Y$ is strong monoidal.
(This follows from the general theory of Day convolutions, see \S3 in Glasman~[\DCIC].)
Given a cocontinuous functor~$F$,
a strong monoidal structure on~$F$ can be constructed on the symmetric monoidal \ii-category of generators,
and then canonically extended to the whole category by cocontinuity.
Observe that the associated graded functor $\Gr$ is cocontinuous.
The associated graded of $\langle m,A\rangle$ is $A$ in degree~$m$ and zero elsewhere, which allows us to lift the structure maps of the strong monoidal structure on~$\cC$ to~$\iFil(\cC)$.
For maps of generators, observe that there is only the zero map $\langle m,A\rangle\to\langle n,B\rangle$ for $m<n$.
If $m\ge n$, then such maps can be identified with map $A\to B$, and for these the functoriality conditions are satisfied because they are satisfied for~$\cC$.

The functor $\Gr$ also plays nicely with the internal homs in both categories.
Let $\Hom_\cC$ denote the internal hom in $\cC$ adjoint to $\otimes$.
Let $\Hom_{\Gr}$ denote the internal hom adjoint to $\otimes_{\Gr}$ on $\prod_\Z \cC$.
Observe the following explicit formula:
$$\Hom_{\Gr}(X,Y)_n=\prod_{m\in\Z}\Hom_\cC(X_m,Y_{m+n}).$$
Similarly, let $\Hom_{\Fil}$ denote the internal hom adjoint to $\ctp$ on $\iFil(\cC)$.
We now obtain a convenient explicit formula for $\Hom_{\Fil}$.

\proclaim Lemma.
For any $X, Y \in \iFil(\cC)$,
$$\Hom_{\Fil}(X,Y)(n) \simeq \int_{m\in\Zs}\Hom_\cC(X(m),Y(m+n))$$
where $\int_\bullet B(-,-)$ denotes the end of a bifunctor~$B$.

\proof Proof.
It suffices to show that both sides represent the same functor.
We want to show that for any object $A$ in $\cC$ we have
$$\cC(A,\Hom_{\Fil}(X,Y)(n)) \simeq \cC(A,\int_{m\in\Zs}\Hom_\cC(X(m-n),Y(m))).$$
Using hom-tensor adjunctions and the fact that ends are limits, we transform this desired equivalence to
$$\cC(\langle n,A\rangle\otimes X,Y)\simeq\int_{m\in\Zs}\cC(A\otimes X(m-n),Y(m)).$$
The equality now expresses the fact that the space of natural transformations between two functors can be computed as an end.
In the last transition we commuted the continuous functor $\cC(A,-)$ past the end, which is a limit,
and used the hom-tensor adjunction.

\label\strongclosed
\proclaim Proposition.
The associated graded functor~$\Gr$ is strong closed, meaning the canonical morphism
$$\Gr(\Hom_{\Fil}(X,Y))\to\Hom_{\Gr}(\Gr X,\Gr Y),$$
which is adjoint to the composition $$\Gr X\otimes_{\Gr}\Gr(\Hom_{\Fil}(X,Y))\to\Gr(X\ctp\Hom_{\Fil} (X,Y))\to\Gr Y,$$
is an equivalence.

\proof Proof.
As in the previous proof, it suffices to verify the case when $X$~is a step-sequence $\langle m,A\rangle$.
Shifting the degree by~$m$, we may assume $m=0$.
We have $\Hom_{\Fil}(\langle0,A\rangle,Y)=Y^A$, where $Y^A$ denotes the objectwise internal hom $n\mapsto\Hom_{\cC}(A,Y(n))$.
Furthermore, $\Gr(\langle0,A\rangle)=A[0]$, which is a copy of $A$ in degree~0 and 0~elsewhere,
and $\Hom_{\Gr}(A[0],\Gr Y)=(\Gr Y)^A$, where the superscript~$A$ again denotes the objectwise internal hom, this time for graded objects.
Combining these observations together, we want to show that $\Gr(Y^A)\to(\Gr Y)^A$ is an equivalence,
which is true because $\Hom_{\Fil}(A,-)$ preserves finite (co)limits and $\Gr$ is defined using a cofiber.

\subsection Dualizability and reflexivity

\proclaim Convention.
In this section $\cC$ denotes a presentable closed symmetric monoidal \ii-category with a monoidal unit~$\bone$.

\subsubsection Abstract theory

Any symmetric monoidal \ii-category has an intrinsic notion of a dualizable object (see Ponto and Shulman~[\Traces] for an expository account and \S4.6 of~[\HA] for a quasicategorical treatment).
When the category is closed, however, the situation is slightly simpler and so we will work in that context.
Let $\cC^\otimes$ be a closed symmetric monoidal \ii-category,
and let $\Hom_\cC$ denote the internal hom adjoint to $\otimes$.
Let $\bone$ denote the unit object in $\cC$.
An object $X$ in $\cC$~is {\it dualizable\/} (or {\it rigid\/})
if and only if the canonical map $X\otimes\Hom_\cC(X,\bone)\to\Hom_\cC(X,X)$ is an equivalence (see~\S2 in~[\Traces]
or Proposition~10.2.(iii) in Boyarchenko and Drinfeld~[\GVdual]).
The dual $X^\vee$ is given by $\Hom_\cC(X,\bone)$.

\proclaim Remark.
It is possible to develop notions of dualizability in monoidal \ii-categories,
but to simplify the exposition below, we always assume~$\cC$ to be {\it symmetric}.
In this setting the square of the dualization functor $\Hom_\cC(-,\bone)$ is canonically isomorphic
to the identity functor once we restrict it to dualizable objects.

This notion of dualizability does not encompass all natural situations.
For example, the sheaf of filtered differential operators on a smooth variety in characteristic~0
is not dualizable as a filtered $\cO$-module, essentially because the associated graded is nonzero in infinitely many degrees.
In many ways, however, the sheaf of filtered differential operators behaves similarly to dualizable objects,
so it is desirable to capture the essence of this situation in a definition.

There is a more general notion of dual, given by a functor of the form $\Hom_\cC(-,D)$, where $D$ is a {\it reflecting object\/} that must be specified in advance.
For instance, consider the Verdier duality functor on a locally compact Hausdorff space~$M$.
Here one works with $\cC$-valued sheaves on such spaces.
For $M = \pt$, the reflecting object $D_\pt$ is simply the monoidal unit $\bone$ in $\cC$,
and for other such spaces $M$, the reflecting object~$D_M$ is taken to be the derived exceptional inverse image of $\bone$ along the map $M\to \pt$.

We want to introduce a more general version of dualization on filtered objects,
so we use this $D$-de\-pen\-dent notion and do not require {\it a priori\/} that $\Hom_\cC(-,D)$ is a contravariant autoequivalence (as is usually done for reflecting objects).
After all, one can always restrict to the full subcategory of objects~$X$ for which the canonical morphism $X\to\Hom_\cC(\Hom_\cC(X,D),D)$ is an equivalence,
and on this subcategory $\Hom_\cC(-,D)$ is a contravariant autoequivalence.

\proclaim Definition.
Let $D$ be an object in a closed symmetric monoidal \ii-category $\cC$.
An object~$X$ in $\cC$ is {\it$D$-reflexive\/} if the canonical map $X\to\Hom_\cC(\Hom_\cC(X,D),D)$ is an equivalence;
the object $\Hom_\cC(X,D)$ is the {\it$D$-reflector\/} of~$X$.

If every object is $D$-reflexive, we call the pair $(\cC,D)$ a symmetric {\it Grothendieck--Verdier category\/}, following Boyarchenko and Drinfeld~[\GVdual].
As the name suggests, there are many examples arising from geometry, such as flavors of bounded derived categories of constructible sheaves.
These were also studied by Michael Barr under the name of {\it *-autonomous categories}, see~[\StarAut].

\proclaim Remark.
Dualizable objects are $\bone$-reflexive, but the opposite is not true.
An elementary example is given by Boyarchenko and Drinfeld in Example~3.3 of~[\GVdual]:
take the closed symmetric monoidal category of finite-dimensional normed real vector spaces
with contractive maps (i.e., maps of norm at most~1) as morphisms.
All objects in this category are $\bone$-reflexive.
However, for any dualizable object~$X$, the trace of the identity map on~$X$ is a contractive map~$\bone\to\bone$
given by the multiplication by the dimension of~$X$, which implies that the dimension must be at most~1.

\subsubsection Graded reflexivity and dualizability

Given a reflecting object $D\in\cC$,
we can construct a reflecting object $D[0]\in\Gr(\cC)$
by putting $D$ in degree~0 and the zero object in all other degrees.
Using the formula
$$\Hom_{\Gr}(X,Y)_n=\prod_{m\in\Z}\Hom_\cC(X_m,Y_{m+n})$$
for the internal hom in graded objects,
one immediately computes an explicit expression for the $D[0]$-reflector:
$$\Hom_{\Gr}(X,D[0])_n=\Hom_\cC(X_{-n},D).$$
Hence we obtain the following.

\proclaim Proposition.
A~graded object is $D[0]$-reflexive if and only if its individual components are $D$-reflexive.

When no ambiguity can arise, we refer to $D[0]$-reflexive objects and $D[0]$-reflectors
as $D$-reflexive and $D$-reflectors, respectively.
Recall that being $\bone$-reflexive is a necessary condition for being dualizable,
and the $\bone$-reflector computes the dual in this case.

\proclaim Proposition.
If a graded object $M = (M_n)$ has finitely many nonzero components, each of which is $\cC$-dualizable, then $M$ is dualizable in $\prod_\Z \cC$.
As a partial converse, if the monoidal unit~$\bone$ is a compact object in~$\cC$,
then a dualizable object in $\prod_\Z \cC$ has only finitely many nonzero components, each of which is $\cC$-dualizable.

\proof Proof.
For the first part, start with a graded object~$M$ with finitely many nonzero components, each of which is $\cC$-dualizable.
Construct the dual~$N$ of~$M$ by placing the dual of~$M_{-n}$ in degree~$n$.
The counit map $M\otimes N\to1[0]$ is given by the counit maps $M_{-n}\otimes N_n\to1$.
The unit map $1[0]\to N\otimes M$ has as its target $\coprod_{n\in\Z}N_n\otimes M_{-n}$,
and only finitely many components are nonzero, hence by stability $\coprod=\prod$ in this case
and we can use the unit maps $1\to N_n\otimes M_{-n}$ as individual components.
The triangle identities follows from the triangle identities of~$M_n$.
\ppar
For the second part, recall that
an object~$X$ is dualizable if and only if the canonical map $$X\otimes\Hom_{\Gr}(X,\bone)\to\Hom_{\Gr}(X,X)$$ is an equivalence.
Via the explicit formulas for tensor products and homs of graded objects, the degree~$n$ component of that canonical map is
$$\coprod_{m\in\Z}\Hom_\cC(X_m,\bone)\otimes X_{m+n}\to\prod_{m\in\Z}\Hom_\cC(X_m,X_{m+n}).$$
For the canonical map to be an equivalence,
we conclude, using the next lemma below, that the individual components $$\Hom_\cC(X_m,\bone)\otimes X_{m+n}\to\Hom_\cC(X_m,X_{m+n})$$ must be equivalences for all $m$ and~$n$.
For $n=0$ we now show that individual components of~$X$ are dualizable.
The counit map in graded objects $$c\colon \bone[0]\to X\otimes\Hom_{\Gr}(X,\bone)$$
has degree~0 component $$c_0 \colon \bone\to\coprod_{m\in\Z}\Hom_\cC(X_m,\bone)\otimes X_m.$$
Consider the map
$$\pi_m \colon \coprod_{m\in\Z}\Hom_\cC(X_m,\bone)\to\Hom_\cC(X_m,\bone),$$
given by the identity map in degree~$m$ and the zero map in other degrees.
The composition $\pi_m \circ c_0$ is the counit map of~$X_m$.
On the other hand, the compactness of~$\bone$ implies that only finitely many of these counit maps are nonzero, and hence only finitely many of~$X_m$ are nonzero.

\proclaim Lemma.
Let $f \colon \coprod_{n\in\Z}A_n\to\prod_{n\in\Z}B_n$ be induced by maps $f_n \colon A_n\to B_n$.
If $f$ is an equivalence, then the components $f_n \colon A_n\to B_n$ are equivalences.

\proof Proof.
Suppose $g \colon \prod_{n\in\Z}B_n\to\coprod_{n\in\Z}A_n$ is an inverse to $f$.
Precomposing and postcomposing with the canonical maps $B_n\to\prod_{n\in\Z}B_n$ and $\coprod_{n\in\Z}A_n\to A_n$
gives us a collection of maps $g_n \colon B_n\to A_n$, which are the inverses of the given maps $f \colon A_n\to B_n$.

\subsubsection Reflexivity and dualizability of sequences

Recall that the monoidal unit~$\bone_s$ in sequences
is the step-sequence $\langle0, \bone \rangle$, which is zero in negative degrees and $\bone$, with identity structure map, in nonnegative degrees.

We start by computing the $\bone_s$-reflector of~$X$, which is the dual of~$X$ whenever $X$ is dualizable in $\iSeq(\cC)$.

\proclaim Lemma.
For any sequence~$X$, the $\bone_s$-reflector $\Hom_{\iSeq}(X,\bone_s)$ of~$X$ is equivalent to the sequence $$n\mapsto \Hom_\cC(X(\infty)/X(-n-1),\bone).$$
More generally, if $D$ is a reflecting object in $\cC$ and $D_s$ denotes the sequence $\langle 0, D \rangle$,
then the $D_s$-reflector of $X$ is equivalent to the sequence $n \mapsto \Hom_\cC(X(\infty)/X(-n-1),D)$.

\proof Proof.
Recall that for any step-sequence $\langle n, A \rangle$ and any sequence $X$, we have $$\iSeq(\langle n, A \rangle, X) \cong\cC(A, X(n)).$$
Thus, by Yoneda, we see that the functor $A \mapsto \iSeq(\langle n, A \rangle, X)$ is represented by $X(n)$.
We will use this observation to compute the $n$th value of the sequence $\Hom_{\iSeq}(X,\bone_s)$.
To do so, however, it is convenient to use the following observation.
\ppar
For a sequence~$Z$, let $Z_{[-1]}$ denote the sequence where $Z_{[-1]}(k)=Z(-1)$ for $k\ge0$, with identity as transition maps, and where $Z_{[-1]}(k)=Z(k)$ for $k<0$, with transition maps inherited from $Z$.
Consider the cofiber $Z^{[-1]}$ of the natural map $Z_{[-1]} \to Z$:
it has $$Z^{[-1]}(k)=\cofib( Z(-1) \to Z(k)) = Z(k)/Z(-1)$$ for $k\ge0$ and $Z^{[-1]}(k)=0$ for $k<0$.
Observe that
$$\iSeq(Z,\bone_s) = \cC(Z(\infty)/Z(-1), \bone) = \iSeq(Z^{[-1]}(\infty),\bone_s),$$
which is the main property we need of this operation $Z \mapsto Z^{[-1]}$.
\ppar
Combining the two observations above, we compute
$$\eqalign{
\iSeq(\langle n, A \rangle, \Hom_{\iSeq}(X,\bone_s))
&\simeq \iSeq(\langle n, A \rangle \otimes_s X, \bone_s) \cr
&\simeq \iSeq((\langle n, A \rangle \otimes_s X)^{[-1]}, \bone_s) \cr
&\simeq \cC((\langle n, A \rangle \otimes_s X)^{[-1]}(\infty),\bone) \cr
&\simeq \cC(A \otimes (X(\infty)/X(-n-1)), \bone) \cr
&\simeq \cC(A, \Hom_{\cC}(X(\infty)/X(-n-1),\bone)).}$$
By Yoneda, we thus see that $\Hom_{\iSeq}(X,\bone_s)(n)$ is given by $\Hom_\cC(X(\infty)/X(-n-1),\bone)$.
No special properties of~$\bone$ were used, so the proof for general~$D$ is the same.

The above computation shows that the monoidal unit~$\bone_s$ is not very convenient as a reflecting object;
both $X(\infty)$ and a shift appear in the expression for the reflector.

\proclaim Definition.
Consider a closed symmetric monoidal stable \ii-category~$\cC$ with a reflecting object~$D$.
Let $D_{\le0}$ denote the sequence that is zero for $n>0$ and is $D$ in {\it nonpositive\/} degrees, with identity maps as the structure maps.
The sequences $D_{\ge0}$ and $D_{>0}$ are defined analogously, with $D$ placed in degree~$n$ where $n\ge0$ or $n>0$, respectively.
Let $\cst(D)$ denote the constant sequence with value $D$.

\proclaim Lemma.
The $D_{\le0}$-reflector of a sequence~$X$ is the sequence $n\mapsto\Hom_\cC(X(-n),D)$.

\proof Proof.
Consider the exact sequence $$D_{>0}\to\cst(D)\to D_{\le0};$$
in other words, $D_{\le0}=\cofib(D_{>0}\to\cst(D))$.
Observe that
$$\Hom_{\iSeq}(X,\cst(D))=\cst(\Hom_\cC(\colim X,D))=\cst(\Hom_\cC(X(\infty),D)).$$
By the previous lemma, we compute
$$\Hom_{\iSeq}(X,D_{>0})(n) = \Hom_\cC(X(\infty)/X(-n),D).$$
(Note the shift because $D_{>0}$ is the sequence $\langle 1, D\rangle$.)
The $n$th value of the $D_{\le0}$-reflector of $X$ is then the cofiber of the map
$$\Hom_{\iSeq}(X,D_{>0})(n) \to \Hom_{\iSeq}(X,\cst(D))(n),$$
which is
$$\Hom_\cC(\fiber(X(\infty)\to X(\infty)/X(-n)),D)=\Hom_\cC(X(-n),D),$$
verifying the claim.

We have $\Gr(D_{\le0})=\Sigma D[1]$, so the associated graded of a $D_{\le0}$-reflector is the $\Sigma D[1]$-reflector, i.e., the suspension and shift of the $D$-reflector.
On the other hand, $\Gr(D_{\ge0})=D[0]$, so the associated graded of a $D_{\ge0}$-reflector is the $D$-reflector itself,
in particular, the associated graded preserves duals, as expected.

\subsubsection Filtered reflexivity and dualizability

For filtered objects we have the following criterion for reflexivity.
Let $R \in \iFil(\cC)$ be a choice of reflecting object.

\proclaim Proposition.
An object in $\iFil(\cC)$ is $R$-reflexive if and only if its associated graded is $\Gr R$-reflexive.
The associated graded functor sends $R$-reflectors to $\Gr R$-reflectors.

\proof Proof.
The associated graded functor is a strong monoidal and strong closed functor,
and so it preserves reflexive objects and sends reflectors to reflectors.
Furthermore, the associated graded functor reflects equivalences,
so that the morphism $X\to\Hom_{\Fil}(\Hom_{\Fil}(X,R),R)$ is an equivalence
if the associated graded of~$X$ is $\Gr R$-reflexive, i.e., $\Gr X\to\Hom_{\Gr}(\Hom_{\Gr}(\Gr X,\Gr R),\Gr R)$ is an equivalence.

\proclaim Proposition.
An object in $\iFil(\cC)$ is dualizable if and only if its associated graded is dualizable.
The associated graded functor preserves duals.

\proof Proof.
As in the previous proof, the associated graded functor is strong monoidal, strong closed, and reflects equivalences.
Dualizability of~$X$ is equivalent to the map $X\otimes\Hom_{\iFil}(X,\bone_s)\to\Hom_{\iFil}(X,X)$ being an equivalence,
and the dual of~$X$ is computed by $\Hom_{\iFil}(X,\bone_s)$, which implies our claims.

\label\filtmodel
\section Filtered objects in the language of stable model categories

We now develop an enhancement of the filtered derived category using model categories.
Our motivation is to provide concrete tools for manipulating and computing filtered objects
in highly structured situations, like symmetric spectra or chain complexes.

Let $\bA$ be a left proper combinatorial stable model category.
(See Chapter~7 of Hovey~[\ModCat] for a review of stable model categories and Appendix~A.2 of Lurie~[\HTT] for a review of combinatorial model categories.)
Let $\bSeq(\bA)$ denote the functor category $\Fun(\Zs,\bA)$ equipped with the projective model structure,
in which weak equivalences are natural transformations that are objectwise weak equivalences and in which fibrations are natural transformations that are objectwise fibrations.

\proclaim Lemma.
If a sequence $$\cdots\To4{x_{n-1}}X(n)\To4{x_n}X(n+1)\To4{x_{n+1}}X(n+2)\To4{x_{n+2}}\cdots$$ is cofibrant, then each map $x_n$ is a cofibration and all objects $X(n)$ are cofibrant.
More generally, for a cofibration $X\to Y$ of sequences, the individual components $X(n)\to Y(n)$ and the maps $$Y(n-1)\sqcup_{X(n-1)}X(n) \to Y(n)$$ are cofibrations.

\proof Proof.
Any cofibration is a retract of a transfinite composition of cobase changes of generating cofibrations,
so it suffices to check the claim on generating cofibrations and under these operations.
Let $\langle m, A\rangle$ denote the sequence given by $0$ for $n < m$ and for $n \ge n$ the object $A \in \bA$ with identity as every structure map.
(It is a kind of step function; these generate all sequences.)
In $\bSeq(\bA)$, a generating cofibration is a map $\langle m, f \rangle \colon \langle m, A\rangle \to \langle m, B\rangle$
which is $0$ for $n < m$ and the same map $f$ for every $n \ge m$, with $f \colon A \to B$ a generating cofibration in $\bA$.
Thus, a generating cofibration $\langle m, f \rangle$ is clearly a levelwise cofibration
and the other maps (involving the pushouts) are also cofibrations.
\ppar
Consider the inclusion $\{n\} \hookrightarrow \Zs$; restricting along this inclusion sends $X$ to $X(n)$.
Retracts, cobase changes, and transfinite composition in $\bSeq(\bA)$ all commute with restriction along this inclusion.
A~generating cofibration $\langle m, f \rangle$ also clearly restricts to a generating cofibration in $\bA$ (possibly trivial).
As cofibrations in $\bA$ are preserved under these three operations, we see that a cofibration $f$ in sequences restricts to a cofibration~$f(n)$ in~$\bA$.
\ppar
Similarly, consider the inclusion $(n-1 \to n) \hookrightarrow \Zs$; restricting along this inclusion sends $X$ to the arrow $X(n-1) \to X(n)$.
Again, generating cofibrations restrict to generating cofibrations along this inclusion,
and the three operations commute with this restriction.
This implies that the pushout map $$Y(n-1)\sqcup_{X(n-1)}X(n) \to Y(n),$$ arising from a cofibration of sequences, is a cofibration.

As a partial converse, for a sequence $X$ in which each $X(n)$ is cofibrant and each structure map $x_n$ is a cofibration,
if $X$ is bounded below---so that $X(k)=0$ for $k\ll0$---then $X$ is a cofibrant sequence.

\label\cofobj
\proclaim Remark.
There is a model structure on sequences in~$\bA$
for which the cofibrant objects are {\it precisely\/} the sequences with cofibrant components and cofibrations as transition maps.
Indeed, the class of morphisms that satisfies the two properties indicated in the above lemma is weakly saturated, as explained in the proof.
Furthermore, it lies between the projective cofibrations and the injective cofibrations
and is cofibrantly generated by Corollary~3.3 in Makkai and Rosick\'y~[\CellCat].
Thus this intermediate model structure with such cofibrations exists by the Smith recognition theorem.
We will not develop or use this model structure because it presents the same \ii-category as the projective model structure,
but the projective model structure is better developed and more convenient to apply.

\proclaim Example.
Consider $\bA = \Ch(\abcat)$, the category of unbounded chain complexes of a Grothendieck abelian category $\abcat$,
equipped with the injective model structure in which the weak equivalences are the qua\-si-iso\-mor\-phisms and the cofibrations are the levelwise monomorphisms.
(See Proposition~1.3.5.3 in Lurie~[\HA], in which $\bA$ is constructed and shown to be left proper, combinatorial, and stable.)
By the lemma, we see that a (projectively) cofibrant sequence is always a filtered chain complex,
and every bounded below filtered cochain complex is a cofibrant sequence.
Moreover, this lemma ensures that every sequence is weakly equivalent to a filtered chain complex, by applying a cofibrant replacement functor.
Hence, one might profitably view filtrations simply as tools for understanding sequences, rather than as intrinsically important.
(In the intermediate model structure introduced in the preceding remark, the filtered chain complexes are precisely the cofibrant sequences.)

Recall that the cofiber of a map $f \colon A \to B$ is the pushout of $0\gets A\to B$, where 0~is the zero object,
or equivalently, it is the coequalizer of $f\colon A\to B$ and the zero map $0\colon A\to B$.
(That is, it is a cokernel.)
The homotopy cofiber is the cofiber after cofibrantly replacing $f$ in the projective structure on the arrow category of $\bA$.

Let $\bSeq(\bA)^c$ denote the full subcategory of cofibrant sequences.
For $X \in \bSeq(\bA)^c$, each structure map~$x_n$ is a cofibration and so
the homotopy cofiber $\hocofib(X(n)\to X(n+1))$ is weakly equivalent to the cofiber $\cofib(X(n)\to X(n+1))$
because the model category~$\bA$ is left proper.

\proclaim Definition.
A~map $f\colon X \to Y$ of cofibrant sequences is a {\it graded equivalence\/} if for every $n$, the induced map $\cofib(X(n-1)\to X(n)) \to \cofib(Y(n-1)\to Y(n))$ is a weak equivalence.

It will be useful to extend this notion to noncofibrant sequences.
Let $W_G$ denote the collection of graded equivalences.

\proclaim Definition.
For $\bA$ a stable model category, the {\it filtered model category of~$\bA$\/} is the left Bousfield localization $\loc_G \bSeq(\bA)$ of the projective model structure along~$W_G$.
We denote it by $\bFil(\bA)$.

We would like a simple characterization of the weak equivalences in $\bFil(\bA)$.
Let $\overline{W}_G$ denote the subcategory of weak equivalences generated under the 2-out-of-3 condition by the
levelwise weak equivalences and graded equivalences between cofibrant sequences.

Fix a cofibrant replacement functor~$\crf$ on $\bSeq(\bA)$.

\proclaim Definition.
A~map $f\colon X \to Y$ of sequences is a {\it derived graded equivalence\/} if the map $\crf f \colon \crf X \to \crf Y$ is a graded equivalence.

\label\barwg
\proclaim Proposition.
A~map $f\colon X \to Y$ of sequences is in $\overline{W}_G$ if and only if it is a derived graded equivalence.
In particular, the class of derived graded equivalences does not depend on the choice of~$\crf$.

\proof Proof.
Given $f \colon X \to Y$, we have a commuting square
$$\cd{
X &\mapright{f}& Y\cr
\mapup{q_X \simeq} && \mapup{q_Y \simeq}\cr
\crf X &\mapright{\crf f} & \crf Y.\cr
}$$
If $f$ is a graded equivalence between cofibrants, then so is $q_Y \circ \crf f = f \circ q_X$, hence $\crf f$ must be as well.
Conversely, if $\crf f$ is a graded equivalence, then $f$ is in $\overline{W}_G$.

\proclaim Remark.
{\it A~priori\/} we know that the weak equivalences in $\bFil(\bA)$ will include the derived graded equivalences but might include more.
In \vlocalwe, however, we will show that weak equivalences in the left Bousfield localization indeed coincide with derived graded equivalences.

We now provide conditions under which this localization exists.
Recall that a model category is {\it left proper\/} if cobase changes along cofibrations are homotopy cobase changes.
Equivalently, one can say that weak equivalences are closed under cobase changes along cofibrations.
For the purposes of this section, a model category is {\it compactly generated\/}
if maps with a right lifting property with respect to all cofibrations between compact objects
are weak equivalences.
This is a very mild condition, for example, any left Bousfield localization of the projective model structure on simplicial presheaves on any small category has this property.
A~crucial property of compactly generated model categories is that weak equivalences are closed under filtered colimits, in particular, transfinite compositions.

\proclaim Theorem.
If a model category~$\bA$ is combinatorial, stable, and left proper, then so is~$\bFil(\bA)$.
In particular, the model structure on~$\bFil(\bA)$ exists.
In addition, if $\bA$ is compactly generated, then so is~$\bFil(\bA)$.

\proof Proof.
The desired model structure is the left Bousfield localization of the projective model structure on $\bSeq(\bA)$
with respect to the derived graded equivalences.
By Theorem~4.7 in Barwick~[\LR] such left Bousfield localizations exist for any left proper combinatorial model category
as soon as the subcategory of weak equivalences is accessible.
(The statement there requires a {\it set\/} of localizing morphisms.
Here one can take any set of morphisms that
generates derived graded equivalences under $\kappa$-filtered colimits,
where the subcategory of derived graded equivalences is $\kappa$-accessible and the model category~$\bSeq(\bA)$ is $\kappa$-combinatorial.
By \vlocalwe\ below, we can simply take maps between constant sequences $0\to K$, where $K$ runs over a set of homotopy generators for $\bA$.)
\ppar
The underlying category $\bA$ is locally presentable by assumption.
The subcategory of levelwise weak equivalences is accessible because so is the subcategory of weak equivalences in~$\bA$.
We wish to show $\overline{W}_G$ is an accessible subcategory,
meaning it is accessible as a category and the embedding is an accessible functor.
To this end we use Theorem~5.1.6 in Makkai and Par\'e~[\AccCat], which shows that the forgetful functor
from the bicategory of accessible categories and accessible functors to the bicategory of all categories and functors
creates products, inserters, and equifiers (the so-called ``PIE-limits'').
We will show that $\overline{W}_G$ arises as such a PIE-limit.
\ppar
Let $\crf$ be an accessible cofibrant replacement functor for~$\bSeq(\bA)$ constructed by the small object argument.
The functor $\mathord{\Gr}\circ\crf$ to the combinatorial model category $\prod_\Z\bA$ is accessible, and
hence the preimage of the class of weak equivalences in $\prod_\Z \bA$ is accessible.
By \vbarwg, we know that $\overline{W}_G$ is precisely this preimage.

From now on, we assume $\bFil(\bA)$ is equipped with the model category structure just discussed.
That is, we work with $\bA$ a left proper combinatorial stable model category.

\subsection Complete filtered objects

In this subsection, we provide a different perspective on where the filtered stable category comes from.
To be more precise, we provide a different class of morphisms, the {\it completion maps},
such that localization along them agrees with localization along the graded equivalences.

\subsubsection Completions

\proclaim Definition.
Given a collection $W$ of arrows on a category $\cC$, let $\overline{W}$ denote the closure of $W$ under the 2-out-of-3 relation.
Then $(\cC, \overline{W})$ is a category with weak equivalences.

The essence of our problem is to compare two different notions of weak equivalence on the category of sequences $\bSeq(\bA)$.
Let $W_L$ denote the collection of {\it levelwise\/} weak equivalences, which is the usual notion of weak equivalence on such a diagram category.
As above, $W_G$ denotes the collection of {\it graded\/} weak equivalences between cofibrant sequences.

The \ii-categories presented by the relative categories $(\bSeq(\bA), W_L)$ and $(\bSeq(\bA), \overline{W_G \cup W_L})$ are {\it not\/} equivalent.
It is immediate that $W_L \subset \overline{W_G \cup W_L}$, but the converse does not hold.
Consider the following classic example.

\label\poly
\proclaim Example.
We filter the polynomial algebra $k[t]$ by the powers of the ideal~$(t)$:
$$A(n) = \cases{k[t], & for $n\ge0$;\cr(t^{-n}), &for $n<0$.\cr}$$
The map $f\colon p(t) \mapsto (1+t)p(t)$ is filtration-preserving, and $\Gr f\colon \Gr A \to \Gr A$ is the identity.
However, $f$~itself is not invertible.

But these \ii-categories are very close to being equivalent, and we want to understand how to enlarge~$W_L$ in a natural way to obtain $\overline{W_G \cup W_L}$.
The key role is played by completion.

Fix a fibrant replacement functor $R$ on $\bSeq(\bA)$, and let $\holim Y = \lim R(Y)$ be its homotopy limit.
There is a natural map $\lim Y \to \holim Y$.
(The choice of~$R$ does not matter as any other choice will have a zigzag of natural weak equivalences to~$R$.)

Similarly, fix a fibrant replacement functor $R'$ on the arrow category $\Fun(\{0\to1\}, \bA)$ and define $$\hofib(\phi) = \fiber R'(\phi)$$ for $\phi\colon A \to B$ a morphism in~$\bA$.
(In a linear setting, the fiber $\fiber(\phi)$ is typically called the kernel.)
There is a natural map $\fiber(\phi) \to \hofib(\phi)$.

Finally, fix a cofibrant replacement functor $\crf$ on $\bSeq(\bA)$.

The following constructions play an important role for us.
\item{(1)} Consider the functor $\colim \colon \bSeq(\bA)^c \to \bA$.
We denote the image of a cofibrant sequence $X$ under this functor by $X(\infty)$,
and we view this object as equipped with a filtration via the sequence.
The object $X(\infty)$ is a representative of the homotopy colimit,
as we are working with the projective model structure.
(We are thus working only with exhaustive filtrations, in classical terminology.)
Note as well that the natural map $X(n)\to X(\infty)$ is a cofibration between cofibrant objects, so its cofiber is also the homotopy cofiber.
\item{(2)} Consider the functor $\quot \colon \bSeq(\bA)^c \to \bSeq(\bA)$ sending a sequence $X$ to the sequence
$$\quot X(n) = \cofib(X(n) \to X(\infty)).$$
We then define a functor $(-)^\wedge \colon \bSeq(\bA)^c \to \bA$
sending a sequence $X$ to $X^\wedge = \holim \quot(X)$.
\item{(3)} Finally, we define a {\it completion functor\/} $\comp \colon \bSeq(\bA)^c \to \bSeq(\bA)$ where the {\it completion of~$X$}, denoted~$\whX = \comp X$, is the sequence
$$\whX(n) = \hofib \left(X^\wedge \to \quot X(n) \right).$$
There is a natural transformation $c \colon \id \Rightarrow \comp$ arising from the natural (in~$X$) map of sequences $c_X\colon X \to \whX$
given by the natural (in $X$ and~$n$) map $X(n) \to \whX(n)$ for every~$n$. We call this the {\it completion map\/} for $X$.

\proclaim Remark.
There is another approach, following our construction with stable quasicategories.
Let $$X(-\infty) = \holim_{n \in \Zs} X(n)$$ and set $$\whX(n) = \hocofib(X(-\infty) \to X(n)).$$
By suitably modifying the arguments around \vxminusinf, one can show that this definition of~$\whX$ agrees, up to weak equivalence, with the completion functor we just defined.
The reason why we stick to the other approach is that the homotopy limits involved in its definition
can often be computed as ordinary limits, as explained below.
This allows for convenient calculation of the completion, without deriving any functors.

\proclaim Definition.
A~cofibrant sequence is {\it complete\/} if the natural map $X(\infty) \to \whX(\infty)$ is a weak equivalence
(equivalently, if $X(-\infty)\simeq0$).
The {\it derived completion\/} of a sequence $X \in \bSeq(\bA)$ is the completion of its cofibrant replacement~$\crf X$.
A~sequence is {\it derived complete\/} if the completion map $c_{\crf X} \colon \crf X \to \whcrfX$ is a weak equivalence.

\proclaim Remark.
In the case that a cofibrant sequence~$X$ is {\it bounded below},
so that there is some $N$ such that $X(n)=0$ for all $n < N$,
$\whX(\infty) = X(\infty)$ as $\quot X(n) = X(\infty)$ for all $n < N$.
Hence a bounded below sequence is complete.

This notion of completeness diverges from the usual notion
in that we work with the {\it homotopy\/} limit of the sequence of quotients $\cofib(X(n) \to X(\infty))$, rather than the limit (i.e., the limit in the sense of ordinary categories).
In certain situations, these coincide.
For instance, let $\abcat$ be an abelian category satisfying AB3 and AB4* and having a generator.
Let $\bA = \Ch(\abcat)$ be the category of unbounded complexes.
Then for a sequence $X$ in $\bA$ where each map $X(i) \to X(i+1)$ is an epimorphism, the homotopy limit $\holim X$ coincides with the limit $\lim X$.
(This result is Theorem~3.1 of Roos~[\DerLim].)

\proclaim Remark.
Our polynomial example~\poly\ is not complete in the traditional or homotopical sense.
The limit of the quotients is $\lim k[t]/(t^n) = k[[t]]$, and it is also the homotopy limit, as just remarked.
However, using the filtration of $k[[t]]$ by powers of the ideal $(t)$, the map $f\colon p \mapsto (1+t)p$ is invertible by the geometric series.
(Note the perhaps-surprising consequence that the homotopy limit of the sequence
$$\cdots \to (t^n) \to (t^{n-1}) \to \cdots \to k[t]$$
is not zero!
Rather, it can be computed as $(k[[t]]/k[t])[-1]$.
See \vpowerseries.)

\subsubsection The comparison result

\proclaim Definition.
Let $W_C$ denote the collection of completion maps for cofibrant sequences.

Recall that $W_G$ only involves morphisms between cofibrant sequences, so it does not contain $W_L$,
which allows arbitrary sequences.

\proclaim Lemma.
On $\bSeq(\bA)$, $\overline{W_G \cup W_L}$ is equal to $\overline{W_C \cup W_L}$.

\proclaim Remark.
In other words, once sequences are identified with their completions, levelwise and graded weak equivalences coincide.
This statement is analogous to the well-known fact that a filtered object and its completion have the same spectral sequence.

We immediately obtain the following consequence of the lemma.

\proclaim Corollary.
The \ii-categories presented by $(\bSeq(\bA), \overline{W_G \cup W_L})$ and $(\bSeq(\bA), \overline{W_C \cup W_L})$ are equivalent.

\proof Proof of lemma.
The proof relies on several lemmas proved below.
First, we show that every graded equivalence arises in some 2-out-of-3 relation with completions or levelwise equivalences.
Let $f\colon X \to Y$ a graded weak equivalence between cofibrant objects.
Applying completion, we obtain a commuting diagram
$$\cd{ X &\mapright{f}& Y\cr
\mapdown{c_X} && \mapdown{c_Y}\cr
\whX &\mapright{\whf} & \whY.\cr}$$
By \vgradedlevel, we see that $\whf$ is a levelwise equivalence.
Hence the composition $\whf \circ c_X = c_Y \circ f$ is in the closure of $W_L \cup W_C$ under composition, which forces $f$ to be in the closure, too.
\ppar
Second, by \vcomplgreq, we know that every completion map $c_X\colon X \to \whX$ factors into a graded equivalence followed by a levelwise equivalence.
Hence it arises in a 2-out-of-3 relation with graded and levelwise equivalences.

To simplify notation, we write $M/N$ for the cofiber of a map $f\colon M\to N$, where the map $f$ will be clear from context.
For a cofibrant sequence $A$, the structure maps $A(m) \to A(n)$ are cofibrations for all $m < n$.
Hence, for any $k < m < n$, the induced map $A(m)/A(k) \to A(n)/A(k)$ is a cofibration because it is the pushout of the cofibration $A(m) \to A(n)$ along the map $A(m) \to A(m)/A(k)$.

\label\finitegraded
\proclaim Lemma.
If $f\colon A \to B$ is a graded equivalence between cofibrant sequences, then the induced maps
$$\bar f_{mn}\colon A(n)/A(m) \to B(n)/B(m)$$
are weak equivalences for all $m < n$.

\proof Proof.
The proof is by induction on $k=n-m$.
As $\Gr f$ is a weak equivalence, we have the base case of~$k=1$.
For the induction step, consider the map of short exact sequences
$$\cd{0 &\mapright{}& A(n-k)/A(n-k-1) &\mapright{}& A(n)/A(n-k-1) &\mapright{} &A(n)/A(n-k)&\mapright{}&0\cr
&&\mapdown{}&&\mapdown{}&&\mapdown{}&\cr
0 &\mapright{}& B(n-k)/B(n-k-1) &\mapright{}& B(n)/B(n-k-1) &\mapright{} &B(n)/B(n-k)&\mapright{}&0.\cr}$$
(Alternatively, view it as a map of pushout squares.)
These short exact sequences are also homotopy exact due to the left pair of horizontal arrows being cofibrations,
so by left properness their cofibers are also homotopy cofibers.
Thus, as we know the leftmost and rightmost vertical maps are weak equivalences, the middle vertical arrow is also a weak equivalence.

\label\infgraded
\proclaim Lemma.
If $f\colon A \to B$ is a graded equivalence between cofibrant sequences,
then the induced maps $$\bar f_{m\infty}\colon A(\infty)/A(m) \to B(\infty)/B(m)$$ are weak equivalences for all~$m$.

\proof Proof.
Because colimits commute and, in this situation of cofibrant objects, the homotopy colimit coincides with the colimit, we find
$$A(\infty)/A(m) = \colim_{n>m} \left(A(n)/A(m)\right) \simeq \hocolim_{n>m} (A(n)/A(m)),$$
and likewise for $B(\infty)/B(m)$.
The preceding lemma then ensures that
$$\hocolim_{n>m} \bar f_{mn}\colon \hocolim_{n>m} (A(n)/A(m)) \to \hocolim_{n>m} (B(n)/B(m))$$
is a weak equivalence.

\label\infeq
\proclaim Lemma.
If $f\colon A \to B$ is a graded weak equivalence of complete cofibrant sequences, then the morphism $f(\infty)\colon A(\infty) \to B(\infty)$ is a weak equivalence.

\proof Proof.
By \vinfgraded, we know that
$$\holim_n A(\infty)/A(n) \to \holim_n B(\infty)/B(n)$$
is a weak equivalence because the morphism of underlying diagrams is an objectwise weak equivalence.
By the definition of completeness, we then see that $f(\infty)$ is a weak equivalence.

\label\gradedlevel
\proclaim Lemma.
If $f\colon A \to B$ is a graded weak equivalence of complete cofibrant sequences, then the map $f(n)\colon A(n) \to B(n)$ is a weak equivalence for all $n$.
Hence, $f$ is a levelwise weak equivalence.

\proof Proof.
Consider the diagram
$$\cd{0 &\mapright{}& A(n) &\mapright{}& A(\infty) &\mapright{} &A(\infty)/A(n)&\mapright{}&0\cr
&&\mapdown{}&&\mapdown{}&&\mapdown{}&\cr
0 &\mapright{}& B(n) &\mapright{}& B(\infty) &\mapright{} &B(\infty)/B(n)&\mapright{}&0.\cr}$$
By \vinfeq, we know $f(\infty) \colon A(\infty) \to B(\infty)$ is a weak equivalence.
We also know the rightmost vertical arrow is a weak equivalence by \vinfgraded.
Because $A(n)\to A(\infty)$ and $B(n)\to B(\infty)$ are cofibrations, the right map is the homotopy cofiber of the left two maps.
But by stability, the left map is the homotopy fiber of the right two maps, which are weak equivalences.
Hence, the leftmost vertical arrow is also a weak equivalence.

\label\complgreq
\proclaim Lemma.
For every cofibrant sequence $X$, the canonical map $X \to \crf(\whX)$ is a graded equivalence.

\proof Proof.
Observe that the completion map $c_X$ admits a functorial factorization $X \cof \crf(\whX) \afib \whX$,
where the first map is a cofibration and the second map is an acyclic fibration (in $\bSeq(\bA)$ with the projective model structure).
We know that
$$\cofib(\crf(\whX)(n) \to \crf(\whX)(n+1)) \simeq \hocofib(\whX(n) \to \whX(n+1)).$$
Since $\bA$ is stable, the homotopy cofiber agrees with the homotopy fiber, up to shift.
Moreover, homotopy (co)fibers commute.
Applying these facts repeatedly, we compute
$$\eqalign{
\hocofib(\whX(n) \to \whX(n+1))
&\simeq \Sigma \hofib( \hofib(X^\wedge \to X(\infty)/X(n)) \to \hofib(X^\wedge \to X(\infty)/X(n+1)))\cr
&\simeq \Sigma \hofib( \hofib(X^\wedge \to X^\wedge) \to \hofib(X(\infty)/X(n) \to X(\infty)/X(n+1)))\cr
&\simeq \Sigma \hofib(0 \to \hofib(X(\infty)/X(n) \to X(\infty)/X(n+1)))\cr
&\simeq \hofib(X(\infty)/X(n) \to X(\infty)/X(n+1)) \cr
&\simeq \Omega \hocofib(\hocofib(X (n) \to X(\infty)) \to \hocofib(X (n+1) \to X(\infty))) \cr
&\simeq \Omega \hocofib(\hocofib(X (n) \to X(n+1)) \to \hocofib(X (\infty) \to X(\infty))) \cr
&\simeq \Omega \hocofib(\hocofib(X (n) \to X(n+1)) \to0) \cr
&\simeq \hocofib(X(n) \to X(n+1)).\cr
}$$
Hence, the canonical map is a graded equivalence, as desired.

\subsubsection Completion and derived graded equivalences

We finish this subsection by identifying the weak equivalences in $\bFil(\bA)$ with derived graded equivalences.
Although derived graded equivalences are by definition contained in the weak equivalences of~$\bFil(\bA)$,
{\it a priori\/} it might be possible that the process of left Bousfield localization introduces new weak equivalences.
We show that this is not the case in \vlocalwe.

To prepare for the proof of that Proposition,
we recall some important machinery from the theory of model categories and then prove some technical lemmas.

In a model category $\bA$, there is a simplicial set known as the {\it derived mapping space\/} $\Map(X,Y)$
(sometimes also denoted RMap)
between any two objects $X$ and $Y$ in $\bA$.
See Chapter~17, and in particular \S17.5.15, of Hirschhorn~[\ModCat], who calls it
the {\it functorial two-sided homotopy function complex}.
This construction encodes important information about weak equivalences:
Theorem~17.7.7 in~[\ModCat] states that a map $f \colon X \to Y$ is a weak equivalence in $\bA$
if and only if the induced map $\Map(W,X) \to \Map(W,Y)$ is a weak equivalence of simplicial sets for every object $W$ in $\bA$,
and likewise if and only the induced map $\Map(Y,Z) \to \Map(X,Z)$ is a weak equivalence of simplicial sets for every object $Z$ in $\bA$.

This construction is also central in Bousfield localization.
Recall that for a morphism $f \colon X \to Y$ in~$\bA$, an object~$Z$ is {\it $f$-local\/} if the associated map $\Map(Y,Z) \to \Map(X,Z)$ is a weak equivalence.
Similarly, for any object~$Z$, a morphism $f \colon X \to Y$ is {\it $Z$-local\/} if the associated map $\Map(Y,Z) \to \Map(X,Z)$ is a weak equivalence.
Given a collection~$S$ of morphisms in a model category~$\bA$,
a map~$f$ is a weak equivalence in the left Bousfield localization $\loc_S \bA$
if it is local with respect to every object~$X$ that is itself local with respect to every morphism in~$S$.

Besides these general tools, we need two lemmas.

The first concerns simple constructions with sequences.
For any object~$K$ in~$\bA$, let $\cst(K)$ denote the constant sequence for $K$: $\cst(K)(n) = K$ for all $n$ and every structure map is the identity.
For any sequence~$X$, let $X(-\infty)=\holim X$ and $X(\infty)=\hocolim X$.

\label\summary
\proclaim Lemma.
For any sequence~$X$,
$$\cst(X(-\infty)) \simeq \hofib(c_X \colon X\to\whX)$$
and
$$\hocofib(c_X \colon X\to\whX) \simeq \cst(\whX(\infty)/X(\infty)).$$

\proof Proof.
The arguments in \vxminusinf\ (and just around it) apply essentially verbatim.

The second verifies that taking associated graded determines a Quillen adjunction.

\label\Grasadjoint
\proclaim Lemma.
The functor $\Gr\colon\bSeq(\bA)\to\prod_\Z\bA$ is a left Quillen functor
whose right adjoint $R$ is given by
$$R(B)=(\cdots\To1{0}B_n\To1{0}B_{n+1}\To1{0}\cdots),$$
where $B$ is the graded object $(B_n)_{n\in\Z}$.

\proof Proof.
It is clear that $R$ preserves fibrations and acyclic fibrations, as those are objectwise,
so it suffices to verify that we have an adjunction.
Let $X$ denote a sequence $$\cdots\To4{x_{n-1}}X(n)\To4{x_n}X(n+1)\To4{x_{n+1}}\cdots$$
and $B$ a graded object $(B_n)_{n\in\Z}$.
To give a map of sequences $f\colon X\to R(B)$ is to give a map $f(n)\colon X(n)\to B_n$ for every~$n$ such that $f(n)\circ x_{n-1}=0\circ f(n-1)=0$.
Hence $f(n)$ factors as $$X(n)\To2{q_n}X(n)/\im x_{n-1}\To2{\bar f(n)}B_n.$$
Hence the map $f$ factors through $\bar f\colon\Gr A\to B$, the associated map between {\it graded\/} objects.
Conversely, any map of graded objects $g\colon\Gr X\to B$ produces a map of sequences by precomposition with the quotient maps $q_n\colon X(n)\to X(n)/\im x_{n-1}$.

We can now prove the goal of this section.

\label\localwe
\proclaim Proposition.
The class of $W_G$-local objects in $\bSeq(\bA)$ coincides with the class of derived complete sequences.
Furthermore, the class of weak equivalences in $\bFil(\bA)$ (i.e., local equivalences in $\bSeq(\bA)$) coincides with derived graded equivalences.

\proof Proof.
A~morphism of sequences $\cst(K) \to Y$ is given by a map $K \to \lim Y$ in $\bA$, and
similarly the derived mapping space $\Map(\cst(K),Y)$ between sequences is weakly equivalent to the derived mapping space $\Map(K,\holim Y)$ in $\bA$.
Observe as well that the morphism from $\cst(0) \to \cst(K)$ is always a derived graded equivalence.
Thus, a sequence~$Y$ is local with respect to the morphism $\cst(0) \to \cst(K)$ if the morphism of derived mapping spaces
$$\Map(K,\holim Y)\to \Map(0,\holim Y) \simeq {\rm pt}$$ is a weak equivalence (by Theorem~17.7.7 in Hirschhorn~[\ModCat]).
Consequently, $Y$ is local for morphisms $\cst(0) \to \cst(K)$ with {\it arbitrary\/} $K$ if and only if $\holim Y$ is weakly equivalent to~$0$,
and so, by \vsummary, if and only if $Y$ is derived complete.
\ppar
It remains to show that derived graded equivalences contain all local equivalences with respect to derived complete sequences.
Let~$f\colon X\to Y$ be local for every derived complete sequence,
and without loss of generality assume $f$ is a cofibration between cofibrant objects.
Let $Z$ be a fibrant object in $\bA$ and let $Z[n]$ denote the sequence that is zero except for $Z$ placed in degree~$n$;
$Z[n]$ can be see as the right adjoint $R$ applied to the graded object with $Z$ in spot $n$.
By construction $Z[n]$ is derived complete.
Using Proposition~17.4.16 in Hirschhorn~[\ModCat] and \vGrasadjoint, we see that the map $\Map(Y,Z[n]) \to \Map(X,Z[n])$ is weakly equivalent to $\Map(Y(n)/Y(n-1),Z)\to\Map(X(n)/X(n-1),Z)$,
which is a weak equivalence by the assumption that $f$ is local.
Since $Z$ is an arbitrary fibrant object in~$\bA$,
Theorem~17.7.7 in~[\ModCat] again shows that $X(n)/X(n-1) \to Y(n)/Y(n-1)$ is a weak equivalence.
The latter expression computes the (derived) associated graded of~$f$ at index $n$.
Varying $n$, we find that $f$~is a derived graded equivalence.
\ppar
We remark that the above proof shows that it suffices to perform left Bousfield localization with respect to morphisms of the form $\cst(0)\to \cst(K)$.

\proclaim Corollary.
The left derived functor of $\Gr\colon\bFil(\bA)\to\prod_\Z\bA$ creates weak equivalences.

\proof Proof.
The derived graded equivalences are, by definition, the maps created by the left derived functor of
$\Gr\colon\bSeq(\bA)\to\prod_\Z\bA$.
\vlocalwe\ shows that the class of weak equivalences of $\bFil(\bA)$ coincides with derived graded equivalences.

\subsection t-structures and spectral sequences

To construct a spectral sequence for a filtered object of a stable model category $\bA$, we need an analog of homology groups.
The standard tool is a t-structure, which provides an abelian category $\bA^\heartsuit$ and a nice functor $\pi_0 \colon \Ho(\bA) \to \bA^\heartsuit$.
The ``homology groups'' of an object $X \in \bA$ are then the collection $\{\pi_0(\Omega^n X)\}_{n \in \Z}$,
where $\Sigma$ denotes suspension and $\Sigma^{-1} = \Omega$ denotes looping.
Once one has a t-structure, there is a natural way to construct a spectral sequence for a filtered object,
as we explain below.

\label\tstructuresmodel
\subsubsection t-structures for stable model categories

In the context of triangulated categories and stable quasicategories, there is a well-known notion of t-structure.
We give a transcription into the language of stable model categories,
which is desirable if one wants to perform concrete computations with spectral sequences.

Recall (Definition~17.8.1 in Hirschhorn~[\ModCat]) that
a {\it homotopy orthogonal pair\/} is a pair of morphism $i\colon A\to B$ and $p\colon X\to Y$
such that the canonical map $$\Map(B,X)\to\Map(A,X)\times_{\Map(A.Y)}\Map(B,Y)$$ is a weak equivalence.
Here $\Map$ denotes the derived mapping space and $\times$ denotes homotopy pullback.
In this case we say that $i$~is {\it left homotopy orthogonal\/} to~$p$
and $p$~is {\it right homotopy orthogonal\/} to~$i$.

Recall (Bousfield, Definition~6.2 in~[\Fact]) that a {\it homotopy factorization system\/} on a model category~$\bA$ is a pair $(E,M)$,
where $E$ is a class contained in the cofibrations of~$\bA$ and
$M$ is a class contained in the fibrations of~$\bA$,
such that
\item{(1)} the class $E$ coincides with the class of cofibrations that are left homotopy orthogonal to all elements of~$M$;
\item{(2)} the class $M$ coincides with the class of fibrations that are right homotopy orthogonal to all elements of~$E$;
\item{(3)} any morphism in~$\bA$ factors as an element of~$E$ composed with an element of~$M$.
\endlist
We additionally demand that the latter factorization is functorial, in complete analogy to the modern definition of a model category.
A~homotopy factorization system is {\it accessible\/} if its functorial factorization is an accessible functor.

\label\factsat
\proclaim Remark.
According to Theorem~7.1 in Bousfield~[\Fact] (modified to use framing instead of simplicial enrichments
in the case of nonsimplicial model categories), any set~$I_E$ of cofibrations in a combinatorial model category
generates an accessible homotopy factorization system whose class~$E$ is the smallest weakly saturated class containing~$I_E$ and
all (generating) acyclic cofibrations and satisfying the condition that for any composable pair~$(f,g)$ of cofibrations such that $gf$ and~$f$ belong to $E$, then so does~$g$.
Typically $I_E$ will be a subset of the set of generating cofibrations.

The functorial factorization provides two natural functors on $\bA$.
Every object~$X\in\bA$ has a unique terminal morphism $X\to1$,
and we let $\tau_{<0}$ denote the functor arising from the factorization $X\to\tau_{<0}X\to1$.
Likewise, factoring $0\to X$ gives us $0\to\tau_{\ge0}X\to X$ for some functor~$\tau_{\ge0}$.
The notations $\tau_{<0}$ and $\tau_{\ge0}$ were chosen in anticipation of the definition below.
The functors $\tau_{<0}$ and $\tau_{\ge0}$ preserve weak equivalences and thus need not be derived.
The maps $\tau_{\ge0}\tau_{\ge0}X\to\tau_{\ge0}X$ and $\tau_{<0}X\to\tau_{<0}\tau_{<0}X$ are isomorphisms,
i.e., $\tau_{\ge0}$ is an idempotent monad and $\tau_{<0}$ an idempotent comonad.
In particular, they induce two adjunctions $\tau_{<0}\dashv\iota_{<0}$ and $\iota_{\ge0}\dashv\tau_{\ge0}$,
where $\iota_{<0}$ and $\iota_{\ge0}$ denote the inclusion functors
on the full subcategories of~$\bA$ consisting of the objects in the image of $\tau_{<0}$ respectively $\tau_{\ge0}$.
The unit of the first adjunction is $X\to\tau_{<0}X$ and the counit of the second adjunction is $\tau_{\ge0}X\to X$,
whereas the other (co)unit is the identity functor.

In terms of functors $\tau_{<0}$ and $\tau_{\ge0}$ one can characterize the classes $E$~and~$M$ as follows:
$E$~consists of those cofibrations~$f$ for which $\tau_{<0}f$ is a weak equivalence
and $M$~consists of those fibrations~$g$ for which $\tau_{\ge0}g$ is a weak equivalence.

\proclaim Definition.
A~{\it t-structure\/} on a combinatorial stable model category~$\bA$
is a homotopy factorization system on~$\bA$ that is {\it normal\/}:
for any object~$X$ the induced sequence $\tau_{\ge0}X\to X\to\tau_{<0}X$ is homotopy exact,
i.e., applying the derived mapping space functor~$\Map(-,Z)$ for any object~$Z$ gives us a homotopy fiber sequence of pointed spaces.
A~t-structure is {\it accessible\/} if its homotopy factorization system is.

We give a detailed comparison with the quasicategorical notion of t-structure in \vcomptstr.

\proclaim Remark.
Passing to the homotopy category of~$\bA$, the normality property reproduces property~(iii) of the original Definition~1.3.1
of a t-structure on a triangulated category from Beilinson, Bernstein, and Deligne~[\Perv].
The other two properties (i)~and~(ii) are supplied by the notion of a homotopy factorization system.

\proclaim Remark.
As shown by Lurie in Proposition~1.2.1.16 of~[\HA], the normality condition
can be reformulated as an additional condition on~$E$:
the class~$E$ is generated by the subclass consisting of morphisms whose domain is the initial object of~$\bA$.
Equivalently, one can say that if $A\to B$ is in~$E$, then so is its homotopy base change along $0\to B$,
as explained in Proposition~4.11 of Fiorenza and Loregi\`an~[\NTT].
Using the latter as an additional saturation condition to the one explained in \vfactsat,
we can speak of {\it t-saturation}.

\proclaim Remark.
Every stable \ii-category~$C$ has two trivial t-structures:
in one, all objects are in $C_{\ge0}$, and in the other, only the zero object belongs to $C_{\ge0}$.
These correspond to picking $I_E$ empty or everything.
We are typically interested in less extremal choices!

In light of the discussion above, we have the following result.

\proclaim Proposition.
For any combinatorial stable model category~$\bA$
and any set of cofibrations~$I_E$,
there is an accessible t-structure on~$\bA$
whose class~$E$ is generated by~$I_E$ under t-saturation.

\proclaim Example.
Consider the symmetric monoidal combinatorial model category~$\Ch$ of unbounded chain complexes of abelian groups equipped with the projective model structure.
A~standard set of generating cofibrations consists of the inclusions $\Z[n]\to(\Z[n+1]\to\Z[n])$.
Take $I_E$ to consist of those generating cofibrations with $n\ge0$, and construct the associated t-structure.
Maps in~$E$ are cofibrations that induce an isomorphism on homology in strictly negative degrees.
Maps in~$M$ are fibrations that induce an isomorphism on homology in nonnegative degrees.

\proclaim Example.
Consider the symmetric monoidal combinatorial stable model category~$\sSp$ of symmetric simplicial spectra, equipped with the stable projective model structure.
(See Hovey, Shipley, Smith~[\SymSpec].)
The generating cofibrations of~$\sSp$ are $\Sph\otimes(\Sigma_n\times(\sb^k\to\ss^k))_+[n]$,
where $\Sph$ denotes the symmetric sequence that gives the sphere spectrum.
Define $I_E$ to consist of those generating cofibrations for which $k\ge n$.
Observe that for this t-structure a morphism $0\to X$ belongs to~$E$ if the $n$th spectral level of~$X$ is {\it strictly $n$-connective}, i.e.,
is a pointed simplicial set whose set of $k$-simplices consists of one element for $k<n$.
In order to identify this t-structure with the standard t-structure on spectra
it suffices to show that the closure under weak equivalences of the codomains of morphisms in~$E$ of the form~$0\to X$
coincides with connective spectra.
In one direction, elements of~$I_E$ are morphisms of connective spectra, hence $\tau_{<0}$ applied to them produces an equivalence.
For the other direction it suffices to show that spectra~$X$ for which the terminal morphism $X\to1$
is right homotopy orthogonal to~$I_E$ are coconnected (meaning $\tau_{<0}$ vanishes).
(Abusing the language, we say that $X$ itself has this property.)
Without loss of generality assume $X$ to be fibrant.
Next, observe that in the above formula $\Sph\otimes-$ is derived, hence we can
use the adjunction property to reduce the problem to identifying symmetric sequences in pointed simplicial sets
that are right homotopy orthogonal to the maps $(\Sigma_n\times(\sb^k\to\ss^k))_+[n]$ for $k\ge n$.
Again, the functor $(\Sigma_n\times-)_+[n]$ is derived in the above formula, and using the same trick again
we reduce the problem to identifying morphisms of simplicial sets that are right homotopy orthogonal to the maps $\sb^k\to\ss^k$ for $k\ge n$.
These are precisely the $(n-2)$-truncated simplicial sets.
Thus we have an $\Omega$-spectrum whose $n$th spectral level is $(n-2)$-truncated (e.g., the 0th level is contractible),
hence its nonnegative homotopy groups vanish, hence it is coconnected.

\proclaim Remark.
Any t-structure on a stable model category gives rise to two different sequences for any object~$X$.
The first sequence is $$\cdots\to\tau_{\ge k+1}X\to\tau_{\ge k}X\to\tau_{\ge k-1}X\to\cdots,$$
which generalizes Whitehead towers of spectra to the setting of stable model categories.
The second sequence is $$\cdots\to\tau_{<k+1}X\to\tau_{<k}X\to\tau_{<k-1}X\to\cdots,$$
which generalizes Postnikov towers
and can be obtained by taking the homotopy cofiber of the canonical map of the Whitehead filtration on~$X$ to the constant sequence on~$X$.
(Thus in the setting of filtered objects, these two objects only differ by suspension.)

\subsubsection Spectral sequences

We now turn to setting up spectral sequences.
As shown by Beilinson, Bernstein, and Deligne in~[\Perv],
the t-structure on $\Ho(\bA)$ picks out a full subcategory of the homotopy category,
known as its {\it heart\/}, which is abelian.
We denote this category by $\bA^\heartsuit$.
By construction, the double truncation functor $\tau_{\le0} \circ \tau_{\ge0}$ on $\Ho(\bA)$ maps to $\bA^\heartsuit$.
We denote by $\pi_0$ the composition $\bA \to \Ho(\bA) \to \bA^\heartsuit$.
Thus, we obtain ``homotopy groups'' with values in $\bA^\heartsuit$ for every object in $\bA$
by the formula $\pi_n = \pi_0 \circ \Sigma^{-n}$.

Each sequence $X \in \bSeq(\bA)$ has an associated spectral sequence taking values in $\bA^\heartsuit$,
which we quickly describe below, following \S1.2.2 of Lurie~[\HA].
We now obtain a generalization of the classical statement that a filtered complex and its completion have the same spectral sequence.

\label\mcspseq
\proclaim Proposition.
A~sequence $X$ and its derived completion $\whcrfX$ have the same spectral sequence.
Hence, the functor from $\bSeq(\bA)$ to the category of spectral sequences in $\bA^\heartsuit$
factors through the filtered category $\bFil(\bA)$.

This proposition follows immediately from the construction of the spectral sequence $E^{*,*}_*(X)$ of a sequence~$X$, which is given by
$$E_r^{p,q}(X) = \im\pi_{p+q}(\hocofib(X(p-r) \to X(p)) \to \hocofib(X(p-1) \to X(p+r-1))).$$
The differential $d_r$ is determined by the connecting maps in the long exact sequences
of homotopy groups arising from the homotopy pushout square
$$\cd{\hocofib(X(i) \to X(j)) &\mapright{} & \hocofib(X(i) \to X(k))\cr
\mapdown{} &&\mapdown{}\cr
\llap{$0={}\!$}\hocofib(X(j) \to X(j))& \mapright{} & \hocofib(X(j) \to X(k))\cr}$$
that holds for any triple $i\le j\le k$.
(For a more extensive discussion, see Construction~1.2.2.6 and Proposition~1.2.2.7 of Lurie~[\HA].)

\proof Proof.
The canonical map $q_X \colon \crf X \to X$ is a weak equivalence,
and hence this map induces an isomorphism of spectral sequences.
By \vcomplgreq, we know $\crf X \to \crf(\whcrfX)$ is a graded equivalence,
and so \vfinitegraded\ tells us the map provides an isomorphism of spectral sequences.
Finally, $q_{\whcrfX}\colon \crf(\whcrfX) \to \whcrfX$ again induces
an isomorphism of spectral sequences.
Altogether, we obtain the proposition.

\subsection Enrichments and enhancements

The following result shows that an enrichment on the underlying model category $\bA$ passes to its category of filtered objects $\bFil(\bA)$.
In particular, if the original model category is dg enhanced (i.e., enriched over chain complexes),
then so is its category of filtered objects.
Thus this result provides one more sense in which to interpret the paper's title.

Recall from Barwick~[\LR] that a model category is {\it tractable\/}
if it is combinatorial and all cofibrations are generated, under weak saturation, by cofibrations with cofibrant source.
Examples of tractable model categories include
chain complexes over a ring (or differential graded algebra), simplicial symmetric spectra, motivic symmetric spectra,
presheaves of simplicial symmetric spectra with the local model structure,
and symmetric spectra in some tractable left proper symmetric monoidal model category.
Indeed, almost every combinatorial model category encountered in practice is tractable;
constructing a nontractable combinatorial model category is nontrivial,
but see the example after Corollary~8.6.4 in Simpson~[\HTHC].

\proclaim Proposition.
If $\bV$ is a symmetric monoidal tractable model category
and $\bA$ is a $\bV$-enriched model category
that is tractable, left proper, and stable,
then the model category $\bFil(\bA)$ is $\bV$-enriched and tractable.

\proof Proof.
The model category $\bSeq(\bA)$ is $\bV$-enriched
because the pushout-product axiom and the unit axiom can be verified indexwise.
We now apply Theorem~4.46 in Barwick~[\LR],
which requires us to verify (as explained in the proof there) that the class of maps
with respect to which we localize (in our case, the derived graded equivalences)
is closed with respect to the derived tensor product with an arbitrary object of~$\bV$.
Observe that derived graded equivalences are defined using a homotopy cofiber,
and homotopy cofibers commute with derived tensor products.
(Derived tensor product with a fixed object $X\otimes^L- $ can be computed as the tensor product with a fixed cofibrant replacement of this object $\crf(X) \otimes -$.
As the left derived functor of the left Quillen functor, it preserves homotopy colimits.)

\proclaim Remark.
A~special case appears frequently in contemporary mathematics:
$\bV=\Ch(R)$, where $R$ is a commutative ring.
In this case $\bV$-enriched model categories are also known as dg model categories.
Given such a model category $\bC$, one can extract the dg category $\bC^\circ$ of bifibrant objects,
which as a dg category has the same homotopical information as the model category~$\bC$:
its dg nerve is equivalent to $\cU\bC$ as a quasicategory.
If $\bC$ satisfies the conditions of the above proposition,
then the model category $\bFil(\bC)$ is a stable tractable dg model category.
In \S\dgcat\ we also construct a filtered dg category associated to any dg category
using methods intrinsic to dg categories.
In \vdgmodelbig\ of~\vdgcat,
we explain how the enriched model-categorical approach relates to that dg categorical construction.

\subsection Symmetric monoidal structures

Many of the most important stable model categories in mathematics have a natural symmetric monoidal structure.
For example, given a commutative ring $R$, chain complexes of $R$-modules can be tensored over~$R$.
We now explain how to incorporate such a structure into the filtered setting.

Let $(\bA,\otimes)$ now denote a symmetric monoidal combinatorial stable model category.
We do not assume the unit to be cofibrant, but we do require Hovey's unit axiom.
Our main result in this section is the following.

\label\modotimes
\proclaim Theorem.
There is a symmetric monoidal product $\otimes_s$ on $\bFil(\bA)$ such that $(\bFil(\bA),\otimes_s)$ is a symmetric monoidal stable model category.

The rest of the section is devoted to the proof.
We begin by equipping $\bSeq(\bA)$ with a symmetric monoidal structure (via Day convolution) and verify it satisfies the pushout-product axiom.
We then define~$\otimes_s$ and prove the main result.

The arguments depend on verifying certain properties of generating cofibrations,
so we need to identify these maps in $\bSeq(\bA) = \Fun(\Zs,\bA)$ and $\Fun(\Zs \times \Zs, \bA)$,
which are both equipped with the projective model structure.
Every $n \in \Zs$ and every $A \in \bA$ determines a ``step-sequence,'' denoted $\langle n,A \rangle$,
that sends~$m$ to~$0$ if $m < n$ and to~$A$ if $m\ge n$ (with identity as structure map).
A~generating cofibration in $\bSeq(\bA)$ is then
a map $f\colon \langle n, A\rangle \to \langle n,A'\rangle$ where for $m\ge n$, the map is always a fixed cofibration $\phi \colon A \to A'$ in $\bA$.
A~generating cofibration in $\Fun(\Zs \times \Zs,\bA)$ is entirely analogous:
a map $F\colon \langle (m,n),A \rangle \to \langle (m,n), A' \rangle$ where the map is a cofibration $\phi \colon A \to A'$ for $(i,j)$ with $i\ge m$ and $j\ge n$, but is otherwise zero.

Now we define the ``external product'', which is the bifunctor
$$\boxtimes\colon\Fun(\Zs,\bA) \times \Fun(\Zs,\bA) \to \Fun(\Zs \times \Zs, \bA)$$
sending $(X,Y)$ to $X \boxtimes Y\colon (m,n) \mapsto X(m) \otimes Y(n)$.
It possesses the following crucial property.

\label\pushoutproduct
\proclaim Lemma.
Let $f \colon X \to Y$ and $g \colon U \to V$ be cofibrations in $\bSeq(\bA)$.
Then the natural map
$$X \boxtimes V \sqcup_{X \boxtimes U} Y \boxtimes U \to Y \boxtimes V$$
is a cofibration in $\Fun(\Zs \times \Zs,\bA)$.
If one of the cofibrations is acyclic, then so is the resulting map.

\proof Proof.
It suffices to check this assertion when $f$ and $g$ are generating cofibrations.
Let $X = \langle m,A\rangle$ and $Y = \langle m,B\rangle$ with $f$ determined
by a cofibration $\phi\colon A \to B$ in~$\bA$, and likewise $U = \langle n,C\rangle$ and $V = \langle n,D\rangle$ with $g$ determined
by a cofibration $\psi\colon C \to D$ in~$\bA$.
Then, for instance, $X \boxtimes U = \langle(m,n),A\otimes C\rangle$.
The pushout $X \boxtimes V \sqcup_{X \boxtimes U} Y \boxtimes U$ is then
$$\langle (m,n), A \otimes D \sqcup_{A \otimes C} B \otimes C \rangle.$$
The natural map from the pushout to the product is then the generating cofibration determined by $(m,n)$ and the cofibration
$$A \otimes D \sqcup_{A \otimes C} B \otimes C \to B \otimes D$$
in~$\bA$.

Next, we define the ``totalization'', which is the functor
$$\tot\colon \Fun(\Zs \times \Zs,\bA) \to \Fun(\Zs,\bA) = \bSeq(\bA)$$
sending $Z$ to $\tot(Z)(n) = \colim_{p+q\le n} Z(p,q)$.
(This colimit is over the full subcategory of $\Zs^2$ given by objects $(p,q)$ such that $p+q\le n$.)

\proclaim Lemma.
$\boxtimes$ is a left Quillen bifunctor, and $\tot$ is a left Quillen functor.

\proof Proof.
Now observe that $\tot$ is the left Kan extension along the addition map $+\colon \Zs \times \Zs \to \Zs$.
By Proposition~A.2.8.7 in Lurie~[\HTT], $\tot$ is thus a left Quillen functor with respect to the projective model structures.
(We note that this proposition is slightly more general than stated:
it holds for a cofibrantly-generated model category, not just a combinatorial model category.)

\proclaim Definition.
For $X$ and $Y$ in $\bSeq(\bA)$, we define $X \otimes_s Y = \tot(X \boxtimes Y)$.

\proclaim Proposition.
The composition $\otimes_s = \tot \circ\; \boxtimes$ makes $\bSeq(\bA)$ into a symmetric monoidal stable model category.

\proof Proof.
Applying $\tot$ to the map in \vpushoutproduct, we find that
$$A \otimes_s D \sqcup_{A \otimes_s C} B \otimes_s C \to B \otimes_s D$$
is a cofibration in $\bSeq(\bA)$ with the projective model structure.
Hence, $\otimes_s$ satisfies the pushout-product axiom with respect to the projective model structure.

\proof Proof of \vmodotimes.
We now want to verify that we obtain a monoidal model category after left Bousfield localizing at the graded equivalences $W_G$.
By Proposition~4.47 of Barwick~[\LR], it suffices to show that the derived graded equivalences are closed
under taking the derived tensor product with an arbitrary object.
In fact, it suffices to verify this only for a collection of weak ({\it aka\/} homotopy) generators.
We will work with the sequences $\langle m, A \rangle$, which are zero for $n < m$ and a cofibrant object $A$
from $\bA$ for $n\ge m$, with the identity as the structure map; these are homotopy generators.
\ppar
Now consider $f \colon X \to Y$ a graded equivalence between cofibrant sequences.
We want to show that
$$\langle m, A \rangle \otimes_s f\colon \langle m, A \rangle \otimes_s X \to \langle m, A \rangle \otimes_s Y$$
is a graded equivalence of cofibrant sequences.
Observe that
$$\left(\langle m, A \rangle \otimes_s X \right)(n)/\left(\langle m, A \rangle \otimes_s X \right)(n-1)
\cong A \otimes \left(X(-m+n)/X(-m+n-1)\right).$$
Thus, as $A$ is cofibrant, the induced map $\Gr \left(\langle m, A \rangle \otimes_s f \right)$ is an indexwise weak equivalence.
\ppar
We now prove the unit axiom.
The monoidal unit is given by $\langle0,1\rangle$ and its cofibrant replacement can be taken to be $\langle0,\crf1\rangle$,
where $\crf1\to1$ is a cofibrant replacement in~$\bA$.
It now suffices to observe that $\langle0,1\rangle\otimes X=X$ and $\langle0,\crf1\rangle\otimes X=\crf1\otimes X$ (here $\crf1\otimes X$ is taken componentwise),
and that a cofibrant object~$X$ in particular has cofibrant components.
Hence the map $\crf1\otimes X\to1\otimes X=X$ is a weak equivalence by the unit axiom of~$\bA$.

\label\modelgrmonoidal
\proclaim Proposition.
The associated graded functor has canonical strong monoidal and strong closed structures.

\proof Proof.
The argument for \vstrongmonoidal\ and \vstrongclosed\ works here with no changes.
The associated graded functor is cocontinuous, and the tensor product is strong monoidal on generators.

\subsection The monoid axiom

Suppose $\bA$ is a symmetric monoidal combinatorial stable model category and $R$ is a monoid in filtered objects in~$\bA$.
We would like to equip the category $\Mod_R$ of left $R$-modules in $(\bFil(\bA),\otimes_s)$ with a model structure
whose fibrations and weak equivalences are transferred along the forgetful to~$\bFil(\bA)$.
Schwede and Shipley~[\MonAx] have provided a simple and conceptual condition for this kind of construction to work:
Theorem~4.1 of~[\MonAx] states that for a cofibrantly generated monoidal model category~$\bC$,
the category of left modules over any monoid in~$\bC$ inherits a cofibrantly generated model structure
if $\bC$ satisfies the {\it monoid axiom\/} (Definition~3.3 in~[\MonAx]).

As explained in~\S3.2 of Pavlov and Scholbach~[\HTSP], it is convenient to split the monoid axiom into simpler properties that are easier to establish separately.
We recall relevant definitions and results from that paper.
An {\it h-cofibration\/} is a map $X\to Y$ such that the pushout functor $X/\bA\to Y/\bA$ on undercategories preserves weak equivalences.
For the case of left proper model categories considered below, h-cofibrations can also be characterized
as those maps for which cobase changes (pushouts) are also homotopy cobase changes.
In these terms, left properness then demands that cofibrations are h-cofibrations.
A~monoidal model category is {\it h-monoidal\/} if the monoidal product of any object with a cofibration is an h-cofibration.
Lemma~3.2.6 of~[\HTSP] implies that a h-monoidal, compactly generated model category satisfies the monoid axiom.
A~monoidal model category is {\it flat\/} if the pushout product of a cofibration and a weak equivalence is a weak equivalence.
The general idea behind these definitions is that h-monoidality ensures the existence of a model structure,
whereas flatness gives rise to Quillen equivalences between model structures constructed using weakly equivalent data.

Flatness and h-monoidality have symmetric analogs, as explained in Definition~4.2.2 and Definition~4.2.7 of~[\HTSP].
These properties play the same role for symmetric structures, such as commutative monoids and algebras over symmetric operads,
as their nonsymmetric analogs for nonsymmetric structures, such as left modules over monoids and algebras over nonsymmetric operads.

\label\filtprops
\proclaim Proposition.
If $\bA$ is a (symmetric) flat, (symmetric) h-monoidal, symmetric monoidal, left proper, compactly generated, combinatorial stable model category, then so are $\bSeq(\bA)$ and $\bFil(\bA)$.

\proof Proof.
We have already established that $\bSeq(\bA)$ exists and is a symmetric monoidal combinatorial stable model category.
It is compactly generated because cofibrations do not change under left Bousfield localizations,
whereas weak equivalences increase.
It is left proper because projective cofibrations are in particular componentwise cofibrations.
Flatness and h-monoidality of~$\bSeq(\bA)$ follow from the same properties of~$\bA$
because they can be established for generating (acyclic) cofibrations (see Theorem~3.2.8 in~[\HTSP] in the nonsymmetric case
and Theorem~4.3.9 there in the symmetric case),
and the latter are monoidal products of representables and generating (acyclic) cofibrations of~$\bA$.
Finally, by Proposition~6.2.1 of~[\HTSP] the same set of properties holds for the left Bousfield localization $\bFil(\bA)$ of the projective structure on $\bSeq(\bA)$.

\proclaim Examples.
Simplicial sets, simplicial modules, chain complexes of modules, simplicial presheaves all satisfy the conditions of the proposition,
as explained in~\S7 of~[\HTSP].
Topological spaces can also be used if one replaces compact generatedness with strong admissible generatedness, as explained there.

The arguments of Schwede and Shipley in~[\MonAx] then imply the following.

\label\modules
\proclaim Corollary.
Let $\bA$ satisfy the hypotheses in the preceding proposition.
If $R$ is a monoid in $\bFil(\bA)$, then the category $\Mod_R$ of left $R$-modules in $(\bFil(\bA),\otimes_s)$
inherits a model category structure via the forgetful functor $\bFil(\bA)$:
weak equivalences and fibrations are checked in $\bFil(\bA)$.

There is a modest generalization.
A~symmetric monoidal category $(\bC,\otimes,1)$ is a commutative monoid
in the (large) bicategory of categories, functors, and natural transformations.
A~{\it left $\bC$-module\/} $\bD$ is a left module over~$\bC$ in this bicategory.
Given a monoid~$R$ in~$\bC$, one can then study left modules in~$\bD$ over~$R$.
Lemma~3.2.6 of~[\HTSP] gives conditions under which the monoid axiom holds in such a left module category,
and hence when the category of left $R$-modules in~$\bD$ admits a compatible model structure.

\proclaim Corollary.
Let $\bA$ be a left module over a monoidal category~$\bC$
such that $\bA$ is a left proper compactly generated combinatorial stable model category.
Let~$R$ be a monoid in~$\bFil(\bC)$, and
let~$\fMod_R(\bA)$ denote the category of left $R$-modules in~$\bFil(\bA)$.
If for all $n$ the functor $R_n\otimes-\colon\bA\to\bA$ sends acyclic cofibrations to acyclic h-cofibrations,
then category~$\fMod_R(\bA)$ admits a model structure
whose weak equivalences and fibrations are inherited via the forgetful functor to $\bFil(\bA)$.

\label\seccomparison
\section Comparison of the quasicategorical and the model categorical constructions

In this section, we compare our constructions in the quasicategorical and model categorical settings.
We start by showing that every stable model category has an ``underlying'' stable quasicategory,
and then show that this construction intertwines with forming the filtered category.
In consequence, we deduce the folk theorem, \vthmstmodcat.
We show next that, conversely, every presentable stable quasicategory lifts to the stable model category setting.
Finally, we verify the compatibility of t-structures in both settings.

\subsection Comparing versions of filtered categories

Let us begin by describing the ``underlying quasicategory'' of a combinatorial model category, following Lurie's treatment in \S1.3.4 of~[\HA].
Recall Definition 1.3.4.15 of~[\HA]: the {\it underlying quasicategory\/} $\cU(\bA)$ of a model category $\bA$
is any quasicategory~$\cC$ equipped with a functor $\N(A^c)\to\cC$ that realizes~$\cC$ as $\N(A^c)[W^{-1}]$,
i.e., a quasicategory obtained from the nerve $\N(\bA^c)$ of the category~$\bA^c$ of cofibrant objects in~$\bA$ by inverting the weak equivalences $W$ up to homotopy.
(For example, take the marked simplicial set $(\N(A^c),W^c)$ and fibrantly replace.)
Lurie produces such an underlying quasicategory of a combinatorial model category in three steps.
First, by Theorem~1.1 in Dugger~[\Pres], every combinatorial model category~$\bA$ is Quillen equivalent
to a simplicial left proper combinatorial model category~$\hat\bA$.
Second, Lemma~1.3.4.21 in~[\HA] shows that a left Quillen equivalence~$\hat\bA\to\bA$ between combinatorial model categories
induces a weak equivalence $\N(\hat\bA^c)[W_{\hat\bA}^{-1}]\to\N(\bA^c)[W_\bA^{-1}]$ of the corresponding marked simplicial sets, which present inverting morphisms up to homotopy.
Third, Theorem~1.3.4.20 in~[\HA] constructs an equivalence $\hat\bA^c[W^{-1}]\to\N(\hat\bA^\circ)$ of quasicategories
for any simplicial model category~$\hat\bA$, where $\hat\bA^\circ$ is the full simplicial subcategory of bifibrant objects in~$\hat\bA$.
Finally, Proposition~A.3.7.6 in~[\HTT] shows that $\N(\hat\bA^\circ)$ is presentable for any combinatorial simplicial model category~$\hat\bA$.

\label\modtoqc
\proclaim Lemma.
Given a combinatorial stable model category $\bA$, the underlying quasicategory $\cU(\bA)$ is a presentable stable quasicategory.

\proof Proof.
Presentability was established above.
Stability amounts to being pointed and the suspension functor being an equivalence of quasicategories.
By Theorem~4.2.4.1 in~[\HTT], the homotopy coherent nerve functor sends homotopy (co)limit diagrams in~$\hat\bA^\circ$
to (co)limit diagrams in~$\N(\hat\bA^\circ)$.
Thus bicartesian squares are mapped to bicartesian squares, and similarly for the zero object.

We now verify that our constructions of the filtered category intertwine via $\cU$.

\label\bFilvsiFil
\proclaim Proposition.
For any left proper combinatorial stable model category~$\bA$,
the canonical functor of quasicategories $$\cU(\bFil(\bA))\to\iFil(\cU(\bA))$$ is an equivalence,
where the functor~$\cU$ takes the underlying quasicategory of a model category.
Furthermore, the above equivalence of categories is compatible with the associated graded functor:
the square $$\cd{\cU(\bFil(\bA))&\mapright{}&\iFil(\cU(\bA))\cr\mapdown{\cU(\Gr_\bA)}&&\mapdown{\Gr_{\cU(\bA)}}\cr \cU(\prod_\Z\bA)&\mapright{}&\prod_\Z \cU(\bA)\cr}$$
commutes up to a homotopy.

\proof Proof.
Proposition~1.3.4.25 in Lurie~[\HA] establishes the result for sequences: the canonical functor
$$\cU(\bSeq(\bA))\to\iSeq(\cU(\bA))$$
is an equivalence of quasicategories.
This proposition also implies that the bottom arrow of the square is an equivalence.
Proposition~A.3.7.8 in Lurie~[\HTT] verifies that the localizations intertwine: it thus promotes the equivalence for sequences to the desired equivalence of completed sequences.
\ppar
It remains to verify the commutativity of the square.
Using the description of the associated graded functor after \videfgraded\ (and its obvious model categorical counterpart),
we present the square as the vertical composition of two squares,
where the middle row is
$$\cU\left(\prod_\Z\Fun(0\to1,\bA)\right)\to\prod_\Z\Fun(0\to1,\cU(\bA)),$$
which is an equivalence of quasicategories by Proposition~1.3.4.25 in~[\HA].
The vertical maps in the top square are the restriction functors along $\coprod_\Z(0\to1)\to\Zs$, where the $m$th index goes to the arrow $m \to m+1$ in~$\Zs$.
The vertical maps in the bottom square are given by taking cofibers;
they are the left adjoints of the left adjoints of the restriction functors along $\Z\to\coprod_\Z(0\to1)$.
(In other words, we are just decomposing the associated graded functors in the natural way.)
\ppar
It now suffices to show the commutativity of the top and bottom square separately.
These claims follow once we show the general fact that for any functor $I\to J$, the square
$$\cd{\cU(\Fun(J,\bA))&\mapright{}&\Fun(\N(J),\cU(\bA))\cr
\mapdown{}&&\mapdown{}\cr
\cU(\Fun(I,\bA))&\mapright{}&\Fun(\N(I),\cU(\bA))\cr}$$
is commutative up to a homotopy,
where $\Fun(-,\bA)$ is equipped with the {\it injective\/} structure so that the functor $\Fun(J,\bA)\to\Fun(I,\bA)$ is a left Quillen functor.
The horizontal arrows are equivalences, by the arguments above, and the vertical arrows are cocontinuous functors.
Observe that the following diagram commutes strictly:
$$\cd{\N(\Fun(J,\bA)^c)&\mapright{}&\Fun(\N(J),\cU(\bA))\cr
\mapdown{}&&\mapdown{}\cr
\N(\Fun(I,\bA)^c)&\mapright{}&\Fun(\N(I),\cU(\bA)),\cr}$$
where the superscript~$c$ denotes the class of cofibrant objects.
Therefore, the diagram obtained by localizing the left column commutes up to a homotopy
because the space of factorizations of a functor~$F\colon C\to D$ through a localization functor $C\to C[W^{-1}]$ (of quasicategories) is either empty or contractible
(i.e., either~$F$ inverts the elements of~$W$ or it does not, and in our case it does).

\label\thmstmodcat
\proclaim Theorem.
For $\abcat$ a Grothendieck abelian category, let $\bA = \Ch(\abcat)$ be the stable model category given by the injective model structure.
The homotopy category of the filtered stable model category $\bFil(\bA)$ is equivalent to the classical filtered derived category $\Dfil(\abcat)$.

\proclaim Remark.
In~[\HA] Lurie constructs the derived \ii-category $\cD(\abcat)$ by taking the dg nerve of the full dg subcategory of
fibrant-cofibrant objects of $\bA$ with the injective model structure.
Our result here thus implies \vmainresult, the \ii-categorical assertion we made earlier.

\proof Proof.
Let us begin with $\bSeq(\bA)$ equipped with the projective model structure.
Consider the subcategory of classically filtered objects in $\bSeq(\bA)$: the sequences of chain complexes in which each structure map is a monomorphism.
The subcategory of cofibrant objects $\bSeq(\bA)^c$ is also a subcategory of these classically filtered objects.
The cofibrant replacement functor on $\bSeq(\bA)$, restricted to the classically filtered objects, exhibits a natural weak equivalence between these categories with weak equivalences.
In particular, their homotopy categories are equivalent.
\ppar
Left Bousfield localization does not change the cofibrant objects and so one can work with the same cofibrant replacement functor.
Thus, even after localization of $\bSeq(\bA)$ at the graded weak equivalences,
the cofibrant replacement functor produces a weak equivalence between sequences and cofibrant sequences,
and hence a weak equivalence between the subcategory of classically filtered objects and the cofibrant sequences.
(Using the intermediate model structure constructed in \vcofobj, the cofibrant objects are precisely the filtered chain complexes in the usual sense.
Thus it is manifest that the homotopy categories are equivalent, because this model structure has the same weak equivalences as $\bFil(\bA)$.)

We now show that the stable quasicategory setting can be lifted to the stable model category setting.

\proclaim Proposition.
Every presentable stable \ii-category is presented by a simplicial left proper combinatorial stable model category~$\bA$.

\proof Proof.
By Proposition~1.4.4.9.(3) in Lurie~[\HA], we observe that a presentable stable \ii-category is presented
by an accessible left exact localization of the \ii-category of presheaves of spectra on a small \ii-category.
The latter small \ii-category can be rectified to a small simplicial category,
and the above localization can be presented by a left exact accessible left Bousfield localization of the category of presheaves of simplicial spectra.
The resulting model category is simplicial, left proper, combinatorial, and stable.
(Simplicial spectra are stable and stability is preserved under taking presheaf categories and left exact left Bousfield localizations.)

\label\comptstr
\subsection Comparing t-structures

\vmodtoqc\ showed that a combinatorial stable model category has an associated presentable stable quasicategory.
We now that this construction respects t-structures.

\proclaim Proposition.
An accessible t-structure on a combinatorial stable model category $\bA$ produces an accessible t-structure on the quasicategory $\cU(\bA)$.

\proclaim Remark.
The notion of an accessible t-structure is Definition~1.4.4.12 in Lurie~[\HA].
Our definition of accessible t-structure was motivated to ensure the compatibility encoded in this proposition.

\proof Proof.
The proof is an assembly of several results.
\item{(1)} A~homotopy factorization system on $\bA$ (Bousfield, Definition~6.2 in~[\Fact])
induces a factorization system on $\cU(\bA)$ (Lurie, Definition~5.2.8.8 in~[\HTT], or Fiorenza and Loregi\`an, Definition~8 in~[\NTT])
by taking the image of the elements of $E$ and $M$ with bifibrant domains and codomains.
\item{(2)} The normality property of a homotopy factorization system descends to the quasicategorical normality property as in Definition~17 in Fiorenza and Loregi\`an~[\NTT],
which follows from the fact that the comparison functor maps homotopy exact sequences to exact quasicategorical sequences, see Propositions 1.3.4.23 and 1.3.4.24 in~[\HTT].
\item{(3)} The accessibility property of a homotopy factorization system descends to the small generation property, as in Remark~5.5.5.2 in Lurie~[\HTT].
Indeed, we can take the image under the comparison functor of some generating set for~$E$, which we can assume to have bifibrant (co)domains, without loss of generality.
Proposition~1.3.4.24 in~[\HTT] then implies that the image is still a generating set for the quasicategorical factorization system.
\ppar
\noindent In consequence, we can invoke Theorem~2 in Fiorenza and Loregi\`an~[\NTT]
to see that we obtain a t-structure on the triangulated category $\Ho(\bA)$ given by the homotopy category of~$\bA$.

Conversely, we can rigidify quasicategorical accessible t-structures to our situation.

\proclaim Proposition.
Any presentable stable quasicategory $\cC$ with an accessible t-structure
can be presented by a combinatorial stable model category with an accessible t-structure.

\proof Proof.
Let $\cC_{\ge0}$ denote the full subquasicategory of $\cC$ consisting of all objects~$X$ with $\tau_{<0}X\simeq0$.
The inclusion $\iota\colon\cC_{\ge0}\to\cC$ of quasicategories is a cocontinuous functor between presentable quasicategories.
Choose a cardinal~$\kappa$ such that both $\cC_{\ge0}$ and $\cC$ are locally $\kappa$-presentable and the functor $\tau_{\ge0}\colon\cC\to\cC_{\ge0}$ is $\kappa$-accessible.
Denote by $S$, respectively~$T$, the full subquasicategory of $\kappa$-compact objects in $\cC_{\ge0}$, respectively~$\cC$.
By construction, $S$ and~$T$ are essentially small quasicategories and the inclusion $\iota\colon\cC_{\ge0}\to\cC$ sends $S$ to~$T$.
The functor~$\tau_{\ge0}$ can now be presented as the composition $\cC\to\P(S)\to\cC_{\ge0}$,
where the first functor is the restricted Yoneda embedding
and the second functor is induced by the inclusion $S\to\cC_{\ge0}$ using the universal cocompletion property of~$\P(S)$.
\ppar
Recall Theorem~5.5.1.1.(5) in Lurie~[\HTT]:
a $\kappa$-presentable quasicategory can be presented as an accessible localization of the quasicategory of presheaves
on the (essentially small) full subquasicategory of $\kappa$-compact objects.
We have the following commutative square of quasicategorical right adjoint functors:
$$\cd{
\cC &\mapright{\tau_{\ge0}}& \cC_{\ge0} \cr
\mapdown{} && \mapdown{}\cr
\P(T) &\mapright{\iota^*}& \P(S),\cr}$$
where the vertical functors are the (restricted) Yoneda embeddings into presheaf quasicategories
and the bottom functor is restriction along the inclusion $\iota\colon S\to T$.
\ppar
Using the left adjoint of the homotopy coherent nerve functor, rigidify~$\iota$ to a fully faithful inclusion $\bar\iota\colon\bar S\to\bar T$ of simplicial categories.
(First rectify the codomain of~$\iota$, then take the full subcategory consisting of objects weakly equivalent to an object in the domain of~$\iota$.)
Take the induced adjunction of combinatorial model categories of simplicial presheaves
and perform a left Bousfield localization with respect to the maps that are inverted by the corresponding quasicategorical localizations.
(As usual, it suffices to take a generating set of such maps, which exists by accessibility.)
\ppar
Now use the above theorem to construct a t-structure generated by
a set of cofibrations~$I_E$ obtained by tensoring a generating cofibration of simplicial sets and the representable functor of an object in the image of~$\bar\iota$.
The domain of~$\iota$ corresponds to $\cC_{\ge0}$, so for any element~$f$ of the resulting class~$E$ the morphism~$\tau_{<0}f$ is an equivalence.
\ppar
Vice versa, any morphism~$g$ with $\tau_{<0}f$ an equivalence is weakly equivalent to a cofibration in~$E$.
Indeed, the homotopy factorization system constructed above factors~$g$ as an element~$q$ of the saturation (in the sense explained above) of the set~$I_E$
(hence a cofibration in~$E$),
and a morphism~$h$ that is right homotopy orthogonal to~$I_E$.
As explained above, $\tau_{<0}q$ is an equivalence, hence so is $\tau_{<0}h$.
Using the definition of~$I_E$, we immediately show that $h$ is a weak equivalence when evaluated on any element of~$\bar S$,
and by definition of~$S$ this means that $\tau_{\ge0}h$ is an equivalence.
Thus $\tau_{<0}h$ and $\tau_{\ge0}h$ are equivalences, which means that $h$ is an equivalence
and the morphism~$g$ can be presented as~$q$, which has the form we need.

\label\dgcat
\section Filtered objects in the language of differential graded categories

Differential graded (dg) categories---i.e., categories enriched over the category of chain complexes over a ring---are a particularly important class of higher categories
due to their wide use.
Here we will describe a construction of the filtered objects in a pretriangulated, idempotent-complete dg category and verify that it agrees with the constructions in the quasicategorical and model-categorical settings.
We do not discuss, however, spectral sequences, symmetric monoidal structures, or dualizability in this dg setting.
Indeed, our comparison result here relies on embedding the dg construction into a model categorical context,
and hence it seems natural to refer to that context for further results.

\subsection First definitions

Throughout this section we fix a commutative ring~$R$.
(Most constructions go through with $R$ a commutative dg ring.)

\proclaim Definition.
A~{\it differential graded category over a ring~$R$\/} is a category enriched over the symmetric monoidal category of (unbounded) chain complexes over~$R$
equipped with the usual tensor product.
A~{\it dg functor\/} between dg categories is an enriched functor.
The {\it underlying category\/} $\udg(\tC)$ of a dg category~$\tC$ has the same objects
and 0-cycles in the hom-complexes as morphisms.

The motivating example of a dg category is $\Ch(R)_\dg$,
where one enriches the ordinary category $\Ch(R)$ of chain complexes of $R$-modules and chain maps over itself.
That is, given two objects $M, N \in \Ch(R)$, the morphisms from $M$ to $N$ in $\Ch_\dg(R)$ form a chain complex of $R$-modules whose zero cycles are the chain maps from $M$ to $N$.
Another class of examples is given by sheaves of chain complexes on a space,
which can be naturally enriched to a dg category.
Such examples arise in topology (e.g., constructible sheaves) and geometry (e.g., dg enhancements of~$\D({\rm Coh}(X))$).

There is a category of dg categories $\dgCat_R$, and many familiar constructions from category theory admit analogs.
For instance, there is a common generalization of presheaves in ordinary category theory and modules in algebra.
A~{\it dg presheaf\/} or {\it module\/} on a dg category~$\tC$ is a dg functor $\tC^\op\to\Ch_\dg(R)$.
The collection of all modules on $\tC$ naturally forms a dg category that we denote~$\tC\dashdgMod$.

\subsection Homotopy limits and a special class of dg categories

We will need homotopy limits in a dg category,
in particular to talk about complete filtered objects.
There are several possible approaches,
but we will give the bare minimum needed to realize the construction of the filtered dg category.
Our approach exploits the enriched Yoneda embedding for dg categories, which preserves limits,
allowing us to use model category techniques to talk about homotopy limits of dg modules.

\label\dgcmod
\proclaim Definition.
Given a small dg category~$\tC$, let $\tC\dashdgMod$ denote the dg category of dg functors from~$\tC^\op$ to~$\Ch_\dg(R)$.
We equip it with objectwise weak equivalences and fibrations,
so that a natural transformation is a weak equivalence if on each object in $\tC$,
it determines a quasi-isomorphism or fibration in $\Ch(R)$.
Thus $\tC\dashdgMod$ becomes an enriched model category over~$\Ch(R)$ by Theorem~4.31 in Guillou and May~[\Enrich],
whose acyclicity condition is satisfied by Remark~4.34 there:
all functors $\tC(b,a)\otimes_R-\colon\Ch(R)\to\Ch(R)$ are left Quillen functors.
By construction, $\tC\dashdgMod$ is a tractable left proper stable model category,
because it inherits these properties from~$\Ch(R)$.

\proclaim Remark.
An astute reader will notice that the above definition involves mapping complexes whose source is not necessarily cofibrant.
One is lead to wonder whether such a notion is homotopically correct.
Proposition~3.2 in To\"en~[\HoDG] shows that Dwyer--Kan equivalences of source categories
induce Dwyer--Kan equivalences of dg modules.
One can also show that the quasicategory $(\Ndg\tC\dashdgMod)[W^{-1}]$ is equivalent to $\Fun(\Ndg\tC,\cU\Ch(R))$
by generalizing the proof of Proposition~1.3.4.25 in [\HA].

\proclaim Definition.
A~dg module on a small dg category~$\tC$ is {\it weakly representable\/} if it is weakly equivalent to a representable dg module.

\proclaim Notation.
Let $\tC\dashdgMod_\wr^\circ$ denote the full subcategory of weakly representable bifibrant dg modules.

Now we turn to defining homotopy limits.

\label\dgholim
\proclaim Definition.
A~dg category~$\tC$ {\it admits $D$-indexed homotopy limits}
if the (objectwise) homotopy limit of a $D$-indexed diagram of representable dg modules
is weakly representable.

This definition says that we view a diagram in $\tC$ as a diagram in $\tC\dashdgMod$,
where we can compute its homotopy limit using model categories.
If that homotopy limit in dg modules is representable up to weak equivalence,
then we know what the homotopy limit should be in $\tC$ itself.

\proclaim Remark.
Other approaches to homotopy limits use weighted colimits
(i.e., replicate a dg analog of the Bousfield--Kan approach to simplicially enriched categories)
or exploit the dg nerve construction, discussed below,
which produces a quasicategory from a dg category.

With these definitions in hand,
we introduce the class of well-behaved dg categories relevant to stable phenomena,
which are the natural analogs of stable quasicategories or stable model categories.
(A~useful reference for these notions and for comparisons with \ii-categories is~[\DGstable].)

\proclaim Definition.
A~dg category is {\it pretriangulated\/} if it admits finite homotopy limits (in the sense of the above definition)
and the suspension of any representable module is weakly representable.

Our definition is equivalent to the usual one by the following observations.
First, there is a terminal object in the sense that the zero module,
which assigns the zero chain complex to every object of $\tC$,
is weakly representable.
In fact, the existence of a terminal object
and of homotopy fibers (computed as objectwise mapping cocones)
implies the existence of equalizers (constructed as the homotopy fiber of the difference)
and therefore arbitrary finite homotopy limits,
constructed as the equalizer of two maps $$\prod_{O\in\rm Ob}O\to\prod_{f\in\rm Mor}{\rm codom}(f).$$

The name ``pretriangulated'' is justified by the fact that the homotopy category
(i.e., $\hodg(\tC)$) of a pretriangulated category~$\tC$ is a triangulated category;
see, for example, Proposition~1 in Bondal and Kapranov~[\Enhanced].
Additionally, the dg nerve of a pretriangulated dg category is a stable \ii-category;
see, for example, Theorem~4.3.1 in Faonte~[\SNIC].
Thus one can see dg categories as a particular strict model for stable \ii-categories
that have generalized Eilenberg--MacLane spectra as spectral mapping complexes.

\proclaim Definition.
A~dg category $\tC$ is {\it idempotent-complete\/} (alias Karoubi complete) if retracts of representable dg modules on $\tC$ are weakly representable.

Examples of idempotent-complete pretriangulated dg categories include all the classical ones:
chain complexes of $R$-modules,
modules of chain complexes on a site,
chain complexes of quasicoherent sheaves in algebraic geometry,
among many others.
(See To\"en~[\LecDG] or simply follow references to Bondal and Kapranov~[\Enhanced] and Keller~[\KelDG] for more examples.)

\proclaim Remark.
Two types of limits will show up in our constructions below:
homotopy fibers (i.e., derived kernels) and homotopy sequential limits.
We recall two concrete models for these notions in the world of chain complexes.
\ppar
The {\it mapping cocone\/} of a chain map $f\colon X\to Y$ in $\Ch(R)$
is given by the complex
$$C(f) = (X \oplus Y[-1], \d_X + \d_{Y[-1]} + f).$$
That is, the differential sends $(x,y)$ to $(\d x, \d y + f(x))$.
A~cocycle $(x,y)$ is precisely a cocycle $x$ for $X$ with the property that $f(x)$ is exact in~$Y$, i.e., $f(x)$ is homologous to zero.
There is a chain map $C(f) \to X$ by projection.
\ppar
Given a functor $F\colon\Zs\to\udg(\Ch(R))$,
the product $\prod_n F(n)$ admits two natural endomorphisms:
the identity $\id$ and the ``shift'' $s_F$,
which is the product of the maps $F(n-1)\to F(n)$.
The {\it mapping cotelescope\/} of a functor $F\colon\Zs\to\udg(\Ch(R))$
is the mapping cocone of the endomorphism $\id-s_F$ of the object~$\prod_n F(n)$.

\subsection Relationship with model and \ii-categories

To work in a homotopically correct fashion with dg categories
(and to compare with other approaches in this paper),
it will be useful to relate them with quasicategories and model categories.

\label\dgnerves
\subsubsection The differential graded nerves

There are several ways to convert a dg category into a quasicategory.
On the one hand, one can piggyback on the Dold--Kan correspondence:
construct a simplicially enriched category by keeping the same objects but applying the Dold--Kan correspondence to (the appropriate truncation of) the hom-complexes, and then take the homotopy-coherent nerve to obtain a quasicategory.
This {\it big dg nerve\/} is defined as the (homotopy coherent) nerve
of the simplicial category of Construction~1.3.1.13 in~[\HA].
On the other hand, Lurie constructs a {\it small dg nerve\/} functor~$\Ndg$ in \S1.3.1 of~[\HA],
which will be our preferred construction and which below we will simply call the dg nerve, following Lurie.
Although the two nerves are not isomorphic,
Proposition~1.3.1.17 establishes a quasicategorical equivalence between them.
(Additional details can be found in [\SNIC].
Faonte introduced the terminology of big/small dg nerve.)

In practice, one must take care when thinking about a dg category as presenting an $\infty$-category,
for exactly the same reasons that one works with projective or injective resolutions in ordinary homological algebra.
For example, it is often the case that the dg category is equipped with a model structure enriched over $\Ch(R)$.
In such setting, where $\bC$ is a dg model category,
Proposition~1.3.1.17 in Lurie~[\HA] shows that the underlying quasicategory of~$\bC$ is equivalent to the dg nerve of the category of {\it bifibrant\/} objects in~$\bC$.

It is important to recognize that the dg nerve functors discard a substantial amount of information:
they only see the {\it connective\/} part of each chain complex of maps in a dg category.
Recall that for an unbounded chain complex~$X$,
its connective cover $\tau_{\geq 0} X \to X$ is obtained by setting
$(\tau_{\geq 0} X)_n=X_n$ for $n>0$, $(\tau_{\geq 0}X)_n=0$ for $n<0$, and $(\tau_{\geq 0} X)_0=\ker(X_0\to X_{-1})$.
(This construction is often called the ``smart truncation''.)
The connective cover functor is symmetric lax monoidal,
so applying it to each hom-complex in a dg category~$\tC$
produces another dg category $\tau_{\geq 0}\tC$.
The dg nerve functors send the canonical functor $\tau_{\geq 0}\tC\to\tC$ to an isomorphism,
and hence lose data about an arbitrary dg category.

For us, however, the dg categories of interest are pretriangulated,
and a pretriangulated dg category can be canonically reconstructed from its connective cover.
This feature is analogous to how stable \ii-categories can be defined as ordinary \ii-categories with certain properties
(so that the homotopy groups of its mapping spaces are concentrated in nonnegative degrees),
even though every stable \ii-category is canonically enriched over the \ii-category of spectra
and its mapping spectra naturally have homotopy groups in both positive and negative degrees.

\subsubsection Comparison results

With these tools in place, we can introduce an analogy that may be helpful
for the reader coming from higher categories or algebraic topology:
dg categories are to stable model categories as chain complexes are to spectra.
Indeed, this analogy can be made precise, as follows.
First, every combinatorial stable model category is Quillen equivalent to
the category of enriched presheaves of spectra on a small category enriched in spectra.
Furthermore, if we restrict to those combinatorial stable model categories whose mapping spectra are Eilenberg--MacLane spectra,
then every such model category is Quillen equivalent to
the model category of dg modules on a small dg category.
Thus, the underlying quasicategory of such a combinatorial stable model category
is equivalent to the dg nerve of the bifibrant objects in the associated dg model category.

This analogy can be pushed further.
The category of dg categories $\dgCat_R$ has an obvious notion of weak equivalence:
a dg functor is a {\it quasi-equivalence\/} if it is always a quasi-isomorphism on hom-complexes and it induces a categorical equivalence of homotopy categories.
But one can enlarge this class of weak equivalences of dg categories to {\it Morita equivalences\/}:
these are the dg functors that induce a Dwyer--Kan equivalence of the associated categories of dg modules.
There is a model category structure on the category of small dg categories,
using Morita equivalences and Dwyer--Kan cofibrations,
and the fibrant objects are precisely small idempotent-complete pretriangulated dg categories.
Furthermore, the underlying \ii-category of this ``Morita'' model category
is equivalent to the \ii-category of small idempotent-complete $R$-linear stable \ii-categories and exact functors (Corollary 5.5 in Cohn~[\DGstable]).
In short, we know how to translate between a nice class of stable \ii-categories and their dg counterparts.

\subsection A useful technical proposition

The following proposition says something quite useful: for a dg category~$\tC$,
the \ii-category obtained by inverting chain homotopy equivalences on the underlying category $\udg\tC$ knows all the same higher data as~$\tC$ itself.
This fact will play a role in proving \vdgmaincf, our main comparison theorem for dg categories.

\label\dgrect
\proclaim Proposition.
Suppose $\tC$ is a small dg category.
Denote by $W$ the set of chain homotopy equivalences in~$\tC$.
The canonical functor $\N\udg(\tC)[W^{-1}]\to\Ndg\tC$ is an equivalence of quasicategories
provided that $\tC$ can be realized as a full dg subcategory of some dg model category $\tC'$
such that the objects in~$\tC$ are cofibrant as objects of~$\tC'$
and are closed under chain homotopy equivalences between cofibrant objects in~$\tC'$.

\proclaim Remark.
Our result is a modest generalization of Proposition~1.3.4.5 of~[\HA],
which applies to a subcategory of a category $\Ch(A)$ of chain complexes.
To prove that proposition, Lurie proves a statement similar to ours for simplicially-enriched categories in Proposition~1.3.4.7 of~[\HA].

\proof Proof.
The proof consists of producing a zigzag of (weak) equivalences between the relevant relative categories, using several intermediate relative categories.
We will also make manifest exactly where the hypothesis on~$\tC$ appears.
\ppar
Consider the model category $\tC\dashdgMod$ (\vdgcmod).
The dg Yoneda embedding identifies $\tC$ with the full dg subcategory $\tC\dashdgMod^\circ_\repr$ of representable dg presheaves
in $\tC\dashdgMod^\circ$, the bifibrant objects of $\tC\dashdgMod$.
Thus, the induced functor $\Ndg\tC\to\Ndg(\tC\dashdgMod^\circ)$ is fully faithful,
so our goal will be to identify the category $\N\udg(\tC)[W^{-1}]$ with this full subcategory.
By \vchainmaps\ below,
the functor $$(\N\udg(\tC\dashdgMod^\circ))[W^{-1}]\to\Ndg(\tC\dashdgMod^\circ)$$ is an equivalence of quasicategories.
This allows us to get rid of the dg nerves and work exclusively with relative categories.
\ppar
To show that $\udg\tC\to\udg(\tC\dashdgMod^\circ)$
is a homotopically fully faithful functor of relative categories,
we recall that such functors satisfy the 2-out-of-3 property.
First, we postcompose with the inclusion $$\udg(\tC\dashdgMod^\circ)\to\udg(\tC\dashdgMod),$$
which is a homotopy equivalence of relative categories; a bifibrant replacement functor determines an inverse.
\ppar
Next, we factor $\udg\tC\to\udg\tC\dashdgMod$ through
the inclusion $\udg(\tC\dashdgMod_\wr)\to\udg\tC\dashdgMod$, which is homotopically fully faithful by \vfullyfaithful\ below.
We then further factor $\udg\tC\to\udg(\tC\dashdgMod_\wr)$
through the inclusion $\udg(\tC\dashdgMod_\wr^\circ)\to\udg(\tC\dashdgMod_\wr)$,
which is a homotopy equivalence of relative categories by virtue of a bifibrant replacement functor.
\ppar
We now deploy some tools from~\vtstructuresmodel\ to show that the inclusion $\udg\tC\to\udg(\tC\dashdgMod_\wr^\circ)$ is a Dwyer--Kan equivalence of relative categories.
Let $I_E$ denote all the cofibrations obtained as a tensor product of a representable dg presheaf with a generating cofibration $\Z[n-1]\to(\Z[n-1]\gets\Z[n])$, for all $n\ge0$, in the category of chain complexes of abelian groups.
The induced homotopy factorization system determines a t-structure on $\tC\dashdgMod$ with a truncation functor $\tau_{\ge0}$.
As usual, this functor fits into an adjunction
$$\iota_{\ge0} \colon\tau_{\ge0}\tC\dashdgMod \longleftrightarrows \tC\dashdgMod \colon \tau_{\ge0},$$
as explained in~\vtstructuresmodel.
\ppar
Consider the inclusion $$\tau_{\ge0}\udg(\tC\dashdgMod_\wr^\circ)\to\udg(\tC\dashdgMod_\wr^\circ),$$
where the left category is the image of $\udg(\tC\dashdgMod_\wr^\circ)$ under $\tau_{\ge0}$
(which differs from dg presheaves with vanishing negative chain degrees because hom-complexes in~$\tC$ need not be concentrated in degree~0).
The bifibrancy condition is preserved by $\tau_{\ge0}$: all objects are fibrant, whereas cofibrancy is preserved by the very definition of~$\tau_{\ge0}$:
we factor $0\to X$ into $0\to\tau_{\ge0}X\to X$, where the first morphism belongs to~$E$, which is a subclass of cofibrations.
The weak representability condition is preserved by $\tau_{\ge0}$ because $\tau_{\ge0}$ preserves weak equivalences and sends representable dg presheaves to themselves.
The unit of the second adjunction is the identity map
and the counit of the second adjunction is the natural transformation $\tau_{\ge0}X\to X$, which is a weak equivalence by assumption.
Thus the second adjunction is a homotopy equivalence of relative categories.
Hence the right inclusion $$\tau_{\ge0}\udg(\tC\dashdgMod_\wr^\circ)\to\udg(\tC\dashdgMod_\wr^\circ)$$ is a Dwyer--Kan equivalence of relative categories.
\ppar
To show that the inclusion $\udg\tC\to\tau_{\ge0}\udg\tC\dashdgMod_\wr^\circ$
is a homotopy equivalence of relative categories,
we construct an inverse relative functor $\tau_{\ge0}\udg\tC\dashdgMod_\wr^\circ\to\tC$.
It is at this point that the special assumptions of the dg category~$\tC$ become important;
the previous discussion applies to any small dg category.
The dg category $\tC\dashdgMod$ is the free dg cocompletion of~$\tC$,
which means that there is a canonical dg cocontinuous (i.e., preserving colimits and tensorings) dg functor $V\colon\tC\dashdgMod\to\tC'$.
The functor $\udg V$ sends generating (acyclic) cofibrations to (acyclic) cofibrations
(since all objects in $\tC\subset\tC'$ are cofibrant and $\tC'$ is a dg model category)
and therefore is a left Quillen functor.
In particular, it preserves weak equivalences between cofibrant objects.
We restrict the domain of~$\udg V$ to $\tau_{\ge0}\udg\tC\dashdgMod_\wr^\circ$ and claim that
the resulting functor preserves weak equivalences and factors through $\tC\subset\tC'$.
Weak equivalences between bifibrant objects are chain homotopy equivalences,
and any dg functor preserves chain homotopy equivalences.
Any object~$X$ in $\tau_{\ge0}\udg\tC\dashdgMod_\wr^\circ$ is chain homotopy equivalent to a representable dg presheaf $X'$ (itself a bifibrant object),
which implies that the images of both $X$ and $X'$ under~$\udg V$ are also cofibrant and chain homotopy equivalent.
Therefore $\udg V(X)$ is in~$\tC$.
\ppar
Thus we have constructed functors $\tC\to\tau_{\ge0}\udg\tC\dashdgMod_\wr^\circ$ and $\tau_{\ge0}\udg\tC\dashdgMod_\wr^\circ\to\tC$
that preserve weak equivalences.
The composition $\tC\to\tC$ is the identity functor.
The other composition is not the identity functor, but becomes the identity if we restrict its domain and codomain to representable dg presheaves.
This shows that $\tC\to\tau_{\ge0}\udg\tC\dashdgMod_\wr^\circ$ is a homotopically fully faithful functor.
This functor is homotopy essentially surjective by definition of weakly representable dg presheaves.
Thus it is a Dwyer--Kan equivalence of relative categories.
\ppar
Our zigzag of fully faithful functors of \ii-categories between $\Ndg\tC$ and $\udg\tC[W^{-1}]$
shows that the functor $\udg\tC[W^{-1}]\to\Ndg\tC$ is fully faithful.
It is also essentially surjective because both \ii-categories have the same set of objects.
Therefore it is an equivalence of \ii-categories.

The proof of the proposition relied on two lemmas that we now prove.
The first lemma uses the Dold--Kan correspondence,
and we denote the Dold--Kan functor by $\Gamma\colon\Ch_{\ge0}\to\sSet$.

\label\chainmaps
\proclaim Lemma.
Suppose $\bC$ is a model category enriched over the model category of unbounded chain complexes of abelian groups with its projective model structure.
Given a cofibrant object~$X$ and a fibrant object~$Y$ in~$\bC$,
the derived mapping space $\Map(X,Y)$ can be computed as $\Gamma\tau_{\ge0}\bC(X,Y)$,
the Dold--Kan functor of the connective cover (\vdgnerves) of the hom-complex~$\bC(X,Y)$.

\proof Proof.
Let $\Mor_\bC$ denote the set of morphisms in the category~$\bC$.
The derived mapping space $\Map(X,Y)$ can be computed as the simplicial set $n\mapsto\Mor_\bC(W_n,Y)$,
where $W_\bullet$ is a Reedy cofibrant replacement of the constant cosimplicial object with value~$X$.
For instance, we can take $W_\bullet=\Z[\ss^\bullet]\otimes X$, with induced cosimplicial maps,
which are all weak equivalences, so $W$ is indeed a replacement.
Furthermore, $W_\bullet$ is Reedy cofibrant because the latching maps $\Z[\sb^n]\otimes X\to\Z[\ss^n]\otimes X$ are cofibrations.
We now compute
$$\eqalign{\Map(X,Y)
&=(n\mapsto\Mor_\bC(W_n,Y))\cr
&=(n\mapsto\Mor_\bC(\Z[\ss^n]\otimes X,Y))\cr
&=(n\mapsto\Mor_{\Ch}(\Z[\ss^n],\bC(X,Y)))\cr
&=(n\mapsto\Mor_{\Ch_{\ge0}}(\Z[\ss^n],\tau_{\ge0}\bC(X,Y)))\cr
&=(n\mapsto\Mor_\sSet(\ss^n,\Gamma\tau_{\ge0}\bC(X,Y)))\cr
&=\Gamma\tau_{\ge0}\bC(X,Y),\cr
}$$
which completes the proof.

\label\fullyfaithful
\proclaim Lemma.
Suppose $\bC$ is a model category and $\bD\subset\bC$ is a full subcategory of~$\bC$ closed under weak equivalences.
Then the inclusion $\bD\to\bC$ is a homotopically fully faithful functor of relative categories,
i.e., the induced functor of quasicategories $\cU\bD\to\cU\bC$ is fully faithful.

\proof Proof.
Proposition~6.2(i) of~[\Calc] gives a simple construction of the mapping space
in a relative category that admits a homotopy calculus of fractions (Definition~6.1(i) in~[\Calc]).
Specifically, the simplicial set $\Map(X,Y)$ in such a relative category can be computed as the nerve
of the category $\bW^{-1}\bC\bW^{-1}$ as defined in~\S5.1 of~[\Calc].
(Take ${\bf m}=\bW^{-1}\bC\bW^{-1}$ there.)
Its objects are zigzags $X\mathbin{\tilde\gets}A\to B\mathbin{\tilde\gets}Y$
and morphisms are pair of weak equivalences $A\mathbin{\tilde\to}A'$, $B\mathbin{\tilde\to}B'$ that make the corresponding diagram commute.
\ppar
We verify below that $\bD$ and $\bC$ admit a homotopy calculus of fractions.
This allows us to conclude that the inclusion $\bD\to\bC$ induces
a fully faithful functor of the associated simplicial categories,
i.e., a homotopically fully faithful functor of the associated hammock localizations.
\ppar
The existence of a homotopy calculus of fractions is granted by Proposition~8.2 in [\Calc],
whose conditions are satisfied for~$\bC$ by Proposition~8.4 there.
The same holds for $\bD$, for which it suffices to observe that all constructions in Proposition~8.4
produce new objects in~$\bC$ that are weakly equivalent
to one of the objects of the given diagram in~$\bD$ and hence are in~$\bD$ themselves.

\subsection The filtered category of a dg category

We now turn to our main purpose: describing the filtered dg category of a well-behaved dg category, namely an idempotent-complete pretriangulated dg category.
In effect, we mimic our constructions with quasicategories and model categories, with minor modifications.
Unfortunately, a technical issue prevents mindless mimicry:
although the category of dg categories $\dgCat_R$ has a closed monoidal structure given by the pointwise tensor product,
this monoidal product is {\it not\/} compatible with Dwyer--Kan equivalences.
In particular, the correct (derived) tensor product or hom of dg categories
cannot be computed using this monoidal structure.
Hence we cannot simply define the dg category of sequences $\dgSeq(\tC)$ of $\tC$
as the hom of dg categories $\Hom_{\dgCat_R}(\Zs_R, \tC)$
where $\Zs_R$ is the dg category whose objects are integers and whose morphisms $m\to n$ ($m\le n$) form the chain complex $R[0]$.
Instead, we take a bare-handed approach.

\proclaim Remark.
Another construction of the derived internal hom of dg categories is given by To\"en [\HoDG]:
if $A$ and $B$ are two dg categories,
their derived internal hom can be computed as the dg category of bifibrant right quasi-representable dg modules on $A^\op\hat\otimes B$,
where $\hat\otimes$ denotes the derived tensor product.
Here a dg module $F\colon A\hat\otimes B^\op\to\Ch(R)$ is {\it right quasi-representable\/} if for any $a\in A$ the functor $F(a,-)\colon B^\op\to\Ch(R)$
is weakly equivalent to a representable dg module on~$B$.
The model structure can be taken to be the projective model structure, so bifibrancy reduces to projective cofibrancy.
We choose not to use this construction because it is not very explicit.

\proclaim Definition.
A~{\it sequence\/} in a dg category~$\tC$ is a functor $F\colon \Zs\to\udg(\tC)$.

The obvious, non-dg notion of a map of sequences $F \to G$ is a collection of maps $F(n) \to G(n)$ that intertwine with the structure maps $F(n) \to F(n+1)$ of the sequences.
Using this definition, we have an ordinary category of sequences,
which should be the underlying category of any reasonable dg category of sequences.
That is, to have a dg category of sequences,
we need a hom-complex between any pair of sequences,
with the property that its zero cycles recover this obvious notion.
To do that, we need to unpack the obvious notion just discussed.

Given a sequence $F$ in $\tC$, let $s_F$ denote the {\it shift},
which is the map $$s_F\colon \prod_{n \in \Z} F(n) \to \prod_{n \in \Z} F(n)$$
that sends an element $(a_n)_{n \in \Z}$ to the family whose $n$th element is $F(n-1)(a_{n-1})$.
That is, it applies the defining structure maps of~$F$.
Thus the obvious notion of map of sequences $\phi \colon F\to G$ is that $s_G \circ \phi = \phi \circ s_F$.
We can weaken this condition in a natural way to get the hom-complex.

\proclaim Definition.
Given sequences $F$, $G$ in~$\tC$,
let $\dgSeq(F,G)$ denote the mapping cocone of the map
$$[-,s]_{FG} \colon \prod_{n\in\Z}\tC(F(n),G(n))\to\prod_{m\in\Z}\tC(F(m),G(m+1))$$
sending a family $(\phi_n)$ of maps to the family $(\phi_n\circ s_F-s_G \circ\phi_{n-1})$.
In other words, it is the homotopy equalizer for precomposition by the shift for~$F$ and postcomposition by the shift for~$G$.
(Here $[-,s]$ is meant to suggest a commutator with~$s$.)

\proclaim Remark.
If one replaces mapping cocone with the kernel
(equivalently, the homotopy equalizer with an ordinary equalizer)
in the definition above,
one recovers the notion of a (strict) natural transformation between two functors $\Zs\to\tC$.
In our case, natural transformation commute up to a given homotopy.

We want $\dgSeq(F,G)$ to be the morphisms in a dg category,
so we need to describe composition as a chain map and verify the associativity of composition.
Explicitly, we define
$$\circ_{\dgSeq}\colon\dgSeq(G,H) \otimes \dgSeq(F,G) \to \dgSeq(F,H)$$
as follows.
An element of $\dgSeq(F,G)$ is a pair
$$(\phi,\psi) \in \prod_{n\in\Z}\tC(F(n),G(n))\oplus\prod_{m\in\Z}\tC(F(m),G(m+1))[-1],$$
and the differential acts as
$$\d(\phi,\psi) = (\d\phi, -\d\psi + [\phi,s]),$$
where on the left side $\d$ means the differential in $\dgSeq(F,G)$,
but on the right side $\d\phi$ means the differential in $\prod_{n\in\Z}\tC(F(n),G(n))$ and $\d\psi$ means the differential in $\prod_{n\in\Z}\tC(F(n),G(n+1))$.
We verify that $\d^2=0$:
$$\d^2(\phi,\psi)=\d(\d\phi, -\d\psi + [\phi,s])=(\d^2\phi,\d^2 \psi - \d [\phi,s] + [\d\phi,s])=(0,0).$$
Then we define the composition by
$$(\phi',\psi') \circ_{\dgSeq} (\phi,\psi) = (\phi' \circ \phi, \psi' \circ \phi + \phi' \circ \psi).$$
Direct computation shows that this composition is a chain map and also associative.

\label\dgseq
\proclaim Definition/Lemma.
There is a {\it dg category of sequences\/} $\dgSeq(\tC)$ in a dg category~$\tC$
whose objects are functors $F\colon\Zs\to\udg(\tC)$
and whose hom-complexes are $\dgSeq(F,G)$,
with composition as just defined.

\proclaim Remark.
Our construction $\dgSeq$ actually defines an endofunctor on the category $\dgCat_R$ of dg categories.
Indeed, it is a functor that preserves Dwyer--Kan equivalences and hence is a functor of relative categories.
To see that $\dgSeq$ induces a quasi-isomorphism on each mapping complex,
note that the construction of $\dgSeq$ uses homotopy limits of mapping complexes and therefore preserves quasi-isomorphisms.
To see that $\dgSeq$ preserves quasi-equivalences, we argue as follows.
The cofibrant replacement functor $\crf$ on dg categories
preserves the set of objects and produces a cofibrant dg category.
If we apply it to a quasi-equivalence $\tC\to\tD$,
then the $\crf\tC\to\crf\tD$ is not only a quasi-equivalence but
is in fact a chain homotopy equivalence (a functor that admits an inverse functor, up to a chain homotopy),
because any quasi-equivalence between cofibrant dg categories is one.
Moreover, any chain homotopy equivalence of dg categories is a Dwyer--Kan equivalence.
By functoriality, $\dgSeq$ preserves chain homotopy equivalences.
Hence, we have a commutative square using the functors $\tC\to\tD$, $\crf\tC\to\crf\tD$, $\crf\tC\to\tC$, and $\crf\tD\to\tD$;
these last two functors are quasi-equivalences that induce identities on objects.
If we apply $\dgSeq$ to this square, we see that all four functors induce quasi-isomorphisms on hom-complexes and
the functor $\dgSeq\crf\tC\to\dgSeq\crf\tD$ is a chain homotopy equivalence.
Thus, it suffices to show that the functors $\dgSeq\crf\tC\to\dgSeq\tC$ and $\dgSeq\crf\tD\to\dgSeq\tD$
are homotopy essentially surjective.
Concretely, an object~$X$ of $\dgSeq\tE$ (where $\tE$ is any dg category)
is a sequence of objects~$X(n)$ of~$\tE$ together with elements $f(n)\in\udg\tE(X(n),X(n+1))$.
The dg functor $\crf\tE\to\tE$ is an acyclic fibration of dg categories,
in particular, every map $\crf\tE(X(n),X(n+1))\to\tE(X(n),X(n+1))$ is an acyclic fibration of chain complexes,
so every element $f(n)$ can be lifted to $\crf\tE(X(n),X(n+1))$.
By the 2-out-of-3 property for Dwyer--Kan equivalences of dg categories,
we thus deduce that $\dgSeq(\tC)\to\dgSeq(\tD)$ is a Dwyer--Kan equivalence.

We now turn to complete filtered sequences.
We choose to identify these as a dg {\it sub\/}category of sequences.

\proclaim Definition.
A~sequence~$F$ in a dg category~$\tC$ is {\it complete\/} if its homotopy limit is weakly equivalent to the zero module on~$\tC$.

This definition says, in essence, that a sequence $X$ is complete if the ``initial'' term $X(-\infty) = \lim X(n)$ is~0,
where limit means the homotopy limit over $\Zs$.
Motivation for this definition can be found in \vicomp,
where we explain how the classical notion of completeness relates to this version.
If $\tC$ has a zero object, then a bounded below sequence (i.e., $F$ for which $F(n)=0$ for sufficiently small~$n$) is complete.

Thus we come to our key definition.

\label\dgdef
\proclaim Definition.
The {\it filtered dg category\/} $\dgFil(\tC)$ is the full dg subcategory of $\dgSeq(\tC)$ consisting of complete sequences.

To justify this definition, we now turn to showing that it agrees with our stable \ii-category constructions
when we work with well-behaved dg categories.

A~technical proposition plays a key role.
Recall that in a model category $\bC$, we denote the full subcategory of bifibrant objects by~$\bC^\circ$.

\label\dgmodelbig
\proclaim Proposition.
Let $\bA$ be a tractable, left proper, stable $\Ch(R)$-enriched model category.
Then the dg category $\bSeq(\bA)^\circ$ is quasi-equivalent to the dg category $\dgSeq(\bA^\circ)$ constructed in \vdgseq.
Likewise, the dg category $\bFil(\bA)^\circ$ is quasi-equivalent to the dg category $\dgFil(\bA^\circ)$ constructed in \vdgdef.

\proof Proof.
We have a dg functor $\bSeq(\bA)^\circ\to\dgSeq(\bA^\circ)$, which is an inclusion on objects.
(Recall that a bifibrant functor $\Zs\to\bA$ has bifibrant levels.)
For the mapping complexes, we use the following map in~$\Ch(R)$:
$$\bSeq(F,G)=\int_{n\in\Z}\bA(F(n),G(n))\to\dgSeq(F,G),$$
where the chain homotopy given by the second summand in $\dgSeq(F,G)$ is set to zero,
which also ensures that the composition is preserved.
Recall that $\dgSeq$ was defined as a homotopy end of the same type.
Due to the fact that $F$ and $G$ are projectively bifibrant,
the end also computes the homotopy end, so the map is a weak equivalence.
Homotopy essential surjectivity is shown by applying the projective cofibrant and fibrant replacement functors.
\ppar
By \vlocalwe, the dg category $\bFil(\bA)^\circ$ is the full subcategory of $\bSeq(\bA)^\circ$ consisting of derived complete sequences.
Homotopy limits of sequences computed in the model category $\bSeq(\bA)^\circ$ map to homotopy limits of sequences in the dg category $\dgSeq(\bA^\circ)$ (\vdgholim),
thus $\bFil(\bA)^\circ)$ maps to the full subcategory $\dgFil(\bA^\circ)$ of $\dgSeq(\bA^\circ)$ consisting of derived complete sequences.
Thus the dg functor $\bFil(\bA)^\circ\to\dgFil(\bA^\circ)$ is a quasi-equivalence.

The dg model category $\tC\dashdgMod$ is a tractable, left proper, stable Ch(R)-enriched model category, as noted in \vdgcmod, so we can apply \vdgmodelbig.
This leads to another useful description of the relevant categories, in the following corollary.

Recall that $\tC\dashdgMod_\wr^\circ$ denotes the full subcategory of weakly representable bifibrant dg modules.
Similarly, let $\bSeq(\tC\dashdgMod_\wr)^\circ$ denote the full subcategory of $\bSeq(\tC\dashdgMod)$
consisting of bifibrant sequences whose entries are weakly representable dg modules.
Let $\bFil(\tC\dashdgMod_\wr)^\circ$ denote the corresponding full subcategory of~$\bFil(\tC\dashdgMod)$.

\label\dgmodel
\proclaim Corollary.
For any small idempotent-complete pretriangulated dg category~$\tC$,
the dg category $\dgSeq(\tC)$ is quasi-equivalent to $\bSeq(\tC\dashdgMod_\wr)^\circ$.
Similarly, $\dgFil(\tC)$ is quasi-equivalent to $\bFil(\tC\dashdgMod_\wr)^\circ$.

\proclaim Remark.
Below we require all these categories to be {\it small\/} in order to apply \vdgrect.
This size condition is achieved via a standard trick:
we work only with those dg modules $M$ such that for each~$X \in \tC$,
the underlying set of the chain complex $M(X)$ is a subset of a set of sufficiently large cardinality.
(Sufficiently large here means that the inclusion of such restricted objects
into all objects is homotopy essentially surjective, i.e.,
any object is weakly equivalent to a restricted object.)
We use this convention till the end of this subsection.

\proof Proof.
The Yoneda embedding $\tC\to\tC\dashdgMod_\wr^\circ$
is a quasi-equivalence of dg categories,
where the target consists of bifibrant weakly representable dg modules.
Therefore, the induced dg functor $$\dgSeq(\tC)\to\dgSeq(\tC\dashdgMod_\wr^\circ)$$ is also a quasi-equivalence.
\ppar
The dg functor constructed in \vdgmodelbig\ restricts to a dg functor $$\bSeq(\tC\dashdgMod_\wr)^\circ\to\dgSeq(\tC\dashdgMod_\wr^\circ),$$
which induces quasi-isomorphisms on hom-complexes---as shown there---and
is homotopy essentially surjective by applying the (projective) bifibrant replacement functor.
For the filtered case we further restrict the functor to $\bFil(\tC\dashdgMod_\wr)^\circ\to\dgFil(\tC\dashdgMod_\wr^\circ)$.

Finally, we obtain our main comparison result for dg categories.

\label\dgmaincf
\proclaim Theorem.
For any small idempotent-complete pretriangulated dg category~$\tC$,
the dg nerve of $\dgFil(\tC)$ is equivalent to the filtered \ii-category of the dg nerve of~$\tC$.
Similarly the dg nerve of $\dgSeq(\tC)$ is equivalent to the \ii-category of sequences in the dg nerve of~$\tC$.

\proof Proof.
We establish the statement for sequences first.
We have
$$\eqalign{\Ndg\dgSeq(\tC)
&\simeq\Ndg\bSeq(\tC\dashdgMod_\wr)^\circ\cr
&\simeq\cU(\bSeq(\tC\dashdgMod_\wr)^\circ)\cr
&\simeq\cU(\bSeq(\tC\dashdgMod_\wr))\cr
&\simeq\iSeq(\cU(\tC\dashdgMod)_\wr)\cr
&\simeq\iSeq(\cU(\tC\dashdgMod_\wr))\cr
&\simeq\iSeq(\Ndg(\tC)).\cr
}$$
Here $\dgMod_\wr$ denotes the full dg subcategory of dg modules weakly equivalent to representable dg modules.
The first step was established in \vdgmodel.
The second step follows from \vdgrect:
take $\bSeq(\tC\dashdgMod_\wr)^\circ$ as $\tC$
and $\bSeq(\tC\dashdgMod)$ as $\tC'$.
Then we have an equivalence of quasicategories $$\cU(\bSeq(\tC\dashdgMod_\wr)^\circ)\simeq\Ndg(\bSeq(\tC\dashdgMod_\wr)^\circ).$$
\ppar
The third step is induced by the inclusion $\bSeq(\tC\dashdgMod_\wr)^\circ\to\bSeq(\tC\dashdgMod_\wr)$,
which has an inverse given by the bifibrant replacement.
\ppar
The third to the last step follows from \vbFilvsiFil\ combined with the fact that the equivalence
between $\cU(\bSeq(\tC\dashdgMod))$ and $\iSeq(\cU(\tC\dashdgMod))$
respects the full subcategories of sequences of representables.
\ppar
The second to the last step is a simple consequence of the fact that $\tC\dashdgMod_\wr\to\tC\dashdgMod$ is a full dg subcategory,
hence so is $\cU(\tC\dashdgMod_\wr)\to\cU(\tC\dashdgMod)$.
\ppar
For the last step we must establish an equivalence of quasicategories
$\cU(\tC\dashdgMod_\wr^\circ)\to\Ndg(\tC)$.
Observe that $\tC\to\tC\dashdgMod_\wr^\circ$ is a quasi-equivalence of dg categories.
Now \vdgrect\ shows the desired equivalence:
take $\tC\dashdgMod_\wr^\circ$ as~$\tC$
and $\tC\dashdgMod$ as~$\tC'$.
Then we have an equivalence of quasicategories $\cU(\tC\dashdgMod_\wr^\circ)\simeq\Ndg(\tC\dashdgMod_\wr^\circ)\simeq\Ndg(\tC)$.
\ppar
We extend the statement from sequences to filtered objects
by observing that the full subcategories of complete sequences
in $\dgSeq(\tC)$ and $\iSeq(\Ndg(\tC))$ correspond to each other under the above chain.

\subsection Uniqueness of the dg enhancement of the filtered derived category

Let $\abcat$ be a Grothendieck abelian category.
Our work above shows that its filtered derived category $\Dfil(\abcat)$ (see \vfildercat)
admits a dg enhancement,
i.e., a pretriangulated dg category whose homotopy category is equivalent to $\Dfil(\abcat)$ as a triangulated category.
We now turn to showing that this enhancement is unique in the sense of Orlov and Lunts (see [\Uniq, \EnhGroth]),
meaning that any other pretriangulated dg category whose homotopy category is equivalent to $\Dfil(\abcat)$
is quasi-equivalent (via a zigzag of quasi-equivalences) to the dg category $\dgFil(\Ch(\abcat)^\circ)$.

One convenient way to ensure that one makes homotopically meaningful constructions with $\Ch_\dg(\abcat)$
is to equip it with a combinatorial dg model structure whose weak equivalences are quasi-isomorphisms and fibrations are degreewise epimorphisms.
Then we work with the full dg subcategory $\Ch_\dg(\abcat)^\circ$
of bifibrant objects.

\proclaim Proposition.
The homotopy categories of the dg categories $\dgSeq(\Ch(\abcat)^\circ)$ and $\dgFil(\Ch(\abcat)^\circ))$
admit unique enhancements (as triangulated categories) to pretriangulated dg categories.

\proof Proof.
Our approach will exploit the connection with model categories.
As noted, $\Ch(\abcat)$ admits a natural combinatorial model structure.
Further, the model category $\bSeq(\Ch(\abcat))$ is the same as the model category $\Ch(\Seq(\abcat))$,
where $\Seq(\abcat):=\Fun(\Zs,\abcat)$ is again a Grothendieck abelian category,
so $\Ch(\Seq(\abcat))$ provides a model category, as just described.
\ppar
By \vdgmodelbig, we know that the dg category $\dgSeq(\Ch(\abcat)^\circ)$
has the same underlying homotopy category as the model category $\bSeq(\Ch(\abcat))$.
But the homotopy category of $\bSeq(\Ch(\abcat))$ is precisely the derived category of the Grothendieck abelian category $\Seq(\abcat)$:
chain complexes of sequences in~$\abcat$ are the same as sequences of chain complexes in~$\abcat$,
and quasi-isomorphisms of chain complexes of sequences in~$\abcat$ are the same as sequences of quasi-isomorphisms of chain complexes in~$\abcat$.
Theorem~A in Canonaco and Stellari [\EnhGroth]
then tells us that the derived category of any Grothendieck abelian category
(including $\Seq(\abcat)$) admits a unique dg enhancement.
\ppar
The case of $\dgFil(\Ch(\abcat)^\circ))$ is more subtle,
as we cannot directly apply Theorem~A\null.
Instead, note that $\dgFil(\Ch(\abcat)^\circ)$
has the same underlying homotopy category as $\bFil(\Ch(\abcat))$
by \vdgmodelbig.
We then apply Theorem~C of [\EnhGroth] to $\bFil(\Ch(\abcat))$, as follows.
In the notation of that theorem,
let the category $\bA$ be the full subcategory of $\Seq(\abcat)$ on some set of generators of $\Seq(\abcat)$,
and let $\cD(\bA)$ be the derived category of $\bA$ viewed as a dg category.
Let the localizing subcategory~$\bL_0$ of $\cD(\bA)$
be such that the (Verdier) localization of $\cD(\bA)$ at~$\bL_0$ is equivalent to $\cD(\Seq(\abcat))$ as a triangulated category
(and the left Bousfield localization of $\Fun(\bA^\op,\Ch(\Z))$ at~$\bL_0$ is Quillen equivalent to $\bSeq(\Ch(\abcat))$ as a model category).
Now construct a new localizing subcategory~$\bL$ by adding to $\bL_0$ all the homotopy constant sequences and then taking the closure under shifts and small homotopy colimits.
Hence, the new localization is the homotopy category of $\bFil(\Ch(\abcat))$
and is well generated because $\bFil(\Ch(\abcat))$ is a combinatorial model category.
\ppar
Condition~(b) in the cited work requires that in the homotopy category of $\bFil(\Ch(\abcat))$,
there is only one morphism from $A[0]$ to $\coprod_I A_i[k_i]$, where $A$ and $A_i$ are objects in~$\bA$, $I$ is a set, and $k_i<0$.
This condition is satisfied because the objects involved are complete and for degree reasons every map must be the zero map
(for instance, the projective resolution will replace $A[0]$ with something concentrated in degrees~0 and below,
so clearly there are no maps to $A_i[k_i]$ if $k_i<0$, since this lives in a positive degree).

\label\filtoperads
\section Application: filtered operads and algebras over them

In this section we very briefly indicate how to obtain a homotopy theory of filtered operads, and (filtered) algebras,
by combining the abstract machinery of this paper with that in Pavlov and Scholbach~[\Operads].

Given a monoidal model category $\bC$ and a set $W$ of colors,
let $\Oper_W(\bC)$ denote the category of $W$-colored symmetric operads with values in $\bC$.
(See Kelly~[\OpMay] for a definition in the single-colored case.)
Given an operad $P \in \Oper(\bC)$, and $\bM$ a left module category over $\bC$,
let $\Alg_P(\bM)$ denote the category of $P$-algebras in $\bM$.
(Often we take $\bM=\bC$.)

\proclaim Definition.
Let $\bA$ be a symmetric monoidal stable model category.
The category of $W$-colored (symmetric) {\it filtered operads\/} in~$\bA$
is the category of $W$-colored (symmetric) operads in the symmetric monoidal category $\bFil(\bA)$ of filtered objects in~$\bA$.

\proclaim Definition.
Let $\bM$ be a stable model category that is a left module over $\bA$.
The category of {\it filtered algebras\/} in $\bM$ over a filtered operad $O$ is $\Alg_O(\bFil(\bM))$.

\proclaim Remark.
We emphasize that a filtered operad is not the same as a sequence of operads,
and we will not develop here the notion of a filtered object in operads.
In a sequence $\{O_k\}$ of operads, one can only compose operations at a fixed index $k$,
whereas in a filtered operad, composition mixes stages of the filtration.
A~similar issue is visible at the level of associative algebras.
In a filtered associative algebra, the product of an element from the $m$th stage and $n$th stage lives in the $(m+n)$th stage.
In a sequence of associative algebras, by contrast, supposing $m \le n$,
an element in the $m$th algebra maps to an element of the $n$th algebra,
but any products then remain in the $n$th algebra.

\proclaim Remark.
When $\bA$ is the category of chain complexes,
the notion of a filtered operad considered here recovers that defined by Kimura, Stasheff, and Voronov
in~\S5 of~[\HMC] by requiring the $n$-ary operations of an operad to be a cofibrant sequence.

\proclaim Theorem.
Let $\bA$ be a (symmetric) flat, (symmetric) h-monoidal, symmetric monoidal, compactly generated, combinatorial stable model category.
For a fixed set of colors~$W$ the category of filtered $W$-colored (symmetric) operads in~$\bA$ has a compactly generated combinatorial model structure,
as does the category of filtered algebras over any fixed filtered colored (symmetric) operad.
Furthermore, weak equivalences of operads induce Quillen equivalences of model categories of algebras.

\proof Proof.
Apply Theorem~5.10 and Theorem~7.10 of~[\Operads] to the symmetric monoidal model category $\bFil(\bA)$,
which satisfies the assumptions of these theorems by \vfiltprops.

\proclaim Example.
If $\bA$ is the category of rational chain complexes (or over any characteristic zero field $k$),
then the conditions of the above theorem are satisfied in the symmetric case.
Hence filtered algebras over filtered operads in chain complexes over~$k$ admit a model structure,
and quasi-isomorphisms of operads induce Quillen equivalences.

\proclaim Nonexample.
The conditions of the above theorem are satisfied only in the nonsymmetric case
when $\bA$ is the category of chain complexes of abelian groups, vector spaces over a field of characteristic~$p$,
or, more generally, modules over a ring.
For such categories, the symmetric case does not apply.
Recall, for instance, that the transferred model structure does not exist on commutative differential graded algebras in characteristic~$p$.

For symmetric operads one needs the rather strong condition of symmetric flatness, which is not satisfied, for instance, by simplicial sets or chain complexes of abelian groups.
However, for the case when $\bA$ is the category of symmetric spectra, one can obtain much easier criteria to verify, as explained in the work of the second author and Scholbach~[\Spectra].

\proclaim Theorem.
Let $\bC$ be a flat, h-monoidal, symmetric monoidal, compactly generated, combinatorial model category,
$R$ be a commutative monoid in symmetric sequences in~$\bC$,
and $W$ be a fixed set of colors.
The category of filtered $W$-colored (symmetric) operads in symmetric $R$-spectra in~$\bC$
admits a compactly generated combinatorial model structure,
as does the category of filtered algebras over a fixed such operad.
Furthermore, weak equivalences of filtered operads induce Quillen equivalences of categories of filtered algebras.

\proof Proof.
This follows from Theorem~3.4.1 and Theorem~3.4.4 in~[\Spectra]
and \vfiltprops.

Often, we are interested in the case where $R$ is a free commutative monoid on a symmetric sequence concentrated in degree~1,
where it is given by an object $R_1\in\bC$.

\proclaim Example.
If $\bC$ is pointed simplicial sets and $R_1=S^1_*$, the pointed simplicial circle, then symmetric $R$-spectra are symmetric simplicial spectra.
We thus obtain a model structure on filtered algebras in simplicial symmetric spectra over filtered operads in simplicial symmetric spectra, in particular, filtered simplicial operads.

\proclaim Example.
For natural sequences of Lie groups, such as the orthogonal groups $\O(n)$, the Thom spectrum $\M\O$ is filtered by the sequence $\M\O(n)$.
Similarly, $\M\U$ is filtered by the sequence $\M\U(n)$.
These are important examples of filtered $\E_\infty$ ring spectra.
The (symmetric) Thom spectrum $\M\O(m)$ has as its $k$th spectral level the simplicial set $\B\O(m)^{\R^k\oplus V_m}$,
where $V_m\to\B\O(m)$ is the universal rank $m$ vector bundle over~$\B\O(m)$
and the superscript denotes the Thom space of a vector bundle.
We have inclusion maps $\M\O(m)\to\M\O(m+1)$
and multiplication maps $\M\O(m)\wedge\M\O(n)\to\M\O(m+n)$.
The latter are induced on individual spectral levels $k$~and~$l$ by the maps $\B\O(m)^{\R^k\oplus V_m}\wedge\B\O(n)^{\R^l\oplus V_n}\to\B\O(m+n)^{\R^{k+l}\oplus V_{m+n}}$
constructed using the universal property of~$\B\O(m+n)$.
Combined together, this data gives us a filtered $\E_\infty$-ring spectrum~$\M\O$.

\proclaim Example.
The Adams--Novikov spectral sequence of an $\E_\infty$-algebra spectrum~$E$ over an $\E_\infty$-ring spectrum~$S$
is obtained from the filtered object given by applying the quasicategorical Dold--Kan correspondence to the bar construction of~$E$ over~$S$,
a cosimplicial spectrum (the $n$th term is $E\wedge_S\cdots\wedge_SE$, where $E$ repeats $n+1$ times).

\proclaim Example.
Consider the category of pointed motivic spaces, i.e., simplicial presheaves on the Nisnevich site localized at Nisnevich covers and the maps ${\bf A}^1\times X\to X$.
If $R_1={\bf P}^1_*$, then symmetric $R$-spectra are motivic symmetric spectra
and we obtain a model structure on filtered algebras over filtered motivic operads.

\proclaim Remark.
As an example of a situation involving more complicated ambient categories than chain complexes or spectra,
one can cite the {\it slice filtration\/} in motivic homotopy theory, which gives rise to the slice spectral sequence.

It can be convenient to consider the associated graded operad of a filtered operad (similarly, the associated graded algebra of a filtered algebra).
For the case of chain complexes, these notions were considered by Dotsenko~[\DQPB] and Griffin~[\OpCom], for example.

\proclaim Proposition.
For any symmetric monoidal stable model category $\bC$ such that the categories of filtered and graded (symmetric) operads admit a transferred model structure,
the associated graded functor $\Gr$ induces a left Quillen functor from filtered operads to graded operads.
For a fixed filtered operad~$O$ such that filtered algebras over~$O$ and graded algebras over~$\Gr(O)$ admit a transferred model structure,
the associated graded functor induces a left Quillen functor from filtered algebras over~$O$ to graded algebras over~$\Gr(O)$.

\proof Proof.
As established in \vmodelgrmonoidal, the associated graded functor $\Gr$ is strong monoidal.
In particular, it induces functors from filtered operads to graded operads
and from algebras over a filtered operad to algebras over the associated graded operad.
We treat both cases simultaneously; recall that operads are themselves algebras over a certain colored operad.
Recall that $\Gr$, on the level of $\bC$, has a right adjoint $R$, which inserts zeros as transition maps.
For trivial reasons, this right adjoint is also strong monoidal.
Then $R$ induces a functor on the level of algebras (or operads), and it is a right Quillen functor because both model structures are transferred.
The functor $\Gr$, on the level of algebras, is left adjoint to this functor and is thus a left Quillen functor.

\section Application: filtered $\cD$-modules

In this section we apply the abstract theory above to construct a model structure on filtered $\cD$-modules.
First, though, we need to introduce some machinery for sheaves. Throughout, $X$ will denote a complex manifold.

\proclaim Remark.
We could also take $X$ to be a smooth variety (with the Zariski topology) over~$\C$,
though in such a situation one typically also wants log poles, which would make the exposition
here considerably more involved.
See Pavlov and Scholbach~[\MdR] for an account of this situation.

Denote by $\PSh(X,\Ch(\C))$ the category of presheaves on~$X$ of chain complexes of complex vector spaces equipped with the opens-wise monoidal product, which we denote $\otimes_\C$.
There is always the projective model structure on this category,
where a map of presheaves is a weak equivalence if it is a quasi-isomorphism on every open and
a map of presheaves is a fibration if it is a fibration on every open.
This model structure does not care about the topology on $X$, however, for which we need the following definitions.
A~map $f \colon F \to G$ of presheaves is a {\it local weak equivalence\/} if every map on stalks $f_x \colon F_x \to G_x$ is a quasi-isomorphism.
In fact, every presheaf is actually locally weakly equivalent to a sheaf.

\proclaim Proposition.
The left Bousfield localization at the local weak equivalences of $(\PSh(X,\Ch(\C)),\otimes_\C)$ with the projective model structure
is a symmetric flat, symmetric h-monoidal, symmetric monoidal, left proper, compactly generated, combinatorial stable model category $\bMod_\C(X)$.

\proof Proof.
The relevant machinery on Bousfield localizations in this context can be found around Definition~6.1.1 in~[\HTSP].
The projective model structure is itself transferred along the forgetful functor given by the evaluation on all objects,
and the projective model structure on unbounded chain complexes has all the properties in the statement, as explained in \S7.4 of~[\HTSP].
As established in \S5 and \S6 of~[\HTSP], all these properties survive the constructions mentioned above.
\proclaim Remark.
Because we work with presheaves localized along local weak equivalences, many constructions become simpler.
For instance, we will never apply sheafification in defining some functor:
as an example, we work directly with the tensor product of presheaves, which is formed open by open, and do not sheafify.

We now turn our attention to $\cD$-modules.
Let $\cO$ denote the structure sheaf of $X$: for instance, if $X$ is a complex manifold, then $\cO$ denotes the sheaf of holomorphic functions.
The sheaf $\End_\C(\cO)$ of $\C$-linear endomorphisms admits a filtration, defined as follows:
for $n<0$, $\D_n=0$ and for $n\ge0$,
$$\D_n=\{P\in\End_\C(\cO)\mid[P,\cO]\in\D_{n-1}\}.$$
Then $\cD = \colim_n \D_n$ is the sheaf known as the {\it sheaf of differential operators\/} on $X$.
The resulting filtered object admits a monoid structure given by composition of endomorphisms,
which amounts to saying that the composition map $\D_m\otimes\D_n\to\D_{m+n}$ satisfies the associativity and unitality conditions.
Thus $\cD$ is a monoid in $(\bFil(\bMod_\C(X)),\otimes_\C)$.
(We have to use complex vector spaces instead of $\cO$-modules because $\cD$ has two different actions of~$\cO$.)

Let $\fMod_\cD(X))$ denote the category of filtered left $\cD$-modules over~$X$, i.e., the category of left modules over the monoid~$\cD$ in filtered presheaves of chain complexes.

\proclaim Proposition.
The category $\fMod_\cD(X)$ admits a left proper, compactly generated, combinatorial stable model structure transferred from filtered presheaves of chain complexes of complex vector spaces.

\proof Proof.
Apply \vmodules.

This stable model category is a useful enhancement of the classical approaches to working with filtered $\cD$-modules,
such as that developed by Laumon in~[\DmodF].
Historically, people have worked with filtrations that are eventually bounded below.
This category is a full subcategory of the category defined here.
Moreover, filtered weak equivalences among such filtered $\cD$-modules are weak equivalences in our category.
Our stable model category then provides functorial homotopy limits and colimits for arbitrary diagrams,
which enlarges the class of manipulations beyond those available in the derived category itself.

\proclaim Remark.
In~[\MdR], filtered $\cD$-modules are equipped with a more pliable model structure, known as the {\it flat\/} model structure,
which enlarges the cofibrations so as to produce a {\it monoidal\/} model category.
The subtle Koszul duality between $\cD$-modules and $\Omega$-modules can then be concretely realized as
a Quillen equivalence of combinatorial stable monoidal model categories between {\it filtered\/} $\cD$-modules and {\it filtered\/} $\Omega$-modules.
In~[\MdR], this duality is used to develop the functoriality of filtered $\cD$-modules in an efficient and clean way.
In particular, given a map $f \colon X \to Y$ of varieties, a Quillen adjunction is constructed
whose right derived functor recovers Laumon's derived direct image functor~[\DmodF],
the key ingredient in his proof of a Riemann--Roch theorem for $\cD$-modules.
The approach chosen in~[\MdR] essentially constructs the pullback-pushforward adjunction for filtered $\cD$-modules
using the pullback-pushforward Quillen adjunction for $\Omega$-modules (defined in a straightforward way),
pre- and postcomposed with the Quillen equivalences between $\cD$-modules and $\Omega$-modules.
This highly structured approach does not impose any conditions (such as smoothness) on~$f$
and allows one to define pullback and pushforwards of arbitrary $\cD$-modules, without additional restrictions
such as holonomicity or boundedness of the filtration.
Moreover, this formalism provides a powerful, efficient toolkit for mixed Hodge modules, a key aim of~[\MdR], by allowing log poles.

\section References

\refs

\endinput